%% file: main.tex
\documentclass[a4paper,onecolumn, superscriptaddress,10pt,shorttitle=papers]{compositionalityarticle}
\pdfoutput=1
\usepackage[utf8]{inputenc}
\usepackage[english]{babel}
\usepackage[T1]{fontenc}
\title{Data-Parallel Algorithms for String Diagrams}
\author{Paul Wilson}
\affiliation{Independent}
\author{Fabio Zanasi}
\affiliation{University College London and University of Bologna}
\input{preamble}

\begin{document}
\maketitle

\begin{abstract}
  \label{section:abstract}
  \input{abstract}

\end{abstract}

\section{Introduction}
\label{section:introduction}
\input{introduction}

\section{Background}
\label{section:background}
\input{background}

\section{Representing Diagram Wirings with $\ACSet$s}
  \label{section:diagram-wirings}
  \input{representing-diagram-wirings}

  \subsection{$\Wires$: The Category of Multi-Sorted Finite Functions}
  \label{section:adding-colour}
  \input{adding-colour}

  \subsection{Bipartite Multigraphs}
  \label{section:bipartite-multigraph}
  \input{bipartite-multigraphs}

\section{String Diagrams as Structured Cospans}
\label{section:string-diagrams}
\input{string-diagrams}

\section{Diagrams as Combinatorial Syntax} %
\label{section:isomorphism}
\input{isomorphism}

\section{Symmetric Monoidal Case}
\label{section:symmetric-monoidal}
\input{symmetric-monoidal}

\section{A Faster Functor}
\label{section:fast-functor}
\input{fast-functor}

\section{Fast Wiring}
\label{section:fast-wiring}
\input{fast-wiring}

\section{Applying Functors to Diagrams}
\label{section:applying-functors}
\input{applying-functors}

\section{Optic Composition using Frobenius Structure}
\label{section:optic-composition}
\input{optic-composition}

\section{Discussion and Future Work}
\label{section:conclusion}
\input{conclusion}

\bibliographystyle{plainnat}
\bibliography{main}

\appendix

\section{Sorting Permutations}
\label{section:sorting}
\input{sorting}

\section{Proofs for \Cref{section:isomorphism}}
\label{section:proofs-isomorphism}
\input{proofs-isomorphism}

\section{Frobenius Spiders as Cospans}
\label{section:cospan-proofs}
\input{cospan-proofs}

\section{Frobenius Spider Lemmas}
\label{section:spider-lemmas}
\input{spider-lemmas}

\end{document}

%% file: preamble.tex
\usepackage[numbers]{natbib}
\usepackage{amsmath, amssymb, amsthm}
\usepackage{mathrsfs} %
\usepackage{hyperref}
\usepackage{float}
\usepackage{xcolor}
\usepackage{enumitem}
\usepackage{stmaryrd}
\usepackage{cleveref}
\usepackage{mdframed}
\usepackage{booktabs}
\usepackage[frozencache=true]{minted}

\usepackage{tikzit}
\input{figures.tikzstyles}

\newcommand{\gen}[1]{\scalebox{0.8}{\input{#1.tikz}}}
\newcommand{\sgen}[1]{\scalebox{0.6}{\input{#1.tikz}}}

\newcommand{\derivecong}[1]{\ensuremath{\overset{(#1)}{\cong}}}

\newtheorem{observation}{Remark}[section]

\newtheorem{definition}[observation]{Definition}
\newtheorem{example}[observation]{Example}
\newtheorem{remark}[observation]{Remark}

\newtheorem{proposition}[observation]{Proposition}
\newtheorem{corollary}[observation]{Corollary}

\newcommand{\deftext}[1]{\textbf{#1}}
\newcommand{\defeq}[0]{\ensuremath{:=}}
\newcommand{\cp}[0]{\ensuremath{\fatsemi}} %
\newcommand{\id}[0]{\ensuremath{\mathsf{id}}}
\newcommand{\twist}[0]{\ensuremath{\sigma}}
\newcommand{\bigtensor}[0]{\ensuremath{\bigotimes}}
\newcommand{\tensor}[0]{\ensuremath{\otimes}}
\newcommand{\unit}[0]{\ensuremath{I}}
\newcommand{\counit}[0]{\ensuremath{0}}
\newcommand{\coeq}[0]{\ensuremath{\mathsf{c}}}

\newcommand{\cat}[1]{\ensuremath{\mathscr{#1}}}
\newcommand{\Cat}[1]{\ensuremath{\mathsf{#1}}}
\newcommand{\functor}[1]{\ensuremath{\mathsf{#1}}}
\newcommand{\Functor}[1]{\ensuremath{\mathsf{#1}}}

\newcommand{\source}[0]{\mathsf{source}}
\newcommand{\target}[0]{\mathsf{target}}
\newcommand{\sources}[0]{\mathsf{sources}}
\newcommand{\targets}[0]{\mathsf{targets}}
\newcommand{\values}[0]{\mathsf{values}}

\newcommand{\cod}[0]{\mathsf{cod}}
\newcommand{\dom}[0]{\mathsf{dom}}

\newcommand{\inj}[0]{\iota}
\newcommand{\proj}[0]{\pi}

\newcommand{\sort}[0]{\ensuremath{\mathsf{sort}}}

\newcommand{\terminal}[0]{\ensuremath{\mathsf{!}}}
\newcommand{\terminalObj}[0]{\ensuremath{1}}
\newcommand{\initial}[0]{\ensuremath{\mathsf{?}}}

\newcommand{\Nat}[0]{\ensuremath{\mathbb{N}}}

\newcommand{\ord}[1]{\ensuremath{\overline{#1}}}
\newcommand{\typer}[0]{\ensuremath{\tau}}

\newcommand{\Set}[0]{\Cat{Set}}
\newcommand{\FinFun}[0]{\Cat{FinFun}}
\newcommand{\FinOrd}[0]{\Cat{FinOrd}}

\newcommand{\Free}[0]{\Cat{Free}}
\newcommand{\FreeFrob}[0]{\Free_{\Signature + \Frob}}
\newcommand{\FreeSMC}[0]{\Free_{\Signature}}

\newcommand{\Diagram}[0]{\Cat{Diagram}}
\newcommand{\DiagramFrob}[0]{\Diagram_{\Signature + \Frob}}
\newcommand{\DiagramSMC}[0]{\Diagram_{\Signature}}
\newcommand{\ToDiagram}[0]{\Functor{ToDiagram}}
\newcommand{\ToDiagramSlow}[0]{\mathcal{S}}
\newcommand{\ToDiagramFast}[0]{\mathcal{F}}
\newcommand{\FromDiagram}[0]{\mathcal{D}}

\newcommand{\OpticObj}[2]{\ensuremath{\left({#1 \atop #2}\right)}}
\newcommand{\ToOptic}[0]{\ensuremath{\mathcal{O}}}
\newcommand{\fwd}[1]{\ensuremath{\overrightarrow{#1}}}
\newcommand{\rev}[1]{\ensuremath{\overleftarrow{#1}}}

\newcommand{\singleton}[0]{\ensuremath{\mathtt{singleton}}}

\newcommand{\Signature}[0]{\ensuremath{\Sigma}}
\newcommand{\Obj}[0]{\ensuremath{\Signature_0}}
\newcommand{\Arr}[0]{\ensuremath{\Signature_1}}
\newcommand{\Eq}[0]{\ensuremath{\Signature_2}}

\newcommand{\CMon}[0]{\mathbf{CMon}}
\newcommand{\Frob}[0]{\mathbf{Frob}}

\newcommand{\bl}[0]{\langle}
\newcommand{\br}[0]{\rangle}

\newcommand{\tl}[0]{\langle}
\newcommand{\tr}[0]{\rangle}

\newcommand{\ACSet}{\ensuremath{\Cat{ACSet}}}

\newcommand{\indegree}[0]{\ensuremath{\mathsf{indegree}}}
\newcommand{\outdegree}[0]{\ensuremath{\mathsf{outdegree}}}

\newcommand{\SchemaBPMG}[0]{\ensuremath{\mathcal{S}}}
\newcommand{\BPMG}[0]{\ensuremath{\mathsf{BipartiteMultigraph}}}

\newcommand{\SchemaWires}[0]{\ensuremath{\mathcal{W}}}
\newcommand{\Wires}[0]{\Cat{Wires}}

\newcommand{\WiresFrob}[0]{\Wires_{\Signature}}

\newcommand{\MkWiresObj}{\ensuremath{\SchemaWires}}
\newcommand{\MkWiresArr}{\ensuremath{\SchemaWires}}
\newcommand{\MkD}[1]{\ensuremath{\mathsf{D}(#1)}}

\newcommand{\len}[1]{\ensuremath{|#1|}}

\newcommand{\Ei}[0]{\ensuremath{E_i}}
\newcommand{\Eo}[0]{\ensuremath{E_o}}
\newcommand{\Ki}[0]{\ensuremath{K_i}}
\newcommand{\Ko}[0]{\ensuremath{K_o}}

\newcommand{\swi}[0]{\ensuremath{w_i}}
\newcommand{\swo}[0]{\ensuremath{w_o}}
\newcommand{\sxi}[0]{\ensuremath{x_i}}
\newcommand{\sxo}[0]{\ensuremath{x_o}}

\newcommand{\sporti}[0]{\ensuremath{\mathsf{port}_i}}
\newcommand{\sporto}[0]{\ensuremath{\mathsf{port}_o}}
\newcommand{\swn}[0]{\ensuremath{w_n}}
\newcommand{\sxn}[0]{\ensuremath{x_n}}

\newcommand{\wNode}[0]{\blacksquare}
\newcommand{\xNode}[0]{\circ}

\newcommand{\Ob}[0]{\ensuremath{\mathsf{Ob}}}
\newcommand{\Sig}[0]{\ensuremath{\mathsf{Sig}}}

\newcommand{\Port}[0]{\ensuremath{\mathsf{Port}}}

\newcommand{\SCL}[0]{\ensuremath{\functor{L}}}
\newcommand{\SCR}[0]{\ensuremath{\functor{R}}}
\newcommand{\SCRL}[1]{\ensuremath{\MkD{#1}}}
\newcommand{\HalfSpider}[0]{\ensuremath{\functor{S}}}

\newcommand{\iW}[0]{\ensuremath{\epsilon}}

\newcommand{\parent}[0]{\ensuremath{\mathtt{parent}}}
\newcommand{\depth}[0]{\ensuremath{\mathtt{depth}}}
\newcommand{\isRoot}[0]{\ensuremath{\mathtt{root}}}
\newcommand{\isRight}[0]{\ensuremath{\mathtt{isRight}}}
\newcommand{\isLeft}[0]{\ensuremath{\mathtt{isLeft}}}
\newcommand{\isOdd}[0]{\ensuremath{\mathtt{odd}}}
\newcommand{\isEven}[0]{\ensuremath{\mathtt{even}}}
\newcommand{\ancestor}[0]{\ensuremath{\mathtt{ancestor}}}

\newcommand{\R}[0]{\ensuremath{\mathsf{R}}}
\newcommand{\smcinterleave}[0]{\ensuremath{\mathsf{interleave}}}

\newcommand{\Csp}[0]{\ensuremath{\mathsf{Csp}}}
\newcommand{\cospan}[5]{\ensuremath{#1 \xrightarrow{#2} #3 \xleftarrow{#4} #5}}

\newcommand{\codelink}[0]{\href{https://yarrow.id}{https://yarrow.id}}

%% file: figures.tikzstyles
\tikzstyle{edge}=[fill=white, draw=black, shape=circle, inner sep=2.5pt]
\tikzstyle{node}=[fill=black, draw=black, shape=circle, inner sep=2.5pt]
\tikzstyle{morphism}=[fill=white, draw=black, shape=rectangle]
\tikzstyle{negate}=[fill=red, draw=red, shape=rectangle]
\tikzstyle{blacksquare}=[fill=black, draw=black, shape=rectangle, scale=0.5]

\tikzstyle{pointy}=[->]
\tikzstyle{dashy}=[-, dashed]
\tikzstyle{dashpoint}=[dashed, ->]
\tikzstyle{bluefill}=[-, fill={rgb,255: red,190; green,240; blue,255}]
\tikzstyle{greyfill}=[-, fill={rgb,255: red,240; green,240; blue,240}]
\tikzstyle{rededge}=[-, draw=red, fill=none]

%% file: abstract.tex
We give parallel algorithms for string diagrams represented as
\emph{structured cospans of ACSets}.
Specifically, we give \emph{linear} (sequential) and \emph{logarithmic}
(parallel) time algorithms for composition, tensor product, construction of
diagrams from arbitrary $\Signature$-terms, and application of functors to
diagrams.

Our datastructure can represent morphisms of both the free symmetric monoidal
category over an arbitrary signature as well as those with a chosen Special
Frobenius structure.
We show how this additional (hypergraph) structure can be used to map diagrams
to \emph{diagrams of optics}.
This leads to a case study in which we define an algorithm for efficiently
computing symbolic representations of gradient-based learners based on reverse
derivatives.

The work we present here is intended to be useful as a \emph{general purpose}
datastructure.
Implementation requires only \emph{integer arrays} and well-known algorithms,
and is data-parallel by constuction.
We therefore expect it to be applicable to a wide variety of settings, including
embedded and parallel hardware and low-level languages.

%% file: introduction.tex
String diagrams are a formal graphical syntax~\cite{selinger} for representing
morphisms of monoidal categories which is now widely used
(see for example
\cite{genovese2021categorical,genovese2021nets,compositional-digital-circuits,compositional-signal-flow-theory}).
The purpose of this paper is to make string diagrams not just a convenient
notation for algebraic reasoning, but also an efficient general-purpose tool in
computing with graphical structures in a compositional manner.
To that end, the datastructures and algorithms we define satisfy the following
desiderata.
\begin{description}
  \item[Fast and data-parallel.]
    Our algorithms are data-parallel by construction, and have \emph{linear}
    (sequential) and \emph{logarithmic} (parallel) time complexities.
  \item[Minimal primitives.]
    Our datastructures are defined in terms of simple integer arrays.
    Moreover, we assume only a small number of simple, well-known primitive
    operations (e.g., prefix sum).
    This makes it possible to implement our algorithms in a wide variety of
    settings, such as embedded and parallel (i.e., GPU) hardware.
  \item[Simple to implement correctly.]
    Key parts of our datastructure are defined in terms of the recent
    construction of \emph{ACSets}~\cite{acsets}.
    Consequently, implementations are essentially the same as their categorical
    definitions, making it easier to ensure correctness.
\end{description}

A number of representations of string diagrams have been explored in the
literature, such as the wiring diagrams of Catlab.jl~\cite{catlab} and the
`hypergraph adjacency representations' of~\cite{coc}.
Our goals most closely align with the latter:
we aim to make string diagrams useful as a general purpose `scalable
combinatorial syntax'.
For example, we hope that our implementation serves as an alternative in cases
where a programmer would currently use a tree or directed graph.

However, the primary motivating application for our work is in representing
gradient-based learners as optics, as described in \cite{cfgbl}.
In particular, this motivates perhaps the most significant extension
to~\cite{coc}: our datastructures can `natively' represent morphisms of the free
symmetric monoidal category over a signature with a chosen Special Frobenius
monoid.
This equips categories with hypergraph structure, which we show can be used to
simulate \emph{diagrams of optics}.
In turn, this allows for a large number of applications modelling `bidirectional
information flow' based on optics
such as~\cite{cfgbl, bayesian-optics, extensive-form-games}.
In our specific example, it allows us to define to an efficient algorithm for
taking reverse derivatives~\cite{rdc} and modeling gradient-based learners in
general.

Our main contributions are as follows:
\begin{itemize}
  \item A representation of morphisms of the free symmetric monoidal category
    over a signature $\Signature$ as \emph{structured cospans of ACSets}.
  \item Proof of correspondence between this representation and
    free symmetric monoidal category
  \item Data-parallel algorithms with \emph{linear} (sequential) and
    \emph{logarithmic} (parallel) time complexity for...
    \begin{itemize}
      \item Composition and tensor product of diagrams
      \item Construction of a diagram from an arbitrary $\Signature$-term
      \item Application of functors to diagrams
    \end{itemize}
  \item An algorithm for mapping diagrams to \emph{diagrams of optics} using
    hypergraph structure, and consequently an algorithm for taking reverse
    derivatives of diagrams in linear (sequential) and logarithmic (parallel)
    time.
\end{itemize}

The structure of the paper is as follows.
In \Cref{section:background} we give necessary background, including
string diagrams, presentations by generators and equations,
structured cospans, and ACSets.
We also recall the bipartite graph representation of hypergraphs introduced in
\cite{coc}, and give a detailed account of the representation of finite
functions as arrays which is the foundation of our implementation.
Our contributions begin in \Cref{section:diagram-wirings} where we show how the
`internal wiring' of diagrams can be represented using ACSets.
This is built upon in \Cref{section:string-diagrams}, where we give the main
definition of the paper: our datastructure for string diagrams as structured
cospans of these wirings.
In \Cref{section:isomorphism} we prove the correspondence between these
structured cospans and the free symmetric monoidal category on a given signature
plus a chosen Special Frobenius monoid.
In \Cref{section:symmetric-monoidal}, we translate a combinatorial condition
first introduced in \cite{rmsms} to our datastructure, allowing for the
representation of morphisms of the free symmetric monoidal category
\emph{without} the additional Frobenius structure.
Sections \ref{section:fast-functor} and \ref{section:fast-wiring} together
define an efficient algorithm for constructing diagrams from $\Signature$-terms.
Finally, in \Cref{section:applying-functors} we define an algorithm for applying
functors to diagrams,
and then in \Cref{section:optic-composition} show how it can be used to map
diagrams to diagrams of optics, leading to an efficient algorithm for taking
reverse derivatives.
We conclude the paper in \Cref{section:conclusion} with some directions for
future work.
A reference implementation for some of the algorithms described in the paper can
be found at \codelink{}.

%% file: background.tex
We introduce the necessary background to describe our contributions.
\Cref{section:background-string-diagrams} recalls string diagrams,
monoidal signatures, and the free symmetric monoidal category presented by a
signature.
In \Cref{section:finite-functions}, we give the details of two isomorphic
categories of finite functions, and give the computational complexities of basic
algorithms for composition, coproduct, and tensor product.
We also list the basic primitive array operations assumed by our implementation.
Finally, in \Cref{section:hypergraphs} we discuss the combinatorial encoding of
string diagrams as hypergraphs introduced in \cite{rmsms},
the encoding of hypergraphs as bipartite graphs used in \cite{coc},
and the definition of $\ACSet$s, in terms of which we will define our
datastructure.

\subsection{String Diagrams, Monoidal Signatures, and $\FreeSMC$}
\label{section:background-string-diagrams}
\input{background/string-diagrams}
\input{background/special-frobenius-monoids}

\subsection{Finite Functions and their Representations}
\label{section:finite-functions}
\input{background/finite-functions}

\subsection{Cospans and Structured Cospans}
\label{section:cospans}
\input{background/cospans}

\subsection{Hypergraphs and Bipartite Graphs for String Diagram Representation}
\label{section:hypergraphs}
\input{background/hypergraphs}

%% file: background/string-diagrams.tex
String diagrams are a two-dimensional graphical syntax for representing
morphisms of categories.
Informally, this syntax consists of widgets placed on the page,
and connected with labeled wires.
The example below has wire labels in capital letters, and is constructed from the widgets
$f, g, h, \sgen{frob/g-join}, \sgen{frob/g-split}$, and $\sgen{frob/g-counit}$.
\begin{equation}
  \label{equation:example-string-diagram}
  \input{examples/string-diagram.tikz}
\end{equation}
The choice of `widgets' and `wire labels' correspond to the \emph{generating
morphisms} and \emph{generating objects} of a \emph{monoidal signature}.

\begin{definition}[Monoidal Signature]
  \label{definition:signature}
  A monoidal signature $\Signature$
  consists of
  \begin{itemize}
    \item $\Obj$ the \deftext{generating objects}
    \item $\Arr$ the \deftext{generating morphisms} or \deftext{operations}
    \item $\Eq$ the \deftext{equations}
    \item A \deftext{typing relation} $\typer : \Arr \to (\Obj^* \times \Obj^*)$
  \end{itemize}
\end{definition}

\begin{remark}
  Note that we include a typing \emph{relation} in order to allow for
  `polymorphic' generators.
  This means that a presentation of a symmetric monoidal category having, for
  example, a generating morphism $f : A \to A$ for all $A$ can be represented
  with just a single label in $\Arr$.
  We have already seen an example of where this is useful in
  \eqref{equation:example-string-diagram}: the generator $\sgen{frob/g-split}$
  is pictured as having both the type $A \to A \tensor A$ and $C \to C \tensor
  C$.
  Although this distinction is useful, we will frequently assume there is a
  chosen typing for a given generator, and therefore speak of `the' type of a
  generator.
\end{remark}

Given a presentation--the basic widgets that define a category--we can now
define $\Signature$-terms: the `syntax trees' representing diagrams built
inductively by tensor and composition of generators.

\begin{definition}
  Given a monoidal signature $\Signature$,
  a \deftext{$\Signature$-term} is a binary tree
  whose leaves are labeled $\Arr \cup \{ \id, \twist \}$,
  and whose nodes are labeled either $\tensor$ or $\cp$.
\end{definition}

Intuitively, a $\Signature$-term is built from building blocks the $\Arr$-operations and the `structural' morphisms $\id$ (identity) and $\twist$ (symmetry), composed sequentially ($\cp$) and in parallel ($\tensor$). We revisit our example string diagram from
\eqref{equation:example-string-diagram}
as a $\Signature$-term in the following example.

\begin{example}
  Equation \eqref{equation:example-string-diagram} can be represented as the $\Signature$-term
  $
    (
      (\sgen{frob/g-split} \otimes (\id \otimes \sgen{frob/g-split} ))
      \cp
      ((f \otimes g) \otimes (h \otimes \id))
    )
    \cp
    (\sgen{frob/g-counit} \tensor (\id \tensor \sgen{frob/g-join}))
  $,
  which we render as a binary tree below.
  \[
    \scalebox{0.7}{\begin{tikzpicture}[level distance=7mm, baseline=-4mm, level 1/.style={sibling distance=7cm},  level 2/.style={sibling distance=5cm}, level 3/.style={sibling distance=3cm}, level 4/.style={sibling distance=1cm}]
      \node {$\cp$}
        child {node {$\cp$}
          child {node {$\tensor$}
            child {node {$\sgen{frob/g-split}$}}
            child {node {$\tensor$}
              child {node {$\id$}}
              child {node {$\sgen{frob/g-split}$}}
            }
          }
          child {node {$\tensor$}
            child {node {$\tensor$}
              child {node {$f$}}
              child {node {$g$}}
            }
            child {node {$\tensor$}
              child {node {$h$}}
              child {node {$\id$}}
            }
          }
        }
        child {
          child {node {$\tensor$}
            child {node {$\sgen{frob/g-counit}$}}
            child {node {$\tensor$}
              child {node {$\id$}}
              child {node {$\sgen{frob/g-join}$}}
            }
          }
        };
    \end{tikzpicture}}
  \]
\end{example}

There are many distinct $\Signature$-terms that represent the diagram in
\eqref{equation:example-string-diagram}.
To represent morphisms of the free symmetric monoidal category, we will
therefore need to consider isomorphism classes of $\Signature$-terms.

\begin{figure}
	\caption{Laws of strict symmetric monoidal categories.}
	\label{fig:smc}
    \begingroup
    \renewcommand{\arraystretch}{3} %
    \begin{tabular}{llc}
      Axiom        & $\Signature$-term                                  & Diagrammatic Syntax \\
      \toprule
      $\tensor$ Functoriality
        & $(f_1 \tensor g_1) \cp (f_2 \tensor g_2) = (f_1 \cp f_2) \tensor (g_1 \cp g_2)$
        & $\input{bg/smc-law-interchange.tikz} = \input{bg/smc-law-interchange.tikz}$ \\
      $\tensor$ Functoriality
        & $\id_A \otimes \id_B = \id_{A \tensor B}$
        & $\input{bg/smc-law-bifunctor-identity-lhs.tikz} \quad = \quad \input{bg/smc-law-bifunctor-identity-rhs.tikz}$ \\
      $\alpha$ Naturality
        & $f \tensor (g \tensor h) = (f \tensor g) \tensor h$
        & $\input{bg/smc-law-naturality-assoc.tikz} \quad = \quad \input{bg/smc-law-naturality-assoc.tikz}$ \\
      $\rho$ Naturality
        & $f \tensor \id_{\unit} = f$
        & ${ \input{bg/g-empty.tikz} \atop \input{bg/g-f.tikz} } \quad = \quad \input{bg/g-f.tikz}$ \\
      $\lambda$ Naturality
        & $\id_{\unit} \tensor f = f$
        & ${ \input{bg/g-f.tikz} \atop \input{bg/g-empty.tikz} } \quad = \quad \input{bg/g-f.tikz}$ \\
      \midrule
      $\twist$ Naturality
        & $(f \tensor g) \cp \twist_{A', B'} = \twist_{A, B} \cp (f \tensor g)$
        & $\scalebox{0.8}{\input{bg/smc-law-twist-naturality-lhs.tikz}} = \scalebox{0.8}{\input{bg/smc-law-twist-naturality-rhs.tikz}}$ \\
      Unit Coherence
        & $\twist_{A, \unit} = \id_A$
        & $\input{bg/smc-law-unit-coherence-lhs.tikz} \quad = \quad \input{bg/g-obj-id.tikz}$ \\
      Symmetry Coherence
        & $(\twist_{A, B} \tensor \id_C) \cp \id_B \tensor \twist_{A, C} = \twist_{A, B \tensor C}$
        & $\input{bg/smc-law-associator-coherence-lhs.tikz} \quad = \quad \input{bg/smc-law-associator-coherence-rhs.tikz}$ \\
      Inverse Law
        & $\twist_{A, B} \cp \twist_{B, A} = \id_A \tensor \id_B$
        & $\input{bg/smc-law-inverse-law-lhs.tikz} \quad = \quad \input{bg/smc-law-inverse-law-rhs.tikz}$ \\
      \bottomrule
    \end{tabular}
    \endgroup
\end{figure}

\begin{definition}[$\FreeSMC$]
  \label{definition:free-smc}
  Given a monoidal signature $\Signature$,
  the \deftext{strict symmetric monoidal category freely generated by $\Signature$}
  is denoted $\FreeSMC$.
  Its objects are lists of generating objects $\Obj^*$, with the unit object the
  empty list.
  Morphisms are $\Signature$-terms quotiented by the laws of symmetric monoidal
  categories (Figure~\ref{fig:smc}).
  Identities, symmetry, and composition and monoidal product are given by their
  corresponding $\Signature$-terms.
\end{definition}

%% file: background/special-frobenius-monoids.tex
\subsubsection{Special Frobenius Monoids and Hypergraph Categories}
\label{section:special-frobenius-monoids}

One monoidal signature particularly important to us is that of Special Frobenius Monoids.

\begin{definition}[Special Frobenius Monoids]
  \label{definition:frobenius}
  The theory of Special Frobenius monoids is denoted $\Frob$, and consists of generators
  \begin{equation}
    \label{equation:frobenius-generators}
    \input{frob/g-split.tikz}
    \qquad
    \qquad
    \input{frob/g-unit.tikz}
    \qquad
    \qquad
    \input{frob/g-join.tikz}
    \qquad
    \qquad
    \input{frob/g-counit.tikz}
  \end{equation}
  and equations
  \begin{equation}
    \label{equation:special-frobenius-axioms}
    \begin{aligned}
      \input{frob/s-split-twist.tikz} & \quad = \quad \input{frob/g-split.tikz}
      \qquad \qquad
      \input{frob/s-split-assoc-l.tikz} & \quad = \quad \input{frob/s-split-assoc-r.tikz}
      \qquad \qquad
      \input{frob/s-split-unit.tikz} & \quad = \quad \input{g-identity.tikz} \\ \\
      \input{frob/s-twist-join.tikz} & \quad = \quad \input{frob/g-join.tikz}
      \qquad \qquad
      \input{frob/s-join-assoc-l.tikz} & \quad = \quad \input{frob/s-join-assoc-r.tikz}
      \qquad \qquad
      \input{frob/s-counit-join.tikz} & \quad = \quad \input{g-identity.tikz} \\ \\
      & \input{frob/s-split-join.tikz} \quad = \quad \input{g-identity.tikz}
      \qquad & \qquad
      \input{frob/s-snake-l.tikz} \quad = \quad \input{frob/s-snake-r.tikz} &
    \end{aligned}
  \end{equation}
\end{definition}

A category in which each object is equipped with a
special frobenius monoid compatibly with the monoidal product is often called a hypergraph category~\cite{fong-thesis} (also known as a well-supported compact closed category~\cite{carboniwalters}).

\begin{definition}[Hypergraph Category]
  A \deftext{Hypergraph Category} is a symmetric monoidal category in which
  every object $A$ is equipped with a Special Frobenius Monoid
  \begin{equation*}
    \label{equation:typed-frobenius-generators}
    \input{frob/gt-A-split.tikz}
    \qquad
    \qquad
    \input{frob/gt-A-unit.tikz}
    \qquad
    \qquad
    \input{frob/gt-A-join.tikz}
    \qquad
    \qquad
    \input{frob/gt-A-counit.tikz}
  \end{equation*}
  satisfying the equations
  \eqref{equation:special-frobenius-axioms}
  and compatible with the tensor product, i.e. so that
  \begin{align*}
    \label{equation:hypergraph-category-coherence}
    \gen{frob/s-coherence-split-lhs} \quad = \quad \gen{frob/s-coherence-split-rhs}
    \qquad \qquad
    &
    \qquad \qquad
    \gen{frob/s-coherence-join-lhs} \quad = \quad \gen{frob/s-coherence-join-rhs} \\ \\
    \gen{frob/s-coherence-unit-lhs} \quad = \quad \gen{frob/s-coherence-unit-rhs}
    \qquad \qquad
    &
    \qquad \qquad
    \gen{frob/s-coherence-counit-lhs} \quad = \quad \gen{frob/s-coherence-counit-rhs}
  \end{align*}
\end{definition}

The hypergraph structure will play a special role: it will be used to represent
the \emph{wires} of the string diagram.
Specifically, it will represent the \emph{wires} of the string diagram.
It will therefore be useful to refer to morphisms constructed exclusively from
generators in $\Frob$.
Such morphisms are called \emph{Frobenius spiders}: we define them now.

\begin{definition}[Frobenius Spider]
  \label{definition:frobenius-spider}
  A \deftext{Frobenius Spider} in a hypergraph category
  is any morphism built by tensor and composition of generators
  in $\Frob$, i.e.,
  \[
    \gen{g-identity}
    \quad
    \gen{g-twist}
    \quad
    \gen{frob/g-split}
    \quad
    \gen{frob/g-join}
    \quad
    \gen{frob/g-unit}
    \quad
    \gen{frob/g-counit}
  \]
\end{definition}

%% file: frob/g-unit.tikz
\begin{tikzpicture}
	\begin{pgfonlayer}{nodelayer}
		\node [style=none] (0) at (0, 0) {};
		\node [style=none] (1) at (0.5, 0) {};
		\node [style=none] (2) at (0, 0.25) {};
		\node [style=none] (3) at (0, -0.25) {};
	\end{pgfonlayer}
	\begin{pgfonlayer}{edgelayer}
		\draw (1.center) to (0.center);
		\draw (2.center) to (3.center);
	\end{pgfonlayer}
\end{tikzpicture}

%% file: frob/g-counit.tikz
\begin{tikzpicture}
	\begin{pgfonlayer}{nodelayer}
		\node [style=none] (0) at (0, 0) {};
		\node [style=none] (1) at (-0.5, 0) {};
		\node [style=none] (2) at (0, 0.25) {};
		\node [style=none] (3) at (0, -0.25) {};
	\end{pgfonlayer}
	\begin{pgfonlayer}{edgelayer}
		\draw (1.center) to (0.center);
		\draw (2.center) to (3.center);
	\end{pgfonlayer}
\end{tikzpicture}

%% file: g-identity.tikz
\begin{tikzpicture}
	\begin{pgfonlayer}{nodelayer}
		\node [style=none] (1) at (0.75, 0) {};
		\node [style=none] (2) at (-0.75, 0) {};
	\end{pgfonlayer}
	\begin{pgfonlayer}{edgelayer}
		\draw (2.center) to (1.center);
	\end{pgfonlayer}
\end{tikzpicture}

%% file: background/finite-functions.tex
Because our work is defined in terms of ACSets~\cite{acsets}, the datastructures
presented in this paper are ultimately all expressed in terms of the category
$\FinOrd$ of \emph{finite sets and functions}.
We will need two different (but isomorphic) `encodings' of this category:
\begin{itemize}
  \item $\Free_{\CMon}$: a `mathematician-friendly' encoding, freely generated by a signature and equations.
  \item $\FinFun$: a `programmer-friendly' encoding, defined in terms of
    \emph{arrays}.
\end{itemize}

Even though both of these representations are well-studied (see e.g.,
\cite{acsets}), it is useful to record them for future use.

\subsubsection{Finite Functions in Terms of Commutative Monoids}
\label{section:finord}

\begin{definition}[Commutative Monoids]
  Given a chosen set of generating objects $\Obj$,
  the theory of commutative monoids $\CMon(\Obj)$
  has generating objects $\Obj$,
  and for each $A \in \Obj$ the generating arrows
  \begin{equation}
    \label{equation:commutative-monoid}
    \gen{cmon/gt-A-add}
    \qquad
    \qquad
    \gen{cmon/gt-A-zero}
  \end{equation}
  and equations
  \begin{equation}
    \label{equation:commutative-monoid}
    \begin{aligned}
      \gen{cmon/st-A-comm} & \: = \: \gen{cmon/gt-A-add}
      \qquad
      \gen{cmon/st-A-unit} & \: = \: \gen{gt-A-identity}
      \qquad
      \gen{cmon/st-A-assoc-lhs} & \: = \: \gen{cmon/st-A-assoc-rhs}
    \end{aligned}
  \end{equation}
\end{definition}

The equations $\CMon_2$ of the theory of commutative monoids are sufficient to deduce the \emph{naturality} equations (\Cref{proposition:finord-add-naturality} and \ref{proposition:finord-zero-naturality}). It is important that these equations are \emph{derivable} (and not axioms) in order to relate finite functions to Special Frobenius Monoids in \Cref{section:special-frobenius-monoids}.

\begin{proposition}
  \label{proposition:finord-zero-naturality}
  For all morphisms $f : A \to B$ in $\CMon(\Sigma_0)$,
  \begin{equation}
    \label{equation:finord-zero-naturality}
    \gen{cmon/st-A-B-f-zero}
    \qquad = \qquad
    \gen{cmon/gt-B-zero}
  \end{equation}
\end{proposition}
\begin{proof}
  Induction.
  It is clear that \eqref{equation:finord-zero-naturality} holds for all generators.
  In the inductive step, assume
  \eqref{equation:finord-zero-naturality}
  holds for morphisms $f_0$ and $f_1$,
  then it is straightforward to derive that
  \eqref{equation:finord-zero-naturality}
  holds for $f_0 \tensor f_1$ and $f_0 \cp f_1$.
\end{proof}

\begin{proposition}
  \label{proposition:finord-add-naturality}
  For all morphisms $f : A \to B$ in $\CMon(\Sigma_0)$,
  \begin{equation}
    \gen{cmon/st-A-B-add-naturality-lhs}
    \qquad = \qquad \gen{cmon/st-A-B-add-naturality-rhs}
  \end{equation}
\end{proposition}
\begin{proof} Induction as in \Cref{proposition:finord-zero-naturality}.\end{proof}

\begin{corollary} By Fox's theorem~\cite{fox},  Propositions  \ref{proposition:finord-zero-naturality} and  \ref{proposition:finord-add-naturality}  mean that $\CMon(\Sigma_0)$ has finite coproducts.\end{corollary}

Let us write $\CMon$ for $\CMon(\Obj)$ when $\Obj = \{\bullet\}$, ie. there is a single generating object. It is a well-known fact, dating back at least to~\cite{Burroni93} that $\Free_{\CMon}$ is isomorphic to the skeletal\footnote{Skeletal means that we identify isomorphic objects, i.e. sets with the same number of elements.} category $\FinOrd$ of finite sets and functions.

For our purposes it is important to study a slightly different `functional' representation of $\Free_{\CMon}$, called $\FinFun$, which refers more explicitly to the way these structures are implementable. To introduce it, we will first need to give necessary background on arrays and models of parallel computation.

\subsubsection{Finite Functions in Terms of Arrays and Parallel Programs}
\input{background/arrays.tex}

\subsubsection{$\FinFun$: Finite Functions as Arrays}

In addition to the presentation in terms of the category $\Free_\CMon$,
mentioned in \Cref{section:finord},
finite functions may also be represented as \emph{integer arrays}.
This category, which we call $\FinFun$, will be central to our
implementation. Thus we describe it in detail, including complexity results for
several operations.

\begin{definition}[$\FinFun$]
  \label{definition:finfun}
  The category $\FinFun$ has objects
  the natural numbers $\Nat$.
  An arrow $f : A \to B$ is an element $f \in \ord{B}^A$.
  Explicitly, $f$ is an array of values $\bl x_0, x_1, \ldots, x_{A-1} \br$
  where each $x_i \in \ord{B}$.
\end{definition}

We can represent an arrow of $\FinFun$ by its target (codomain) and a table
(array) of elements of its outputs.
The source (domain) of an arrow is the length of its element table.

\begin{minted}{python}
  class FiniteFunction:
    target: int
    table:  array

  @property
  def source(self):
    return len(self.table)
\end{minted}

\begin{proposition}
  \label{proposition:finfun-is-a-category}
  $\FinFun$ forms a category with identities and composition as below.
  \begin{align*}
    \id & : A \to A                       & f \cp g      : A \to C \\
    \id & = \bl 0, 1, \ldots, A-1 \br     & (f \cp g)(i) \mapsto g_{f_i}
  \end{align*}
\end{proposition}

\begin{proof}
  Composition is well-defined:
  $g_{f_i}$ is always defined precisely because $f(i) \in \ord{B}$.
  The identity law is satisfied because $\id_i = i$, so for an arrow $f : A \to B$ we have
  $(\id_A \cp f)_i = \id_{f_i} = f_i$
  and
  $(f \cp \id_B)_i = f_{\id_i} = f_i$
  for all $i$.
  Finally, observe that composition is associative:
  let $A \overset{f}{\to} B \overset{g} \to C \overset{h}{\to} D$
  be arrows. Then we have
  $((f \cp g) \cp h)_i = h_{{(f \cp g)}_i} = h_{g_{f_i}} = (g \cp h)_{f_i} = (f \cp (g \cp h))_i$
\end{proof}

The identity morphism and composition of morphisms can be implemented
as follows.

\begin{minted}{python}
  @staticmethod
  def identity(n: int):
    return arange(0, n) # [0, .., n - 1]

  @staticmethod
  def compose(f, g):
    assert f.target == g.source
    return FiniteFunction(g.target, g.table[f.table])
\end{minted}

\begin{proposition}[Complexity of composition in $\FinFun$]
  \label{proposition:finfun-composition-complexity},
  Let $f : A \to B$ and $g : B \to C$ be morphisms of $\FinFun$.
  Computing the composite $f \cp g$
  has $O(A)$ sequential
  and $O(1)$ PRAM CREW time complexity.
\end{proposition}
\begin{proof}
  The composite $f \cp g : A \to C$ is an integer array of $A$ elements
  in the set $\ord{C}$.
  For each $i \in \ord{A}$,
  Each element $(f \cp g)_i$ is computed by a single memory lookup, for a total of $O(A)$ operations.
  In the parallel (PRAM CREW) setting, each lookup can be performed in parallel,
  giving $O(1)$ time complexity.
\end{proof}

\subsubsection{$\FinFun$ as a strict monoidal category}

The category of finite sets and functions has initial objects and coproducts,
from which it can be made into a strict symmetric monoidal category.

\begin{proposition}
  \label{proposition:finfun-initial}
  $0$ is the initial object in $\FinFun$.
\end{proposition}
\begin{proof}
  Let $B$ be an object in $\FinFun$.
  Then there is a unique morphism $\initial : 0 \to B$:
  the empty array $\bl \br$.
\end{proof}

In code, the initial map $\mathtt{initial} : 0 \to B$ returns the empty array.
\begin{minted}{python}
@classmethod
def initial(cls, B):
  return FiniteFunction(B, zeros(0))
\end{minted}

\begin{proposition}[Coproducts in $\FinFun$]
  \label{proposition:finfun-coproducts}
  Let $f : A_0 \to B$ and $g : A_1 \to B$ be arrows of $\FinFun$.
  The coproduct of objects $A_0$ and $A_1$ is given by addition $A_0 + A_1$,
  with injections defined as
  \begin{align*}
    & \inj_0 : A_0 \to A_0 + A_1             \qquad & \qquad \inj_1 & : A_1 \to A_0 + A_1 \\
    & \inj_0 = \bl 0, 1, \ldots A_0 - 1 \br  \qquad & \qquad \inj_1 & = \bl A_0, A_0 + 1, \ldots A_1 - 1 \br
  \end{align*}

  Then the coproduct $f + g : A_0 + A_1 \to B$, denoted $f + g$, is given by
  array concatenation:
  \begin{align*}
    f + g = \bl f_0, f_1 \ldots f_{A_0 - 1}, g_0, g_1 \ldots g_{A_1 - 1} \br
  \end{align*}
\end{proposition}
\begin{proof}
  This choice of $f + g$ commutes with the injections:
  \[ (\inj_0 \cp (f + g)_i = (f + g)_{\inj_{0_i}} = [f_0, f_1, \ldots, f_{A_0 - 1}]_i = f_i \]
  and
  \[ (\inj_1 \cp [f, g])_i = [f, g]_{\inj_{1_i}} = [g_0, g_1, \ldots, g_{A_1 - 1}]_i = g_i \]
  Moreover, this choice must be unique: if even an entry in the array $f + g$ is
  not as specified above, then the diagram does not commute.
\end{proof}

The coproduct $f + g : A_0 + A_1 \to B$ of maps
$f : A_0 \to B$ and $g : A_1 \to B$ is implemented as array concatenation.
\begin{minted}{python}
def coproduct(f, g):
  assert f.target == g.target
  target = f.target
  table = concatenate(f.table, g.table)
  return FiniteFunction(target, table)
\end{minted}

Naturally, since $\FinOrd$ is cocomplete and we claimed that $\FinOrd \cong \FinFun$,
we expect $\FinFun$ to have coequalizers making it cocomplete as well.
This is indeed the case; a well-known result is below.

\begin{proposition}[Coequalizers in $\FinFun$]
  \label{proposition:finfun-coequalizers}
  Let $f, g : A \to B$ be parallel arrows in $\FinFun$,
  and $G$ the graph of
  $B$ vertices
  and edges $\{ (f(i), g(i)) \mid i \in \ord{A} \}$.
  If $Q$ is the number of connected components of $G$,
  and $q : B \to Q$ is the function labeling a vertex with its connected component,
  Then $q = \coeq(f, g)$ is the coequalizer of $f, g$
\end{proposition}

Coequalizers can be computed directly using the connected components algorithm.
For the purposes of implementation, we assume the existence of a primitive
function
$\mathtt{connected\_components} : B^A \times B^A \to Q \times Q^B$,
which computes connected components from a graph encoded as an adjacency list.
That is, its two arguments are an array of edge sources and targets, respectively.
Then we can implement coequalizers of finite functions as follows.

\begin{minted}{python}
def coequalizer(f, g):
  assert f.source == g.source
  assert f.target == g.target
  Q, q = connected_components(f.table, g.table)
  return FiniteFunction(Q, q)
\end{minted}

\begin{proposition}[Coequalizer Complexity]
  \label{proposition:finfun-coequalizers-complexity}
  Computing the coequalizer of finite functions
  $f, g : A \to B$
  has $O(A + B)$ sequential
  and $O(\log (A + B))$ PRAM CRCW time complexity.
\end{proposition}
\begin{proof}
  Clearly the complexity of computing coequalizers is the same as computing
  connected components.
  In the sequential case, connected components can be labeled in
  $O(V + E)$ time for a graph $G = (V, E)$
  (see e.g. \cite[Chapter 22]{introduction-to-algorithms}),
  and $O(\log V)$ time
  (see~\cite[p. 218]{introduction-to-parallel-algorithms} or \cite{logn-connected-components})
  in the parallel (PRAM CRCW) case.
  Since the graph $G$ has $B$ vertices and $A$ edges, it then follows that
  computing the coequalizer of $f, g : A \to B$ in $\FinFun$ has
  $O(A + B)$ sequential time complexity
  and $O(\log(A + B))$ parallel (PRAM CRCW) time complexity.
\end{proof}

Initial objects and coproducts give $\FinFun$ the structure of a strict monoidal
category.
For morphisms $f : A_0 \to B_0$ and $g : A_1 \to B_1$,
the tensor product $f \tensor g : A_0 \tensor A_1 \to B_0 \tensor B_1$
is given by $(f \cp \inj_0) + (g \cp \inj_1)$
However, tensor products can be written more directly.

\begin{proposition}[Tensor Product]
  \label{proposition:finfun-tensor}
  Let $f : A_0 \to B_0$ and $g : A_1 \to B_1$
  be morphisms in $\FinFun$.
  The tensor product $f \tensor g$ is the array
  $\bl f_0, f_1, \ldots f_{A_0}, (B_0 + g_0), (B_0 + g_1), \ldots (B_0 + g_{A_1} \br$
\end{proposition}
\begin{proof}
  It is straightforward that $f \cp \inj_0$ is an array with the same entries as $f$,
  and that $g \cp \inj_1$
  is an array with entries
  $\bl B_0 + g_0, B_0 + g_1, \ldots B_0 + g_{B_1 - 1}$.
  It is then immediate that the tensor product is the coproduct
  (concatenation) of these two arrays.
\end{proof}

The tensor product of morphisms can be implemented as follows, with time
complexity the same as for the coproduct.
\begin{minted}{python}
def tensor(f, g):
  table = concatenate([f.table, g.table + f.target])
  return FiniteFunction(f.target + g.target, table)
\end{minted}

\begin{proposition}
  \label{proposition:finfun-strict-monoidal}
  $\FinFun$ is a strict symmetric monoidal category.
\end{proposition}
\begin{proof}
  Symmetry is given in the usual way with $\twist \defeq (\inj_1 + \inj_0) = \bl 1, 0 \br$,
  which is evidently self-inverse.
  Strictness follows because concatenation of arrays is strictly associative,
  and concatenation with the empty array is the identity.
\end{proof}

One can easily show (see e.g.,~\cite{acsets}) that there is an isomorphism $\FinFun \cong \FinOrd$, which also respects the symmetric monoidal structure of the two categories.

%% file: background/arrays.tex
The main primitive of our implementation is the \emph{array}, so we begin with
notation.

\begin{definition}[Array]
  An \deftext{array $x$ of length $N$ with elements in $S$}
  is an element of $N^S$.
  The $i^{\mathsf{th}}$ element of an array $x$ is denoted $x_i$,
  defined as long as $x$ has length $N > i$.

  Explicitly, an array $x$
  is an ordered list of $N$ elements
  $\bl x_0, \ldots, x_{N-1} \br$.
  with each $x_i \in S$.
  The empty array is denoted $\bl \br$,
\end{definition}

It will also be useful to have dedicated notation for finite sets.

\begin{definition}[Finite Set]
  Let $A \in \Nat$ be a natural number.
  Then the \deftext{finite set} $\{ 0, 1, \ldots, A - 1 \} $ is denoted $\ord{A}$.
\end{definition}

In the sequential RAM model of computation, a single operation takes a single
timestep.
For example, reading or writing to a memory location.
However, there are a variety of models of \emph{parallel}
computation~\cite{introduction-to-parallel-algorithms}.
In this paper, we will use two variants of the Parallel RAM (PRAM)
model~\cite[p.11, p.15]{introduction-to-parallel-algorithms}.
The first of these is the PRAM CREW (Concurrent Read / Exclusive Write) model,
which we will simply write as PRAM.
This assumes that, while many processors can read a memory location in parallel,
conflicting writes by multiple processors to the same location are forbidden.

A slightly stronger model is the PRAM CRCW (Concurrent Read / Concurrent Write),
which allows multiple processors to write to the same memory location in
parallel.
When conflicting writes occur, an arbitrary processor succeeds.
We will need the full power of the PRAM CRCW only rarely.
Unless specified explicitly, we hereafter only assume the PRAM CREW model.

The algorithms presented in this paper are all in terms of a small number of
primitive integer array operations.
We give a table of these with their sequential and PRAM CREW complexities in
terms of the size of the input $n$ below.
Note that two operations (\texttt{repeat} and \texttt{segmented arange}) have
complexities in terms of their input \emph{values}, and we use $\len{s}$ to
denote the size of the array $s$.

\begin{equation*}
  \label{table:array-primitives}
  \begin{tabular}{lcl}
    Primitive                      & Complexity (Sequential)  & Complexity (PRAM) \\
    \toprule
    \texttt{arange}                & $O(n)$                   & $O(1)$      \\
    \texttt{zeros}                 & $O(n)$                   & $O(1)$      \\
    \texttt{sum}                   & $O(n)$                   & $O(\log n)$ \\
    \texttt{prefix sum}            & $O(n)$                   & $O(\log n)$ \\
    \texttt{dense integer sort}    & $O(n)$                   & $O(\log n)$ \\
    \texttt{concatenate}           & $O(n)$                   & $O(\log n)$ \\
    \texttt{connected components}  & $O(n)$                   & $O(\log n)$ (CRCW) \\
    \midrule
    \texttt{repeat(x, s)}          & $O(sum(s))$              & $O(\log \len{s})$ \\
    \texttt{segmented arange(s)}   & $O(sum(s))$              & $O(\log \len{s})$ \\
    \bottomrule
  \end{tabular}
\end{equation*}

We include under the banner of \texttt{sum} and \texttt{prefix sum} other
operations like \texttt{max} and \texttt{all} which can be implemented with
parallel scans and folds.
For a more in-depth explanation of these primitives, we direct the reader to our
implementation, which can be found at \codelink{}.
However, we give a brief overview now.
The \texttt{dense integer sort} operation refers specifically to the subset of
sorting algorithms operating on positive integer arrays of length $n$ whose
largest element is $O(n)$.
Such sorts can be computed in $O(n)$ sequential time by counting sort, and in
$O(\log n)$ parallel time by radix sort.
The \texttt{concatenate} operation simply copies multiple arrays to a single
contiguous memory location.
There are some subtleties to its use in the parallel case, which we
discuss further in \Cref{section:applying-functors}.

The \texttt{repeat(x, s)} operation takes two equal-length arrays, and outputs
the array whose entries are those of $x$ repeated a number of times indicated by
$s$.
So for example,
$\mathtt{repeat}(\bl a, b, c \br, \bl 0, 1, 2 \br) = \bl b, c, c \br$.
The \texttt{arange} primitive outputs a length $n$ array of indices
$\bl 0, 1, \ldots n - 1 \br$,
and \texttt{segmented arange} computes a concatenation of such arrays whose
lengths are specified by the input argument $s$.
Note also that \texttt{segmented arange} can in fact be expressed in terms of
the other operations, and so is not required to be a primitive.

Finally, note that the time complexity of most operations is at most $O(n)$
(sequential) and $O(\log n)$ (PRAM CRCW).
Most of the algorithms we present later will be in terms of a constant number of
each of such operations, thus guaranteeing linear (sequential) and logarithmic
(parallel) time complexity.

With our model of parallel computation chosen, we can now give basic complexity
results for operations on finite functions defined as arrays.

%% file: background/cospans.tex
Cospans (and more recently \emph{structured} cospans~\cite{structured-cospans})
are now commonly used to model `open' systems.
In particular, cospans are used in the combinatorial representation of string
diagrams introduced in \cite{rmsms}, and illustrated in
Section~\ref{section:hypergraphs} below.
We therefore recall these concepts now.

\begin{definition}
  A \deftext{cospan} in a category $\cat{C}$ is a pair of morphisms
  $\cospan{X}{s}{W}{t}{Y}$.
  We call the maps $s$ and $t$ the \deftext{legs} of the cospan, the objects $X$
  and $Y$ the \deftext{feet}, and $W$ the \deftext{apex}.
\end{definition}

\emph{Structured} cospans are a recently introduced~\cite{structured-cospans}
double-categorical framework for open systems.
The work we present here can be interpreted as a re-examination of the
combinatorial characterisation of string diagrams in \cite{rmsms} through the
lens of this framework.

\begin{definition}[from \cite{structured-cospans}]
  Let $\SCL : \cat{C} \to \cat{D}$ be a functor,
  which we call the \deftext{structuring functor}.
  A \deftext{structured cospan} is a cospan in $\cat{D}$
  of the form
  $\cospan{\SCL(A)}{s}{W}{t}{\SCL(B)}$.
\end{definition}

We will not require the double-categorical structure of structured cospans here.
However, by~\cite[Corollary 3.11]{structured-cospans} structured cospans form a
symmetric monoidal category when $\SCL$ is left-adjoint and $\cat{C}$ and
$\cat{D}$ have finite colimits.
Composition in this category is by pushout, and tensor product is the pointwise
tensor product of the `legs' of the cospan.

\begin{remark}
  Closely related, but not used in this paper, is the construction of
  \emph{decorated} cospans~\cite{decorated-cospans}.
  The relationship between decorated and structured cospans is examined in
  \cite{structured-vs-decorated}.
\end{remark}

%% file: background/hypergraphs.tex
The authors of \cite{rmsms} provide the a combinatorial characterisation of
string diagrams as cospans of hypergraphs.
Although we will not require the technical details in this paper, we now give a
quick sketch of the basic idea, which is useful to keep in mind for later developments.
A string diagram (left) and its hypergraph representation (right) are depicted
below.
\[ \input{examples/string-diagram.tikz} \qquad \qquad \qquad \qquad \input{examples/string-diagram-as-hypergraph.tikz} \]
In this representation, \emph{hyperedges} (depicted as white squares) are
labeled with generators, and \emph{nodes} (depicted as $\wNode$) represent the
wires of the string diagram.
Importantly, the hyperedges in this representation have \emph{ordered lists} of
nodes as their sources and targets. %
Being a cospan, the representation actually consists of three distinct hypergraphs. The hypergraph in the middle (on grey background) encodes the internal structure of the string diagram. The outermost hypergraphs (on blue background) are discrete, i.e. they only containing nodes. The two morphisms from the outermost hypergraphs to the middle one, whose assignments are indicated with dotted arrows, indicate which nodes constitute the left and the right interface of the hypergraph. This corresponds to the `dangling wires' of the string diagrams, which allow sequential composition on the left and on the right with other string diagrams.

As illustrated in~\cite{rmsms}, this representation is known to be an isomorphism of categories when the string diagrams come from a hypergraph category --- so that there is a Special Frobenius structure on each object. When considering string diagrams in a symmetric monoidal category, without extra structure, then the representation is an isomorphism only if we restrict to so-called \emph{monogamous} cospans. Monogamicity is a requirement on the way interfaces, nodes, and hyperedges interact with each other; we refer to~\cite{rmsms} for a full definition,
and give an equivalent one in \Cref{section:symmetric-monoidal}.

However, in order to work with such hypergraphs on a computer, one must choose
how their data should be represented.
For example, one might choose to represent each hyperedge as a pair of ordered
integer lists.
However, in order to define parallel algorithms suitable for e.g., GPUs, we will
need a `flat' array representation.
The issue of how to \emph{efficiently} encode these hypergraphs is addressed in
\cite{coc}, where the authors show that encoding (monogamous) hypergraphs as
bipartite graphs leads to an efficient representation as sparse matrices.
However, \cite{coc} is limited to the \emph{monogamous} case, and so is only
suitable for representing string diagrams of symmetric monoidal categories
without the additional Special Frobenius structure.
Further, \cite{coc} only accounts for categories with a single generating objects, and the
representation described does not explicitly represent diagrams as cospans.

Our work generalises the approach of \cite{coc} to the case of string diagrams equipped with
a chosen Special Frobenius structure.
This addition allows us to define in \Cref{section:optic-composition} an
efficient algorithm for taking reverse derivatives of large diagrams,
Returning to our running example, the addition of this structure allows for the
representation of string diagrams such as the one below left as cospans of
bipartite multigraphs (below right).
\[ \input{examples/string-diagram.tikz} \qquad \qquad \qquad \qquad \input{examples/string-diagram-combinatorial.tikz} \]

\subsection{$ACSet$s}

One of our main contributions is to encode the hypergraph structure of \cite{rmsms}
in such a way that allows for \emph{parallel} algorithms on the datastructure.
This is possible thanks to the construction of $\ACSet$s, first described in
\cite{acsets}.
Informally, $\ACSet$s are a class of category which represent `array-based data
with attributes'.
For example, a graph with node and edge labels has data (the adjacency
information of the graph) and attributes (the labels).
We now informally recall the $\ACSet$s of~\cite{acsets}, beginning with the
definition of a \emph{schema}.

\begin{definition}[Schema \cite{acsets}, Informal]
  A \deftext{schema} is a finitely presented (small) category $S$
  which is `bipartite': objects of $S$ are either in $S_0$ or $S_1$,
  respectively.
\end{definition}

Intuitively, a schema specifies the relationships between data and attributes of
a category.
Concretely, the objects in $S_0$ (the data) will map to finite sets, and those in
$S_1$ (the attributes) to some chosen typing.

\begin{definition}[$\ACSet$s \cite{acsets}]
  Given a schema $S$ and a typing map $K : S_1 \to \Set$,
  the category $\ACSet^S_K$ has objects
  functors $F : S \to \Set$
  which restrict to $F|_{S_1} = K$
  and arrows $\alpha : F_0 \to F_1$
  the natural transformations where $\alpha|_{S_1} = \id$.
\end{definition}

The \emph{typing map} $K$ defines the type of attributes of data.
For example, a graph with vertices labeled in $\Nat$ might be modeled as an
$\ACSet$ whose schema has three objects: $S_0 = \{V, E\}$ and $S_1 = \{ L \}$,
and whose typing map is $K(L) = \Nat$.
The requirement that arrows $\alpha$ must have $\alpha|_{S_1} = \id$
ensures that morphisms of $\ACSet$s cannot arbitrarily modify attribute data.

In this paper, we will only need to consider \emph{finite} $\ACSet$s: those for
which the objects $F$ restrict to functors $F|_{S_0} \to \FinFun$.
Notice that this means that our datastructure consists of dense integer arrays:
this is critical for making our implementation suitable for parallel hardware
such as GPUs.
In the next section, we will see two examples of $\ACSet$s which form the basis
of our datastructure.

%% file: representing-diagram-wirings.tex
We can now begin to define our datastructure for representing string diagrams as
structured cospans.
In this section, we define two categories, $\Wires$ and $\BPMG$, which model the
`internal wiring' of a string diagram.
These categories are related by an adjunction $\SCL : \Wires \to \BPMG$ which
will serve as the `structuring functor'
when we come to define string diagrams as structured cospans
$\cospan{\SCL(A)}{s}{G}{t}{\SCL(B)}$
in \Cref{section:string-diagrams}.
The objects of $\Wires$ will serve as the feet of the cospan and are analogous
to the discrete hypergraphs of \cite{rmsms}.

To make the role of the two categories more clear, we begin with examples.
Pictured below left is a string diagram,
and below right the bipartite multigraph corresponding its
`internal wiring'.
\[
  \gen{examples/string-diagram}
  \qquad \qquad \qquad \qquad
  \gen{examples/string-diagram-bpmg}
\]
The bipartite structure means that there are two kinds of node.
The $\wNode$-nodes represent the \emph{wires} of the string diagram, and are
labeled with the generating objects of $\Obj$.
Meanwhile, the $\xNode$-nodes represent operations, and are labeled in $\Arr$.

Edges of the graph are labeled with a natural number:
An edge $\wNode \to \xNode$ labeled $i$ denotes that a given wire connects to
the $i^{\text{th}}$ source port of an operation.
Similarly, an edge $\xNode \to \wNode$ denotes that a wire connects to a given
\emph{target} port of an operation.

Notice however that this information is missing.
Namely, it is not specified which of the $\wNode$-nodes correspond to the
\emph{boundary} of the string diagram; that is, which are the left and right
`dangling wires'.
This is the purpose of the category $\Wires$, whose objects are thought of as
labeled finite sets, and whose morphisms are `label preserving maps'.
Pictured below is an $\SCL$-structured cospan.
 \[ \gen{examples/string-diagram-combinatorial} \]
The center grey box depicts the apex of the cospan: a bipartite multigraph.
At the edges are two objects of $\Wires$ in blue boxes.
Dashed arrows between them correspond to the legs of the structured cospan, and
can be thought of as morphisms of $\Wires$.

To formalise this construction, we now give the precise details of the
categories $\Wires$ and $\BPMG$ as $\ACSet$s, as well as the structuring functor
$\SCL$.

%% file: adding-colour.tex
We now describe $\Wires$, the category of `multi-sorted' finite sets and
functions.
As with $\FinOrd$, these functions are presented by $\CMon(\Sigma_0)$.
The only difference is that the set of generating objects $\Obj$ is now
arbitrary, instead of being restricted to a single object.
For example, given generating objects $A$ and $B$, we can construct morphisms
like the one below.
\[ \input{examples/multi-sorted-function.tikz} \]

The \emph{cospans} of such morphisms (below right) correspond to string diagrams
with generators of $\Frob$ like the one below left.
\[ \input{examples/multi-sorted-free-hypergraph.tikz} \qquad \qquad \qquad \qquad \input{examples/multi-sorted-free-hypergraph-combinatorial.tikz} \]
Note carefully that the legs of the cospan are now also
\emph{label preserving}.
That is, the source and target of each dashed arrow have the same label.
This property arises by defining $\Wires$ as a category of $\ACSet$s where
labels are attributes.

\begin{definition}
  Let $\Signature$ be an arbitrary monoidal signature.
  The category of `labeled wires' is denoted $\WiresFrob$
  and defined as the category $\ACSet^{\SchemaWires}_{K_\Signature}$
  whose schema $\SchemaWires$
  is defined as follows:
  \[ \SchemaWires \defeq \input{schema-wires.tikz} \]
  and whose typing map $K_\Signature$ is defined
  \[ K_\Signature(\Ob) \defeq \Obj \]
  Given a `labelling function' $l : B \to \Obj$,
  we write $\MkWiresObj(l)$ for the object of $\Wires$
  with $\MkWiresObj(l)(W) = B$
  and $\MkWiresObj(l)(\swn) = l$.
  Given a finite function $f : A \to B$
  we write $\MkWiresArr(f, l)$
  for the \emph{morphism} of $\Wires$
  $f : \MkWiresObj(f \cp l) \to \MkWiresObj(l)$
\end{definition}

\begin{remark}
  Arrows $\alpha : F_0 \to F_1$ in $\Wires$ are natural transformations with a single component,
  $\alpha_W$.
  To reduce the notational burden, we will sometimes omit this subscript.
  However, one should think of these maps as `label-preserving finite functions'
  between labeled finite sets.
\end{remark}

\begin{example}
  Let $F$ and $G$ be objects of $\WiresFrob$.
  Concretely, $F$ consists of a set $F(W)$ and a morphism $F(\swn) : F(W) \to \Obj$.
  Similarly, $G$ is a set $G(W)$ and morphism $G(\swn) : G(W) \to \Obj$.
  A label preserving map $\alpha$ between $F$ and $G$ is a natural
  transformation with a single non-identity component, $\alpha_W$.
  An example of a label-preserving map for $F(W) = 3$ and $G(W) = 3$
  is given below.
  \[ \input{examples/label-preserving-map.tikz} \]
  Note that the finite function $\alpha_W$ can be depicted formally as the string
  diagram below, where we have labeled incoming wires according to $F(\swn)$ and
  outgoing wires according to $G(\swn)$.
  \[ \alpha_W \quad \defeq \quad \input{examples/label-preserving-map-finite-function.tikz} \]
  This example illustrates how morphisms of $\Wires$ indeed correspond to the
  theory of commutative monoids over an arbitrary set of generating objects $\Obj$.
\end{example}

We can alternatively think of the objects of $\Wires$ as being discrete
vertex-labeled graphs.
When we come to define the structuring functor $\SCL$, an object $A \in \Wires$
will map to the \emph{discrete} bipartite multigraph $\SCL(A) \in \BPMG$ having
only $\wNode$-nodes.

In fact, morphisms in $\Wires$ will form the `Frobenius Half-Spiders' in an
arbitrary Hypergraph category.

\begin{proposition}
  Let $\cat{C}$ be a hypergraph category presented by a signature $\Signature$.
  There is a unique strict symmetric monoidal identity-on-objects functor
  $\HalfSpider : \WiresFrob \to \cat{C}$ which maps the monoid structure of
  $\CMon$ to the monoid structure of $\Frob$.
\end{proposition}
\begin{proof}
  Recall that $\FinOrd$ is presented by $\CMon$,
  and the morphisms of $\Wires$ are natural transformations with a single
  non-identity component in $\FinOrd$.
  It follows that the morphisms of $\WiresFrob$ are also presented by $\CMon$
  with the set of generating objects $\Obj$.
  Thus, the action of $\HalfSpider$ is fixed on all generators, composition, and
  tensor products, and is unique.
\end{proof}

We call such morphisms `Half-Spiders'.

\begin{definition}[Half-Spider]
  \label{definition:half-spider}
  Let $\cat{C}$ be a hypergraph category presented by a signature $\Signature$.
  A \deftext{Frobenius Half-Spider} in $\cat{C}$ is a morphism in the image of
  $\HalfSpider : \WiresFrob \to \cat{C}$.
  We will sometimes write $\HalfSpider(f)$ for an arrow $f \in \FinFun$
  to mean $\HalfSpider(\MkWiresArr(f, l))$ when $l$ is clear from context.
\end{definition}
Alternatively, half-spiders are those morphisms in a hypergraph category built
by tensor and composition from the generators
$
  \gen{g-identity},
  \gen{g-twist},
  \gen{frob/g-join},
$ and $
  \gen{frob/g-unit}
$.
In addition, recall that all hypergraph categories have a
dagger~\cite[Section 1.3.3]{fong-thesis}.
Consequently, there is a second functor mapping morphisms of $\Free_\CMon$ to
the \emph{comonoid} structure of a hypergraph category:
$f \mapsto \HalfSpider(f)^\dagger$.
Thus, these `dagger half-spiders' are morphisms built from
$
  \gen{g-identity},
  \gen{g-twist},
  \gen{frob/g-split},
$ and $
  \gen{frob/g-counit}
$.

Putting these two functors together characterises the Frobenius spiders
in a Hypergraph category.
Namely, Frobenius spiders will be those morphisms of the form
$\HalfSpider(f) \cp \HalfSpider(g)^\dagger$.
This will be stated more formally in
\Cref{proposition:frobenius-spiders-in-diagram}.

We will now describe bipartite multigraphs, accounting for the additional
structure required to represent the \emph{operations} in a string diagram.

%% file: bipartite-multigraphs.tex
Bipartite Multigraphs can be regarded as objects of $\Wires$ decorated with
additional data.
More formally, they are given by the following $\ACSet$.

\begin{definition}
  A \deftext{bipartite multigraph} is an object
  of the category $\BPMG \defeq \ACSet^{\SchemaBPMG}_K$, where
  $\SchemaBPMG$ is the schema functor defined below
  \[ \SchemaBPMG \qquad \defeq \qquad \input{schema.tikz} \]
  and the typing map $K : \SchemaBPMG_1 \to \Set$ is defined as
  \begin{align*}
    K(\Ob)   \defeq \Obj \qquad \qquad
    K(\Sig)  \defeq \Arr \qquad \qquad
    K(\Port) \defeq \Nat
  \end{align*}
\end{definition}

Not every bipartite multigraph represents a valid `internal wiring' of a string
diagram.
Note for example that the schema above does not prevent a bipartite multigraph
with a generator having two edges with the same `port label' from a $\wNode$ to
a $\xNode$-node.
\begin{example}
  The following bipartite multigraph is `ill-formed' for two reasons.
  \[ \gen{examples/bipartite-multigraph-ill-formed} \]
  First, there are two $\wNode \to \xNode$ edges labeled $0$.
  This conflicts with the interpretation of edge labels as defining which `port'
  of a generator is connected to a wire.
  Secondly, there is no output edge $\xNode \to \wNode$ labeled $0$: we need
  every output port of $g$ to be accounted for.
\end{example}

In order to rule out objects like the above, we will need to speak of
`well-formed' bipartite multigraphs.

\begin{definition}[Well-Formed Bipartite Multigraph]
  \label{definition:well-formed}
  A bipartite multigraph $G$ is
  \deftext{well-formed with respect to a signature} $\Signature$
  there are chosen typings %
  $(a : X \to \Obj^*, b : X \to \Obj^*)$
  so that for each $x \in G(X)$
  so that $(a(x), b(x)) \in \typer(G(\swn)(x))$,
  both $\tl \sxi, \sporti \tr$ and $\tl \sxo, \sporto \tr$ are mono,
  and
  \[
    \forall e \in \Ei \quad \swn(\swi(e)) = a(\sxi(e))_{\sporti(e)}
    \qquad \qquad
    \forall e \in \Eo \quad \swn(\swo(e)) = b(\sxo(e))_{\sporto(e)}
  \]
\end{definition}

\begin{remark}
  A consequence of \Cref{definition:well-formed}
  is that in a well-formed bipartite multigraph $G$,
  $G(\Ei)$ and $G(\Eo)$ equal the total arity and coarity
  of generators $G(X)$, respectively.
\end{remark}

\begin{remark}
  The above definition of well-formedness specifically allows for polymorphic
  generators.
  More precisely, since an operation $g \in \Arr$ has a typing \emph{relation},
  we may have for example a bipartite multigraph like the following.
 \[ \gen{examples/bipartite-multigraph-polymorphic} \]
  Notice that there are two $\xNode$-nodes labeled $g$ with different types:
  one $g : A \to \unit$,
  and the other $g : B \otimes C \to D$.
  Such a diagram is considered well-formed as long as both such types exist in
  the typing relation.
\end{remark}

When defining algorithms on structured cospans of well-formed bipartite multigraphs,
we will need the following preservation theorem.
This will guarantee that composites of structured cospans built from well-formed
diagrams are also well-formed.

\begin{proposition}[Preservation of Well-Formedness]
  \label{proposition:bpmg-morphisms-preserve-well-formed}
  If a morphism of bipartite multigraphs $\alpha : F_0 \to F_1$
  has all components permutations except for $\alpha_W$,
  and $F_0$ is well-formed,
  then $F_1$ is well-formed.
\end{proposition}
\begin{proof}
  Observe that the number of edges, generators, and their attribute data is
  preserved.
  Moreover, the attributes of $\wNode$-nodes are also preserved by naturality,
  so $F_1$ is well-formed.
\end{proof}

Finally, before defining string diagrams as structured cospans, we first
formalise the relationship between the categories $\Wires$ and $\BPMG$ as an
adjunction.
Observe that $\Wires$ is a subcategory of $\BPMG$, and therefore embeds into it.

\begin{definition}
  Denote by $\SCL : \Wires \to \BPMG$ the (identity-on-objects) inclusion
  functor.
\end{definition}

As observed in \cite[p.18]{acsets}, $\SCL$ is left-adjoint to the forgetful
functor.

\begin{proposition}[$\SCL$ is left-adjoint]
  \label{proposition:scl-left-adjoint}
  Let $\SCR : \BPMG \to \Wires$ be the forgetful functor
  mapping a bipartite multigraph to its set of wires and their labels.
  Then $\SCL$ is left-adjoint to $\SCR$.
\end{proposition}
\begin{proof}
  Let $\alpha : \SCL(A) \to G$ be an arrow in $\BPMG$,
  and let $\iW_G : \SCL(\SCR(G)) \to G$ be the natural transformation
  whose components are defined as follows.
  \begin{align*}
    \iW_{G_Y} & = \begin{cases}
      \id_W    & \text{if $Y = W$} \\
      \initial & \text{otherwise}
    \end{cases}
  \end{align*}
  where $\initial = \sgen{cmon/g-zero}$ is the initial map.
  Then there is a unique $f$ such that
  $\SCL(f) \cp \iW_G = \alpha$.
  Since composition of morphisms in $\BPMG$ is pointwise, we must have that
  $(\SCL(f) \cp \iW_G)_Y = \alpha_Y$
  for each component $Y$.
  We therefore must have
  $\SCL(f)_W = \alpha_W$,
  and
  $\SCL(f)_Y = \initial : 0 \to 0$,
  and so $f = \MkWiresArr(\alpha_W, G(\swn))$
  by naturality.
\end{proof}

The left-adjointness of $\SCL$ says that any morphism of bipartite multigraphs
of the form $\alpha : \SCL(A) \to G$ is completely determined by its $W$
component
$\alpha_W : \SCL(A)(W) \to G(W)$.
When we define string diagrams in \Cref{section:string-diagrams}, this fact will
allow us to define composition in terms of coequalizers of $\FinFun$.
In addition, the bipartite multigraphs in the image of $\SCL$ are particularly
important in the next section.
We therefore introduce the following notation.
\begin{definition}[Discrete Bipartite Multigraph]
  \label{definition:discrete-bipartite-multigraph}
  Given a labeling $w : A \to \Obj$,
  the \deftext{discrete bipartite multigraph} for $w$
  is denoted $\MkD{w} \defeq \SCL(\MkWiresObj(w))$.
  Given a bipartite multigraph $G$, we overload the same notation and write
  $\MkD{G} \defeq \SCL(\SCR(G))$.
\end{definition}

With these definitions in hand, we can finally define string diagrams as
structured cospans.

%% file: string-diagrams.tex
We can now define our category of string diagrams as structured cospans.

\begin{definition}[Diagram]
  \label{definition:diagram}
  Fix a monoidal signature $\Signature$.
  A \deftext{diagram} over $\Signature$
  is an $\SCL$-structured cospan
  $\cospan{\SCL(A)}{s}{G}{t}{\SCL(B)}$
  where $G$ is well-formed with respect to $\Signature$.
  We say the \deftext{type} of such a diagram is $A \to B$.
\end{definition}

Note that diagrams form a symmetric monoidal category
as described in~\cite[Corollary 3.11]{structured-cospans}.
This follows because $\SCL$ is left-adjoint and therefore preserves colimits.

\begin{definition}
  \label{definition:diagram-frob}
  The symmetric monoidal category of diagrams is denoted $\DiagramFrob$
  and defined to be the subcategory of
  $_{\SCL}\Csp(\BPMG)$
  defined in \cite[Corollary 3.11]{structured-cospans}
  whose objects are those of $\Wires$,
  and whose arrows are isomorphism classes of diagrams.
\end{definition}

In order to define $\DiagramFrob$ as a \emph{sub}category, it is important to
verify that composition and tensor product preserve well-formedness of diagrams.
This follows from \Cref{proposition:bpmg-morphisms-preserve-well-formed}, and we
give a proof in \Cref{proposition:composition-preserves-well-formedness}.

Composition of structured cospans is by pushout, and tensor product is
pointwise; data-parallel algorithms for both will shortly be given in
Propositions
\ref{proposition:composition} and
\ref{proposition:tensor},
respectively.
In fact, the remainder of this section is dedicated to showing
(1) how to represent diagrams efficiently;
(2) how to construct `primitive' diagrams like identity and symmetry;
and
(3) data-parallel algorithms for tensor and composition.
In addition, the aim of this section is to serve as a guide to implementing the
algorithms described.
We therefore include several pseudo-Python code listings in several cases.

\subsection{Representing Diagrams Efficiently}

The naive representation of \Cref{definition:diagram} would require a
full bipartite multigraph for each of the feet of the cospan as well as a
complete morphism of bipartite multigraphs for the legs.
However, in order to represent large diagrams it is important to be economical
with data.
We therefore exploit that $\SCL$ is left-adjoint in order to show that much of
this data is redundant.
This will also simplify the definition of several algorithms in the rest of the
paper.

\begin{proposition}
  \label{proposition:representation}
  There is a bijective correspondence between
  diagrams and
  triples
  $(s, t, G)$
  where
  $\cospan{A}{s}{G(W)}{t}{B}$
  is a cospan in $\FinFun$.
\end{proposition}
\begin{proof}
  We first show that any diagram uniquely defines such a triple.
  Let $\cospan{\SCL(A)}{\sigma}{G}{\tau}{\SCL(B)}$ be a diagram.
  Since $\SCL$ is left-adjoint (\Cref{proposition:scl-left-adjoint}),
  there are unique morphisms $s, t$ of $\WiresFrob$
  such that
  $\SCL(s) \cp \iW_G = \sigma$
  and
  $\SCL(t) \cp \iW_G = \tau$.
  By virtue of being morphisms of $\WiresFrob$, $s$ and $t$ are natural
  transformations with a single non-identity component defined on $W$.
  We may therefore think of them as morphisms of $\FinFun$.

  In the reverse direction, suppose that
  $(s, t, G)$ is a triple where
  $\cospan{A}{s}{G(W)}{t}{B}$ is a cospan in $\FinFun$.
  This defines a diagram
  $\cospan{\SCL(A)}{\sigma}{G}{\tau}{\SCL(B)}$
  in the following way.
  We must take
  $\sigma = \SCL(\MkWiresArr(s, G((\swn)))) \cp \iW_G$
  and
  $\tau = \SCL(\MkWiresArr(t, G(\swn))) \cp \iW_G$
\end{proof}

This correspondence means that diagrams can be written in code as the following
\texttt{Diagram} class.

\begin{minted}{python}
  class Diagram:
    s: FiniteFunction
    t: FiniteFunction
    G: BipartiteMultigraph
\end{minted}

We will now use this alternate representation to define primitive diagrams and
algorithms for tensor and composition.
Note that in the remainder of this section we take the 2-categorical
perspective, and consider specific diagrams (i.e., 1-cells) rather than
isomorphism classes of diagrams.
This is more natural from the perspective of the programmer working with
diagrams, since one typically works with a representative of an isomorphism
class instead of the class itself.

\subsection{Primitive Diagrams}

\begin{proposition}
  \label{proposition:diagram-identity}
  Given an object $\MkWiresObj(w) \in \Wires$,
  The \deftext{identity diagram} is represented by the triple
  $(\id, \id, \MkD{w})$
  where $\MkD{w}$ is the
  discrete bipartite multigraph (\Cref{definition:discrete-bipartite-multigraph}).
\end{proposition}
\begin{proof}
  For an object $A = \MkWiresObj(w)$ of $\Wires$,
  the identity structured cospan is
  $\cospan{\SCL(A)}{\id}{\SCL(A)}{\id}{\SCL(A)}$.
  We must therefore have $s = t = \id$,
  and $\SCL(A) = \MkD{\SCL(A)}$ by definition.
\end{proof}

In code, the identity diagram is constructed as follows.

\begin{minted}{python}
@classmethod
def identity(cls, wn: FiniteFunction):
    s = FiniteFunction.identity(wn.source)
    t = FiniteFunction.identity(wn.source)
    G = BipartiteMultigraph.discrete(wn)
    return Diagram(s, t, G)
\end{minted}

The symmetry of $\DiagramFrob$ as given as follows.

\begin{proposition}[Symmetry/Twist Diagram]
  \label{proposition:diagram-twist}
  Let $a : A \to \Obj$ and $b : B \to \Obj$
  be labelings,
  then the \deftext{symmetry diagram}
  is given by
  $(\id_{A + B}, \twist_{A, B}, \MkD{a + b})$.
\end{proposition}
\begin{proof}
  The symmetry structured cospan is given by
  $\cospan{\SCL(A + B)}{\id_{A + B}}{\SCL(A + B)}{\twist_{A,B}}{\SCL(B + A)}$,
  and so the result follows as in \Cref{proposition:diagram-identity}.
\end{proof}

$\DiagramFrob$ is also a Hypergraph category,
which means that every object is equipped with a Special Frobenius monoid.
This is the reason we denote the category $\DiagramFrob$.
In \Cref{section:symmetric-monoidal}, we will see how the combinatorial
condition of \emph{monogamicity} introduced in \cite{rmsms} allows us to define
the category $\DiagramSMC$, without this additional structure.

\begin{proposition}
  \label{proposition:diagram-is-hypergraph-category}
  $\DiagramFrob$ is a Hypergraph category
\end{proposition}
\begin{proof}
  $\SCL$ is left-adjoint and therefore preserves colimits.
  $\DiagramFrob$ is then a hypergraph category by
  \cite[Theorem 3.12]{structured-cospans}
\end{proof}

The Frobenius spiders
(\Cref{definition:frobenius-spider})
in $\DiagramFrob$
can be explicitly characterised as those cospans
whose apexes are discrete hypergraphs with no operations.
That is, diagrams of the form $(s, t, \SCL(A))$.

\begin{proposition}
  \label{proposition:frobenius-spiders-in-diagram}
  A diagram $d = (s, t, G)$ is a Frobenius spider iff $G = \SCL(A)$ for some $A$.
\end{proposition}
\begin{proof}
  Deferred to \Cref{section:cospan-proofs}
\end{proof}

In code, we provide a constructor for Frobenius spiders as follows.

\begin{minted}{python}
def spider(s, t, w):
  G = BipartiteMultigraph.discrete(w)
  return Diagram(s, t, G)
\end{minted}

Being a hypergraph category, $\DiagramFrob$ is equipped with a dagger in the
following way.

\begin{proposition}
  \label{proposition:hypergraph-category-has-dagger}
  Hypergraph categories have a dagger given by
  \input{frob/s-dagger-f.tikz}
\end{proposition}
\begin{proof}
  See \cite[Section 1.3.3]{fong-thesis}.
\end{proof}

The dagger of diagrams has a particularly convenient form:
it swaps the legs of the cospan.
That is, $(s, t, G)^\dagger = (t, s, G)$.
We give a proof of this fact in \Cref{proposition:dagger-swaps-source-target},
but note here that it leads to the following simple implementation.

\begin{minted}{python}
def dagger(self):
  return Diagram(self.t, self.s, self.G)
\end{minted}

The final `primitive' diagram we need to define is the \emph{singleton},
representing a single generator of $\Arr$.

\begin{definition}[Singleton Diagram]
  \label{definition:diagram-singleton}
  Let $x : 1 \to \Arr$ be an operation, and
  $a : A \to \Obj$ and $b : B \to \Obj$
  a typing of $x$.
  The \deftext{singleton} diagram of $x$
  is denoted
  $\singleton(a, b, x) \defeq (\inj_0, \inj_1, G)$
  where
  \[ G(W) = A + B \qquad G({\Ei}) = A \qquad G({\Eo}) = B \qquad G(X) = 1 \]
  and
  \begin{align*}
    G(\swi) = \inj_0
    \qquad
    G(\swo) = \inj_1
    \qquad
    G(\sxi) = \terminal_{A}
    \qquad
    G(\sxo) = \terminal_{B}
  \end{align*}
  \begin{align*}
    G(\sporti) = \id_{A}
    \qquad
    G(\sporto) = \id_{B}
    \qquad
    G(\swn) = a + b
    \qquad
    G(\sxn) = x
  \end{align*}
  We will sometimes write just $\singleton(x)$ when a chosen typing is assumed.
\end{definition}

\begin{example}
  Given a generator $f : A \to B \tensor C$
  represented as the element $x = \bl f \br$,
  and typing $a = \bl A \br$ and $b = \bl B, C \br$,
  the singleton diagram $\singleton(x)$
  is depicted below.
  \[ \input{examples/singleton-diagram-bipartite-multigraph.tikz} \]
\end{example}

This definition translates directly to code as follows.

\begin{minted}{python}
@classmethod
def singleton(cls, a: FiniteFunction, b: FiniteFunction, xn: FiniteFunction):
    G = BipartiteMultigraph(
        wi = F.inj0(a.source, b.source),
        wo = F.inj1(a.source, b.source),
        xi = F.terminal(a.source),
        xo = F.terminal(b.source),
        wn = a + b,
        pi = F.identity(a.source),
        po = F.identity(b.source),
        xn = xn)

    return Diagram(
        s=FiniteFunction.inj0(a.source, b.source),
        t=FiniteFunction.inj1(a.source, b.source),
        G=G)
\end{minted}

\subsection{Tensor Product of Diagrams}

\begin{proposition}[Tensor Product of Diagrams]
  \label{proposition:tensor}
  The tensor product of diagrams
  $c_0 \defeq (s_0, t_0, G_0)$
  and
  $c_1 \defeq (s_1, t_1, G_1)$
  is the diagram %
  $c_0 \tensor c_1 \defeq (s_0 \tensor s_1, t_0 \tensor t_1, G_0 \tensor G_1)$
\end{proposition}
\begin{proof}
  Let $\cospan{\SCL(A_0)}{\sigma_0}{G_0}{\tau_0}{\SCL(B_0)}$
  and $\cospan{\SCL(A_1)}{\sigma_1}{G_1}{\tau_1}{\SCL(B_1)}$
  be the respective structured cospans corresponding to $c_0$ and $c_1$ by
  \Cref{proposition:representation}.
  Then then tensor product $c_0 \tensor c_1$ is the following structured cospan.
  \[ \cospan{\SCL(A_0 \tensor A_1)}{\sigma_0 \tensor \sigma_1}{G_0 \tensor G_1}{\tau_0 \tensor \tau_1}{\SCL(B_0 \tensor B_1)} \]
  Since colimits are pointwise, we have
  $(\sigma_0 \tensor \sigma_1)_W = s_0 \tensor s_1$ and $(\tau_0 \tensor \tau_1)_W = t_0 \tensor t_1$
  and so the representation of the tensor product $c_0 \tensor c_1$ is given by
  $(s_0 \tensor s_1, t_0 \tensor t_1, G_0 \tensor G_1)$
\end{proof}

Once again, this yields a straightforward implementation:

\begin{minted}{python}
def tensor(c0, c1):
  return Diagram(
      s = c0.s @ c1.s,
      t = c0.t @ c1.t,
      G = c0.G @ c1.G)
\end{minted}

The sequential time complexity of tensor product is linear (sequential) and
logarithmic (parallel) in the diagram components.

\begin{proposition}
  \label{proposition:complexity-tensor}
  Let $d_i = (s_i, t_i, G_i)$ be diagrams of type $A_i \to B_i$
  for $i \in \{0, 1\}$.
  Computing the tensor product $(s, t, G) = d_0 \tensor d_1$
  has sequential time complexity
  \[ O(A_0(W) + A_1(W))
    + O(G(W))
    + O(G(\Ei))
    + O(G(\Eo))
    + O(G(X))
    + O(B_0(W) + B_1(W))
  \]
  and PRAM time complexity $O(1)$.
\end{proposition}
\begin{proof}
  Computing the coproduct of bipartite multigraphs has sequential time complexity
  $O(G(W)) + O(G(\Ei)) + O(G(\Eo)) + O(G(X))$
  since one must essentially just concatenate arrays.
  The PRAM time complexity of the same coproduct is $O(1)$.
  Then computing $s_0 + s_1$
  is $O(A_0(W) + A_1(W))$ sequential and $O(1)$ PRAM time complexity,
  and similarly for $t_0 + t_1$.
\end{proof}

It will later be convenient to have an efficient method of tensoring a large
number of operations.
In particular, simply repeating the binary tensor operation has poor complexity:
quadratic in the sequential case and linear in the parallel case.
We therefore give the following theorem, which gives a simpler closed form for
the tensor product of $n$ operations.

\begin{proposition}
  \label{proposition:n-fold-tensor-isomorphism}
  Let $x : N \to \Arr$ be a list of $N$ operations,
  and for each operation $x(i)$
  let $(a(i), b(i)) \in \typer(x)$
  be a chosen typing
  so that $a, b : N \to \Obj^*$.
  Let $\len{a(i)}$ and $\len{b(i)}$
  denote the arity and coarity of $x(i)$ respectively.
  Then the $N$-fold tensor product of operations $x$ is isomorphic
  to the diagram $(\inj_0, \inj_1, G)$, where
  \[
    G(W) = \Ki + \Ko
    \qquad
    G(\Ei) = \Ki
    \qquad
    G(\Eo) = \Ko
    \qquad
    G(X) = N
  \]
  and
  \begin{equation}
    \label{equation:n-fold-tensor-isomorphism-h}
    G(\swi)     = \inj_0
    \qquad
    G(\swo)     = \inj_1
    \qquad
    G(\sxi)     = \bigtensor_{i \in N} \terminal_{\len{a(i)}} %
    \qquad
    G(\sxo)     = \bigtensor_{i \in N} \terminal_{\len{b(i)}} %
  \end{equation}
  \[
    G(\sporti)  = \left( \sum_{i \in N} \left(\inj_0 : \len{a(i)} \to \Ki\right) \right) \cp \inj
    \qquad
    \qquad
    G(\sporti)  = \left( \sum_{i \in N} \left(\inj_0 : \len{b(i)} \to \Ko\right) \right) \cp \inj
  \]
  \[
    G(\swn)     = \sum_{i \in N} (j \mapsto a(i)_j)
                + \sum_{i \in N} (j \mapsto b(i)_j) %
    \qquad
    \qquad
    G(\sxn)     = x
  \]
  where:
  \begin{itemize}
    \item $\Ki = \sum_{i \in N} \len{a(i)}$ is the total arity of all generators
    \item $\Ko = \sum_{i \in N} \len{b(i)}$ is the total coarity of all generators
    \item $\terminal_k : k \to \terminalObj$ is the unique terminal map.
    \item $\inj_0 : A \to A + B$ is the first injection of the coproduct in $\FinFun$
    \item $\inj$ denotes the canonical inclusion of a finite set into $\Nat$
  \end{itemize}
\end{proposition}
\begin{proof}
  Induction.
  The empty ($0$-fold) and singleton ($1$-fold) diagrams are already in the form
  above, so it remains to check the inductive step.

  Let $c = (\inj_0, \inj_1, G)$ be the $N$-fold tensor of singleton diagrams $x : N \to \Arr$.
  By inductive hypothesis, $c \cong \bigtensor_{i \in N} \singleton(x(i))$.
  Now let %
  $s = (\inj_0, \inj_1, H)$
  be a singleton diagram with the chosen typing $(a' : 1 \to \Obj^*, b' : 1 \to \Obj^*)$.
  It suffices to construct an isomorphism of structured cospans
  $\alpha : c \tensor s \leftrightarrow d$
  where $d = (\inj_0, \inj_1, J)$ is in the required form.

  Defining $\alpha$ is straightforward.
  Set $\alpha_W \defeq \gen{s-ex}$ with all other components as identities.
  Then one can compute that the image of $\swi$ and $\swo$ are the injections,
  and further that $\alpha$ is an isomorphism of structured cospans.
  For example,
  \[ \alpha_{\Eo} \cp (G(\swo) + H(\swo)) \cp \alpha_W
      = \gen{proof/n-fold-swo-0}
      = \gen{proof/n-fold-swo-1}
      = \inj_1
  \]
  with the other cases holding similarly.

  Finally, we must have by naturality that
  $J(\swn) = \alpha^{-1}_W \cp G(\swn) \tensor H(\swn)$,
  and so it remains to verify that this morphism is in the desired form.

  Observe that both $G(\swn)$ and $H(\swn)$ are coproducts of maps:
  by definition $H(\swn) = a' + b'$
  and $G(\swn) = g_a + g_b$, where
  $g_a = \sum_{i \in N} (j \mapsto a(i)_j)$
  and $g_b = \sum_{i \in N} (j \mapsto b(i)_j)$.
  for some $a, b : N \to \Obj^*$.
  Then calculate as follows using the associativity and commutativity axioms.
  \[
    G'(\swn) = \alpha^{-1}_W \cp (G(\swn) + H(\swn))
    \: = \: \gen{proof/n-fold-swn-0}
    \: = \: \gen{proof/n-fold-swn-1}
    \: = \: \gen{proof/n-fold-swn-2}
  \]
  We then have that
  $g_a + a' = \sum_{i \in N} (j \mapsto a(i)_j) + (j \mapsto a'(0)_j)$
  and
  $g_b + b' = \sum_{i \in N} (j \mapsto b(i)_j) + (j \mapsto b'(0)_j)$,
  and thus
  $J(\swn) = \sum_{i \in N+1} (j \mapsto (a + a')(i)_j) + \sum_{i \in N+1} (j \mapsto (b + b')(i)_j) $, as required.
\end{proof}

Computing the $N$-fold tensor product of operations can be computed in time
linear (sequential) and logarithmic (parallel) in the resulting diagram, as
witnessed by the following proposition.

\begin{proposition}
  \label{proposition:complexity-n-fold-tensor-isomorphism}
  Let $x : N \to \Arr$ be operations, and $d = (\inj_0, \inj_1, G)$ their
  $N$-fold tensor product as defined in
  \Cref{proposition:n-fold-tensor-isomorphism}.
  The sequential time complexity of computing $d$ is
  $O(G(W)) + O(G(\Ei)) + O(G(\Ei)) + O(G(X))$
  and the PRAM CREW time complexity is
  $O(\log G(X))$
\end{proposition}
\begin{proof}
  Computing $\inj_0$ and $\inj_1$ is $O(G(\Ei))$ and $O(G(\Eo))$ respectively
  in the sequential case, and constant in the parallel case.
  Then to compute $G$ in the sequential case is essentially just concatenation
  of arrays, and so the complexity is linear with respect to the size of each
  component, i.e.
  $O(G(W)) + O(G(\Ei)) + O(G(\Eo)) + O(G(X))$.
  Note that the $G(\sxn)$ component does not need to be computed, since it is
  given.

  In the parallel case, assuming the data for typings $a$ and $b$ is provided as
  a segmented array, one can compute
  components using the \texttt{repeat} and \texttt{segmented arange}
  functions, giving $O(\log N) = O(\log G(X))$ time complexity.
\end{proof}

\subsection{Composition of Diagrams}

Before giving the algorithm for composition of diagrams, we first illustrate it
with an example.
Composition in $\DiagramFrob$ is by pushout, as computed via coequalizer.
Concretely, when composing diagrams $d_0 \cp d_1$,
the basic idea is to first take their tensor product, and then coequalize the
wires corresponding to the shared boundary.

\begin{example}
  Let $d_0 = (s_0, t_0, G_0)$ and $d_1 = (s_1, t_1, G_1)$ be the diagrams
  illustrated below as cospans.
  \[ \gen{examples/string-diagram-composite} \]
  The process of composition is illustrated with the following steps.
  We first tensor $d_0$ and $d_1$ (left)
  and then identify the $\wNode$-nodes connected red arrow between them (right).
  \[ \gen{examples/composite-step1-tensor} \qquad \qquad \gen{examples/composite-step2-identify} \]
  Finally, we quotient together all connected components in the graph of $\wNode$-nodes,
  and take the source map of $d_0$, and the target of $d_1$.
  \[ \gen{examples/composite-step3-quotient} \]
  This amounts to coequalizing the maps $t_0 \cp \inj_0$ and $s_1 \cp \inj_1$.
\end{example}

Let us now describe this process formally.

\begin{proposition}[Composition of Diagrams]
  \label{proposition:composition}
  Let $d_0 \defeq (s_0, t_0, G_0) : A_0 \to A_1$
  and $d_1 \defeq (s_1, t_1, G_1) : A_1 \to A_2$
  be diagrams,
  and denote the coequalizer
  $q \defeq \coeq\left(\sgen{s-t0-inj0} \quad,\quad\sgen{s-s1-inj1}\right)$
  so that
  $q : G_0(W) + G_1(W) \to Q$.
  Then the composition $d_0 \cp d_1 = (s, t, G)$ is given by
  \begin{align*}
    s = \sgen{s-s0-inj0-q}
    \qquad \qquad
    t = \sgen{s-t1-inj1-q}
    \qquad \qquad
    G = \target(\alpha)
  \end{align*}
  where $\alpha : G_0 + G_1 \to G$ is the $\ACSet$ coequalizer whose $W$
  component is $q$, and all other components are $\id$,
  so that $G$ is the functor with the following data
  \begin{equation}
    \label{equation:composition-bpmg}
    \scalebox{0.8}{\input{composition-G.tikz}}
  \end{equation}
  and $G(\swn) = u$ is the canonical morphism of the coequalizer $q$, i.e. the
  unique function defined so that
  $u(q(i)) = (G_0(\swn) + G_1(\swn))(i)$.
\end{proposition}
\begin{proof}
  Let $\cospan{\SCL(A_0)}{\sigma_0}{G_0}{\tau_0}{\SCL(A_1)}$
  and $\cospan{\SCL(A_1)}{\sigma_1}{G_1}{\tau_1}{\SCL(A_2)}$
  be the respective structured cospans corresponding to $d_0$ and $d_1$ by
  \Cref{proposition:representation}.
  Composition of cospans $d \defeq d_0 \cp d_1$ is given by pushout, which we
  compute by coequalizer $q = \coeq(\tau_0 \cp \inj_0, \sigma_1 \cp \inj_1)$ as below.
  \[ \input{sc-pushout-by-coequalizer.tikz} \]
  Note that this makes $G = G_0 +_{\SCL(A_1)} G_1$.
  Recall that composition and colimits in $\ACSet$s are
  pointwise~\cite[Corollary 4]{acsets}, so we can calculate
  $
    \sigma_W
      = (\sigma_0 \cp \inj_0 \cp \alpha)_W
      = s_0 \cp \inj_0 \cp \coeq(t_0 \cp \inj_0)
  $.
  Moreover, since $\SCL$ is left adjoint, we can factor
  $\sigma = \SCL(s) \cp \iW_G$,
  and since $\iW_{G_W} = \id$, we must have
  $s = s_0 \cp \inj_0 \cp \coeq(t_0, s_1)$.
  By a similar argument, we have that
  $\tau = \SCL(t) \cp \iW_G$
  with $t = t_1 \cp \inj_1 \cp \coeq(t_0, s_1)$.
  We also see that $\alpha$ is indeed the coequalizer:
  \begin{align*}
    \alpha_W = \coeq(\tau_0 \cp \inj_0, \sigma_1 \cp \inj_1)_W & = \coeq(t_0 \cp \inj_0, s_1 \cp \inj_1) = q \\
    \alpha_Y = \coeq(\tau_0 \cp \inj_0, \sigma_1 \cp \inj_1)_Y & = \id
  \end{align*}
  The former holds because colimits of $\ACSet$s are pointwise,
  and the latter because for each $Y \neq W$,
  both $(\tau_0 \cp \inj_0)_Y$ and $(\sigma_1 \cp \inj_1)_Y$
  must equal the initial map by uniqueness, and the coequalizer of initial
  maps is the identity.

  Finally, we must verify that $u$ is indeed the canonical
  morphism of the coequalizer (and is thus well-defined).
  For this to hold, we must have that $q' = G_0(\swn) + G_1(\swn)$ defines a co-fork:
  \[ \input{s-cofork-t0.tikz} \quad = \quad \input{s-cofork-s1.tikz} \]
  This equation holds by the counit law and because $s$ and $t$ are the $W$
  components of label-preserving natural transformations $\sigma$ and $\tau$,
  and so $t_0 \cp G_0(\swn) = A_1(\swn) = s_1 \cp G_1(\swn)$.
  It then follows by the universal property that there exists a unique morphism
  $u : G(W) \to \Obj$
  such that
  $q \cp u = G_0(\swn) + G_1(\swn)$.
  Since $u$ is unique, we are therefore justified in defining it as
  the function $u(q(i)) \defeq q'(i)$.
  More directly, if $q(i) = q(j)$, then by the existence of a unique $u$ we must
  have $q'(i) = q'(j)$, and so $u$ is well-defined.
\end{proof}

\begin{remark}
  The above proposition says that in the case of composition, the coequalizer
  (and therefore pushout) of diagrams can be computed purely by considering
  their wirings.
  That is, there is always a choice of coequalizer which only quotients wire
  nodes of the bipartite graph, leaving other data is unchanged.
\end{remark}

Code for computing the composite of diagrams is straightforward,
\begin{minted}{python}
def compose(f, g):
    assert f.type[1] == g.type[0]
    h = f @ g
    q = f.t.inject0(g.G.W).coequalizer(g.s.inject1(f.G.W))

    return Diagram(
      s = f.s.inject0(g.G.W) >> q,
      t = g.t.inject1(f.G.W) >> q,
      G = h.G.coequalize_wires(q))
\end{minted}

We omit the implementation of \texttt{coequalize\_wires} above,
which is computed as in \eqref{equation:composition-bpmg}.
However, we give an explicit, data-parallel implementation for computing the
universal isomorphism $u$ now.

\begin{minted}{python}
  def universal(q: FiniteFunction, f: FiniteFunction):
    u = zeros(q.target)    # preallocate the table of u
    u[q.table] = f.table   # make the table of u such that q ; u = f
    return FiniteFunction(f.target, u)
\end{minted}

Note that in the above listing, we assume the PRAM CRCW model to achieve $O(1)$
parallel time complexity.
Specifically, the assignment
\texttt{u[q.table] = f.table} will allow for multiple conflicting writes to the
same memory location, with an arbitrary write succeeding.
However, note that since both $q$ and $f$ must be label-preserving maps,
all `conflicts' \emph{must} write the same value, ensuring correctness of the
implementation.

As expected, complexity of composition is linear (sequential) and logarithmic (parallel) in the size of the diagram.
\begin{proposition}
  \label{proposition:complexity-compose}
  Let $d_i = \cospan{\SCL(A_i)}{s_i}{G_i}{t_i}{\SCL(A_{i+1})}$
  be diagrams for $i \in \{0, 1\}$.
  Let $(s, t, G) = d_0 \tensor d_1$,
  then computing the composite $d_0 \cp d_1$ has
  sequential time complexity
  \[ O(A_0(W)) + O(A_1(W)) + O(A_2(W))
    + O(G(W))
    + O(G(\Ei))
    + O(G(\Eo))
    + O(G(X))
  \]
  and PRAM CRCW time complexity $O(\log A_1(W) + G(W))$
\end{proposition}
\begin{proof}
  Time complexity is that of tensor product, plus the additional cost of
  computing coequalizers and applying them to $G$.
  In the sequential case, this additional cost
  is $O(A_1) + O(G(W))$.
  In the PRAM CRCW case, it is
  $O(\log A_1(W)) + O(\log G(W))$.
\end{proof}

Finally, we verify that the tensor product and composite of diagrams with
well-formed apexes are also well-formed.
This ensures that $\DiagramFrob$ indeed forms a (sub)category.

\begin{corollary}
  \label{proposition:composition-preserves-well-formedness}
  If diagrams $d_0$ and $d_1$ have well-formed apexes,
  then so do tensor products $d_0 \tensor d_1$ and composites $d_0 \cp d_1$.
\end{corollary}
\begin{proof}
  The case of tensor products is immediate since the edge data of both $d_0$ and
  $d_1$ is contained in $d_0 \tensor d_1$.
  In the case of composites $d_0 \cp d_1$, observe that the coequalizer
  $\alpha : G_0 + G_1 \to G$ has all components except $W$ equal to the
  identity.
  Thus by \Cref{proposition:bpmg-morphisms-preserve-well-formed}, we can
  conclude that $G$ is well-formed as well.
\end{proof}

%% file: isomorphism.tex
We now construct an explicit isomorphism $\DiagramFrob \cong \FreeFrob$.
This will justify our claim that the arrows of $\DiagramFrob$
are string diagrams over the signature $\Signature + \Frob$.
Concretely, it will guarantee that equivalence classes of
$(\Signature+\Frob)$-terms are the same as equivalence classes of diagrams.
With this isomorphism, we can then think of diagrams as a structure for
representing syntax \emph{combinatorially}.
In \Cref{section:symmetric-monoidal}, we will extend this to an isomorphism
$\DiagramSMC \cong \FreeSMC$ by translating to our setting the notion of
\emph{monogamicity} originally due to \cite{rmsms}.

Our approach to constructing the isomorphism is as follows:

\begin{itemize}
  \item Define the functor $\ToDiagramSlow$ mapping $\Signature$-terms to diagrams.
  \item Define the `Frobenius Decomposition' of a diagram
  \item Define a functor $\FromDiagram$ using the Frobenius decomposition mapping diagrams to $\Signature$-terms
  \item Prove that $\ToDiagramSlow$ and $\FromDiagram$ form an isomorphism
\end{itemize}

Of particular note are the `Frobenius decompositions' which we will shortly
define.
Such decompositions are not only useful for constructing the isomorhpism, but
also in applying \emph{functors} to diagrams.
We describe an algorithm based on this decomposition in
\Cref{section:applying-functors}.

Mapping $\Signature$-terms to diagrams is straightforward.
We simply map generators $\Arr$ to their corresponding singleton diagrams
(\Cref{definition:diagram-singleton}),
and build the diagram inductively.
More formally, we can construct the following functor.

\begin{definition}[$\ToDiagramSlow$]
  \label{definition:to-diagram-slow}
  $\ToDiagramSlow : \FreeFrob \to \DiagramFrob$ is the
  identity-on-objects strict symmetric monoidal hypergraph functor defined
  inductively on operations
  $g \in \Arr$ as $\singleton(g)$.
\end{definition}

Notice that $\ToDiagramSlow$ has poor computational complexity.
Although individual operations of tensor and composition are linear in the size
of the resulting diagram, the inductively-defined functor above has $O(n^2)$ complexity.
This is because we repeatedly `append' a small diagram to an accumulator, which
is analogous to constructing a length-$n$ array by repeated appending.
We will fix this problem in \Cref{section:fast-functor}.
Before that, we first construct an inverse to the functor $\ToDiagramSlow$; a
functor mapping from diagrams to $\Signature$-terms.

\subsection{Diagrams to $\Signature$-terms}

We now define the functor $\FromDiagram : \DiagramFrob \to \FreeFrob$, mapping
diagrams to $\Signature$-terms.
To do so, we first show how each diagram can be factorised into what we call a
`Frobenius Decomposition'.
This composition relies on the \emph{hypergraph category} structure of
$\DiagramFrob$ as witnessed by
\cite[Theorem 3.12]{structured-cospans}
via \Cref{proposition:diagram-is-hypergraph-category}.

The main idea of the Frobenius decomposition is to separate the structure of
wires and generators in a string diagram.
To illustrate this, we begin with an example.

\begin{example}
  Suppose we have a morphism $f \cp g$ in a hypergraph category.
  \[ \input{s-f-cp-g.tikz} \]
  A \deftext{frobenius decomposition} of this morphism is a diagram like the one
  pictured below.
  \[ \input{s-f-cp-g-frob-decomp.tikz} \]
  The idea is to picture all the wires in the diagram at the top as a `bus'
  whose width is the number of internal wires in the diagram.
  Wires can then be connected to operations $g$ by the maps $e_s$ and $e_t$,
  which respectively map wires to inputs and outputs of operations.
  Similarly, the $s$ and $t$ maps specify which of the `internal' wires
  appear on the left and right boundaries, respectively.

  Note however that this decomposition is not unique.
  For example, we may permute the order of generators in the center to obtain an
  isomorphic diagram.
  \[ \input{s-f-cp-g-frob-decomp-alt.tikz} \]
\end{example}

\begin{remark}
  The purpose of Frobenius decompositions is to separate the elements of a diagram
  into frobenius spiders and operations.
  This makes it easier to define a hypergraph functor,
  whose action on Frobenius spiders is fixed, and so it will suffice to define
  its action on the $g$ part of the decomposition.
\end{remark}

In fact, every morphism in an arbitrary Hypergraph Category has a Frobenius
decomposition, which we state formally in the following proposition.

\begin{proposition}[Frobenius Decomposition]
  \label{proposition:frobenius-decomposition-existence}
  Let $f : A \to B$ be a morphism in a hypergraph category.
  Then there is a factorisation of $f$ called the \deftext{Frobenius
  Decomposition} of $f$ with the following form.
  \[ \input{frobenius-decomposition.tikz} \]
  where $\sigma$, $\tau$, $e_s$, and $e_t$ are Frobenius spiders,
  and $g$ is an $n$-fold tensoring of generators.
\end{proposition}
\begin{proof}
  Induction.
  See Appendix \ref{section:proofs-isomorphism-from-diagram} for the full proof.
\end{proof}

Since $\DiagramFrob$ is a hypergraph category
by \cite[Theorem 3.12]{structured-cospans},
every morphism must have a Frobenius decomposition.
In order to define the functor $\FromDiagram$, it will be useful to give a
\emph{specific} Frobenius decomposition in terms of which the functor will be
defined.
This will also serve to show that $\DiagramFrob$ is \emph{presented} by the
generators of $\Signature + \Frob$.

\begin{proposition}
  \label{proposition:diagram-frobenius-decomposition}
  Let
  $d = \cospan{\SCL(A)}{\SCL(s)\cp\iW_G}{G}{\SCL(t)\cp\iW_G}{\SCL(B)}$
  be a diagram over the signature
  $\Signature$.
  Then the following diagram is a Frobenius decomposition of $d$
  \begin{equation}
    \label{equation:diagram-frobenius-decomposition}
    \input{diagram-frobenius-decomposition.tikz}
  \end{equation}
  where:
  \begin{itemize}
    \item $W   = \MkWiresObj(\swn)$,
          $\Ei = \MkWiresObj(G(\swi \cp \swn))$, and
          $\Eo = \MkWiresObj(G(\swo \cp \swn))$
    \item $f_i = \MkWiresArr(G(\swi), G(\swn))$ and $f_o = \MkWiresArr(G(\swo), G(\swn))$
    \item $p = (\sort_{\bl G(\sxi) , G(\sporti) \br}, \id, \SCL(\Ei))$
          and
          $q = (\sort_{\bl G(\sxo) , G(\sporto) \br}, \id, \SCL(\Eo))$
          (where $\sort$ is defined in \Cref{proposition:unique-sort})
    \item $g = \bigtensor_i^X \sxn(i)$
      is the diagram formed by the $X$-fold tensor product of singleton diagrams
  \end{itemize}
\end{proposition}
\begin{proof}
  We defer the full proof to \Cref{section:proofs-isomorphism-from-diagram}, but
  provide a sketch here.
  It is clear that the composite pictured in
  \eqref{equation:diagram-frobenius-decomposition} is of the correct form:
  $\HalfSpider(s)$ and $\HalfSpider(t)^\dagger$ correspond directly to $s$ and $t$ in
  \Cref{proposition:frobenius-decomposition-existence}, and the composites
  $\HalfSpider(G(\swi))^\dagger \cp p$ and $q \cp \HalfSpider(G(\swo))$ to $e_s^\dagger$ and $e_t$.
  It therefore remains to show that this decomposition is indeed equal to $(s, t, G)$.
  This follows because there is a choice of coequalizer such that
  $p \cp g \cp q$
  is isomorphic to a tensoring $g' = (\inj_0, \inj_1, G')$ of the generators in
  $G$.
  Then replacing $p \cp g \cp q$ with $g'$, one can again choose coequalizers
  so that the full composite is equal to $(s, t, G)$.
\end{proof}

\begin{remark}
  Note that the sorting permutations
  $\sort_{\bl \sxi , \sporti \br}$
  and
  $\sort_{\bl \sxo , \sporto \br}$
  (\Cref{proposition:unique-sort})
  are unique
  because their `sorting keys' are monomorphic,
  which follows because $G$ is well-formed.
  In addition, instead of actually computing the sorting permutations using a
  sorting algorithm, one can compute them through arithmetic in linear
  (sequential) and constant (parallel) time.
\end{remark}

Frobenius Decompositions allow us to define a functor from diagrams to
$\Signature$-terms.
We simply need to map each component of the decomposition to a representative
$\Signature$-term.

\begin{proposition}
  \label{proposition:from-diagram}
  There is a strict symmetric identity-on-objects monoidal hypergraph functor
  $\FromDiagram : \DiagramFrob \to \FreeFrob$
  which maps arrows to their frobenius decompositions.
\end{proposition}
\begin{proof}
  In \Cref{section:proofs-isomorphism-from-diagram}
  we show that
  $\FromDiagram$ is
  well-defined (\Cref{proposition:from-diagram-is-well-defined}),
  a functor (\Cref{proposition:from-diagram-is-functor}),
  and strict monoidal
  (\Cref{proposition:from-diagram-is-strict-monoidal}).
  Finally, it preserves the symmetry and Special Frobenius monoid structure by definition,
  since Frobenius spiders are mapped to the corresponding spiders in $\FreeFrob$.
\end{proof}

\subsection{Isomorphism}

We can now show that $\ToDiagramSlow$ and $\FromDiagram$ form an isomorphism.
We use the `Change of Basis' lemma from~\cite[B.1]{coc}, which gives conditions
under which two categories are isomorphic.
We begin by proving one of these conditions.

\begin{proposition}
  \label{proposition:diagrams-generated-by-signature-frob}
  Diagrams are generated by the operations of $\Signature + \Frob$.
  That is, every diagram can be constructed by tensor and composition of $\id$,
  $\twist$, and generators of $\Signature + \Frob$.
\end{proposition}
\begin{proof}
  By \Cref{proposition:diagram-frobenius-decomposition},
  every diagram has a Frobenius decomposition.
  Since each element of the decomposition is either a Frobenius spider or a
  tensoring of operations,
  the full diagram can be constructed by tensor and composition of generators.
\end{proof}

We can now give the main result.

\begin{proposition}
  \label{proposition:isomorphism}
  $\ToDiagramSlow$ and $\FromDiagram$ form an isomorphism.
\end{proposition}
\begin{proof}
  The change of basis lemma from~\cite[B.1]{coc} says that if the following
  conditions hold then $\ToDiagramSlow$ and $\FromDiagram$ form an isomorphism.
  \begin{itemize}
    \item $\ToDiagramSlow$ and $\FromDiagram$ are strict symmetric monoidal
      identity-on-objects hypergraph functors
      (\Cref{definition:to-diagram-slow}, \Cref{proposition:from-diagram-is-strict-monoidal})
    \item $\FreeFrob$ is generated by operations of $\Frob$ (by definition)
    \item $\DiagramFrob$ is generated by singleton diagrams of operations in $\Frob$
      (\Cref{proposition:diagrams-generated-by-signature-frob})
    \item $\FromDiagram(\ToDiagramSlow(g)) \cong g$ for $g \in \Arr$.
  \end{itemize}
  It therefore suffices to show that $\FromDiagram \circ \ToDiagramSlow$ is
  inverse on generators $g \in \Arr$.
  This is straightforward to derive: let $\singleton(g) : A \to B$ be an arbitrary
  generator, then the result follows by applying the unit and counit axioms of
  special frobenius monoids as below.
  \[ g \qquad = \qquad \input{proof/isomorphism-inverse-on-generators-lhs.tikz} \qedhere \]
\end{proof}

%% file: symmetric-monoidal.tex
We now examine the case of symmetric monoidal categories.
Specifically, we will conclude that there is a subcategory $\DiagramSMC$ of
$\DiagramFrob$ which is isomorphic to $\FreeSMC$.

Our approach is based on that of \cite{rmsms}, who characterise
string diagrams over $\Signature+\Frob$ in terms of hypergraphs.
By requiring an additional combinatorial condition--`monogamous acyclicity'--
on their hypergraphs, they recover string diagrams over $\Signature$.
That is, \emph{monogamous acyclic} hypergraphs correspond to string diagrams
\emph{without the additional hypergraph structure} of $\Frob$.

Since the datastructure we describe in \Cref{definition:diagram} is essentially
an \emph{encoding} of the hypergraphs in~\cite{rmsms}, we simply translate this
condition to our setting now.

\begin{definition}[Monogamicity \cite{rmsms}]
  \label{definition:monogamous}
  A diagram $d = \cospan{\SCL(A)}{s}{G}{t}{\SCL(B)}$ is \deftext{monogamous} when
  $s$ and $t$ are mono, and the following hold for all $\wNode$-nodes $i$.
  \begin{align*}
    \indegree(i) = \begin{cases}
      0 & \text{if } \exists j . s(j) = i \\
      1 & \text{otherwise}
    \end{cases}
    \qquad & \qquad
    \outdegree(i) = \begin{cases}
      0 & \text{if } \exists j . t(j) = i \\
      1 & \text{otherwise}
    \end{cases}
  \end{align*}
  where we write
  $\indegree(i) = |G(\swi)^{-1}[i]|$
  and
  $\outdegree(i) = |G(\swo)^{-1}[i]|$
\end{definition}

\begin{definition}[Acyclicity]
  \label{definition:acyclic}
  A bipartite multigraph is \emph{acyclic} when its underlying graph is acyclic.
\end{definition}

Let us now verify that these condition indeed characterise morphisms of $\FreeSMC$.
We begin with a lemma which says that the monogamous acyclic `Frobenius spiders'
are simply permutations.

\begin{proposition}
  \label{proposition:monogamous-acyclic-spiders}
  If $d = (s, t, \SCL(Y))$ is a monogamous acyclic Frobenius spider
  then $s$ and $t$ are permutations.
\end{proposition}
\begin{proof}
  By \Cref{definition:monogamous}, $s$ and $t$ are monomorphisms (i.e., injective).
  Then observe that $\SCL(Y)$ is a discrete graph, so each $\wNode$-node must have
  in- and out-degree $0$.
  For this to hold, it must be that each node is in the image of both $s$ and
  $t$, so they are surjective.
  It then follows immediately that $s$ and $t$ are bijections.
\end{proof}

Now check that those morphisms of $\FreeSMC$ map under $\ToDiagram$ to
monogamous acyclic diagrams.

\begin{proposition}
  \label{proposition:free-smc-to-diagram}
  If $f$ is a morphism of $\FreeSMC$, then $\ToDiagram(f)$ is monogamous acyclic.
\end{proposition}
\begin{proof}
  Induction.
  \begin{description}
    \item[Base Case: Generators]
      By \Cref{proposition:monogamous-acyclic-spiders} generators $\id$
      and $\twist$ are monogamous acyclic.
      It is then straightforward to verify that singleton diagrams
      (\Cref{definition:diagram-singleton}) are monogamous acyclic,
      and so we have for each $g \in \Arr$
      that $\ToDiagram(g)$ is monogamous acyclic.
    \item[Inductive Step: Tensor Product]
      Let $d_i = (s_i, t_i, G_i)$ be monogamous acyclic diagrams
      for $i \in \{0, 1\}$.
      Then
      $
        (s, t, G)
        \defeq d_0 \tensor d_1
        = (s_0 \tensor s_1, t_0 \tensor t_1, G_0 \tensor G_1)
      $
      by definition.
      $G$ is acyclic because it is the disjoint union of acyclic graphs,
      and $s$ and $t$ are mono because $s_i$ and $t_i$ are.
      It then remains to check the in/outdegree conditions.
      We verify the former and omit the latter for brevity.

      Suppose $i$ is a $\wNode$-node in $G_0$.
      Then if $\indegree(i) = 1$, by assumption $i$ is not in the image of $s_0$,
      and thus not in the image of $s_0 \tensor s_1$.
      Moreover, if $\indegree(i) = 0$, then there is some $j$ with $s_0(j) = i$,
      and then $(s_0 \tensor s_1)(j) = i$.

      Now suppose $i$ is a $\wNode$-node in $G_1$, and write $m \defeq G_0(W)$.
      If $\indegree(i) = 1$ then by assumption $i$ is not in the image of $s_1$,
      and $m + i$ is not in the image of $s_0 \tensor s_1$.
      Further, if $\indegree(i) = 0$, then there is some $j$  with $s_1(j) = i$,
      and then $(s_0 \tensor s_1)(m + j) = i$.
    \item[Inductive Step: Composition]
      Let
      $d_i = \cospan{\SCL(A_i)}{\SCL(s_i) \cp \iW_{G_i}}{G_i}{\SCL(t_i) \cp \iW_{G_i}}{\SCL(A_{i+1})}$
      be diagrams for $i \in \{0, 1\}$.
      Their composite is given in terms of the coequalizer
      $q = \coeq(t_0 \cp \inj_0, s_1 \cp \inj_1)$.
      Since $t_0 \cp \inj_0$ and $s_1 \cp \inj_1$ are mono,
      if $q(i) = q(j)$ for $i \neq j$,
      then there is some unique $a \in A_1$ with
      $(t_0 \cp \inj_0)(a) = i$ and
      $(s_1 \cp \inj_1)(a) = j$ (or vice-versa).
      But then $\outdegree(i) = 0$ and
      and $\indegree(j) = 0$ because $d_0$ and $d_1$ are monogamous acyclic,
      so $\indegree(q(i)) = \indegree(i)$ and $\outdegree(q(i)) = \outdegree(j)$,
      and the quotiented graph is monogamous as well.

      Finally, note that acyclicity is preserved for the same reason:
      if $G_0$ and $G_1$ are acyclic, then only $\wNode$-nodes in $G_1$ are
      reachable from $q(i)$, and thus no cycles are possible.
  \end{description}
\end{proof}

We now check the reverse: that all monogamous acyclic diagrams correspond to
morphisms of $\FreeSMC$.
This fact is first proven in \cite{rmsms}, and is thus somewhat expected.
Moreover, our argument is essentially the same as \cite[Proposition B.3]{coc},
so we only sketch a proof here.

\begin{proposition}
  \label{proposition:monogamous-acyclic-diagrams-smc}
  If $d$ is a monogamous acyclic diagram, then
  $\FromDiagram(d)$ is isomorphic to a morphism of $\FreeSMC$.
\end{proposition}
\begin{proof}
  Let $d = (s, t, G)$ be a diagram.
  Using the acyclicity property, $d$ can be decomposed into the form
  $s \cp \ToDiagramSlow\left(\sgen{s-layer-g}\right) \cp d'$.
  Then, since $s$ is mono, the result holds by induction on the number of
  generators $d$.
\end{proof}

By these two lemmas, there exists a category of monogamous acyclic diagrams
whose morphisms are isomorphic to those of the free symmetric monoidal category
over the same signature.

\begin{proposition}
  \label{proposition:diagram-smc}
  There is a subcategory $\DiagramSMC$ of $\DiagramFrob$ whose morphisms are
  \emph{monogamous acyclic} diagrams, and $\DiagramSMC \cong \FreeSMC$.
\end{proposition}
\begin{proof}
  The property of monogamous acyclicity is closed under tensor and composition
  as guaranteed by Propositions \ref{proposition:free-smc-to-diagram} and
  \ref{proposition:monogamous-acyclic-diagrams-smc}.
  Moreover, the isomorphism $\DiagramSMC \cong \FreeSMC$ is given by the
  functors $\ToDiagramSlow$ and $\FromDiagram$ restricted to monogamous acyclic
  diagrams.
\end{proof}

\begin{remark}
  Note that the functor $\FromDiagram$ does \emph{not} directly give a
  $\Signature$-term, but `only' a $(\Signature+\Frob)$-term equivalent to one.
  If one wishes to recover a bona fide $\Signature$-term from a diagram, it is
  straightforward to use a decomposition based on e.g., topological sort as
  described in \cite[Proposition B.13]{coc}.
\end{remark}

%% file: fast-functor.tex
In the previous section, we gave algorithms for tensor and composition of diagrams
which had $O(n)$ sequential and $O(\log n)$ PRAM complexity.
This allowed us to define a functor
$\ToDiagramSlow : \FreeFrob \to \DiagramFrob$.
However, the complexity of applying this functor naively is $O(n^2)$.

We now give a faster functor
$\ToDiagramFast : \FreeFrob \to \DiagramFrob$
which computes the resulting diagram directly.
Importantly, the algorithm presented here will have linear (sequential) and
logarithmic (PRAM CREW) time complexity in the size of the resulting diagram.
Given a $\Signature$-term $t$, the basic idea is to tensor all the generators of $t$
and then wire them all together in a single step.
We illustrate this process with an example.

\begin{example}
  Consider the $\Signature$-term below (rendered as a tree) and its corresponding string diagram.
  \[ \input{wiring/example-expression.tikz} \qquad = \qquad \input{wiring/example-string-diagram.tikz} \]
  The idea of this section is to build the string diagram on the right by first
  tensoring all the diagrams $d_0 \tensor d_1 \tensor d_2$,
  and then `wiring up' the diagram in a single step.
  More visually, we can picture this as the string diagram below, with the
  target of $d_0$ `wired up' to the source of $d_1$.
  \[ \input{wiring/example-string-diagram-wiring.tikz} \]
\end{example}

\subsection{A Faster Functor}
In order to `wire up a diagram in a single step', we must specify which
$\wNode$-nodes are to be quotiented together,
and which will belong on the boundaries.
This specification comes in the form of a $4$-tuple, which we call the `wiring
maps'.

\begin{definition}[Wiring Maps]
  \label{definition:wiring-maps}
  Let $e$ be a binary tree with $n$ leaves
  the diagrams
  $d_i = (s_i, t_i, G_i)$,
  and whose nodes are labeled either $\cp$ or $\tensor$.
  The four \deftext{wiring maps} of $e$ have types
  \begin{equation}
    e_s  : A_S \to G(W) \qquad
    e_t  : A_T \to G(W) \qquad
    e_s' : A_I \to G(W) \qquad
    e_t' : A_I \to G(W)
  \end{equation}
  where $G = \bigtensor\limits_{i \in \ord{n}} G_i(W)$.
  and $A_S$, $A_T$, and $A_I$ are computed recursively as follows.
  When $e$ is a binary tree consisting of a single leaf
  $e = \cospan{\SCL(A)}{\SCL(s) \cp \iW_G}{G}{\SCL(t) \cp \iW_G}{\SCL(B)}$,
  the wiring maps are the $\Wires$ morphisms given below.
  \[
    e_s \defeq \gen{wiring/gt-A-GW-s}
    \qquad
    e_t \defeq \gen{wiring/gt-B-GW-t}
    \qquad
    e_s' \defeq \gen{wiring/gt-GW-zero}
    \qquad
    e_t' \defeq \gen{wiring/gt-GW-zero}
  \]
  The wiring maps of the tree
  $e = \scalebox{0.7}{\begin{tikzpicture}[level distance=7mm, sibling distance=7mm, baseline=-4mm]
    \node {$\tensor$}
      child {node {$l$}}
      child {node {$r$}};
  \end{tikzpicture}}$
  are given by
  \begin{equation}
    e_s  \defeq \gen{wiring/ls-tensor-rs} \qquad
    e_t  \defeq \gen{wiring/lt-tensor-rt} \qquad
    e_s' \defeq \gen{wiring/ls-prime-tensor-rs-prime} \qquad
    e_t' \defeq \gen{wiring/lt-prime-tensor-rt-prime}
  \end{equation}
  and those for
  $e = \scalebox{0.7}{\begin{tikzpicture}[level distance=7mm, sibling distance=7mm, baseline=-4mm]
    \node {$\cp$}
      child {node {$l$}}
      child {node {$r$}};
  \end{tikzpicture}}$
  are defined
  \begin{equation}
    e_s  \defeq \gen{wiring/ls-inj0} \qquad %
    e_t  \defeq \gen{wiring/rt-inj1} \qquad
    e_s' \defeq \gen{wiring/ls-prime-tensor-rs-plus-rs-prime} \qquad
    e_t' \defeq \gen{wiring/lt-prime-plus-lt-tensor-rt-prime}
  \end{equation}
\end{definition}

The maps $e'_s$ and $e'_t$ together serve to identify the internal
$\wNode$-nodes to be quotiented together;
the final diagram is given by coequalizing these maps.
Meanwhile the maps $e_s$ and $e_t$ identify the $\wNode$-nodes corresponding to
the boundaries of the diagram.
To illustrate this, we now give the 

\begin{example}
  Recall the example $\Signature$-term and its corresponding string diagram
  below.
  \[ \input{wiring/example-expression.tikz} \qquad = \qquad \input{wiring/example-string-diagram.tikz} \]
  Writing $s_i$ and $t_i$ for the source and target maps of $d_i$,
  we have
  \[
    e_s = \gen{wiring/example-wiring-maps-es}
    \qquad
    e_t = \gen{wiring/example-wiring-maps-et}
    \qquad
    e'_s = \gen{wiring/example-wiring-maps-es-prime}
    \qquad
    e'_t = \gen{wiring/example-wiring-maps-et-prime}
  \]
  Notice that coequalizing `internal wiring maps' $e'_s$ and $e'_t$ will
  identify the targets of $d_0$ with the sources of $d_1$.
\end{example}

Using the wiring maps, we can now define the `fast' functor $\ToDiagramFast$.

\begin{definition}
  \label{definition:fast-functor}
  $\ToDiagramFast : \FreeFrob \to \DiagramFrob$
  is the identity-on-objects functor defined on arrows ($\Signature$-terms) as
  follows.
  Let $e$ be the binary tree formed by replacing the leaves of a $\Signature$-term
  with the corresponding singleton diagrams,
  and let
  \begin{equation*}
    s  : A_S \to G(W) \qquad
    t  : A_T \to G(W) \qquad
    s' : A_I \to G(W) \qquad
    t' : A_I \to G(W)
  \end{equation*}
  be the wiring maps of $e$.
  Define $q \defeq c(t', s')$ to be the coequalizer of internal wiring maps.
  \begin{equation}
    \label{equation:cd-f-of-tree}
    \input{cd-f-of-tree.tikz}
  \end{equation}
  As in \Cref{proposition:composition}, $q$ lifts to a coequalizer $\alpha$ in $\BPMG$
  whose $W$ component is $q$ and all other components are $\id$.
  $\ToDiagramFast(e)$ is then the diagram $(s \cp q, t \cp q, \target(\alpha))$.
\end{definition}

\begin{remark}
  Notice that not only does $\ToDiagramFast$ define a functor from $\FreeFrob$,
  but it can also be used to arbitrarily `wire up' existing diagrams,
  since it can be applied to any tree whose leaves are diagrams- not just singleton diagrams.
\end{remark}

In Proposition \ref{proposition:fast-isomorphic-to-slow} we will verify that
$\ToDiagramFast$ and $\ToDiagramSlow$
give isomorphic diagrams when applied to the same expression tree.
This will also serve to verify the well-definedness of Definition
\ref{definition:fast-functor}.
A key lemma in this proof will be the functoriality of $\ToDiagramFast$, which
we now prove.

\begin{proposition}[Functoriality of $\ToDiagramFast$]
  \label{proposition:fast-functor-preserves-composition}
  $\ToDiagramFast(\id) \cong \id$
  and
  $
    \ToDiagramFast\left(\scalebox{0.7}{\begin{tikzpicture}[level distance=7mm, sibling distance=7mm, baseline=-4mm]
      \node {$\cp$}
        child {node {$l$}}
        child {node {$r$}};
    \end{tikzpicture}}\right)
    \cong \ToDiagramFast(l) \cp \ToDiagramFast(r) \\
  $
\end{proposition}
\begin{proof}
  We have $\ToDiagramFast(\id) = \id$ by definition,
  so it remains to show that
  $
    \ToDiagramFast\left(\scalebox{0.7}{\begin{tikzpicture}[level distance=7mm, sibling distance=7mm, baseline=-4mm]
      \node {$\cp$}
        child {node {$l$}}
        child {node {$r$}};
    \end{tikzpicture}}\right)
    \cong \ToDiagramFast(l) \cp \ToDiagramFast(r)
  $.
  This means checking that $Q$ in \eqref{equation:cd-f-of-tree}
  is isomorphic to $Q'$ in \eqref{equation:cd-f-composed},
  and that this isomorphism extends to an isomorphism of cospans.
  \begin{equation}
    \label{equation:cd-f-composed}
    \input{cd-f-composed.tikz}
  \end{equation}
  The proof is in four parts.
  (1) Construction of the `forward' map $u : Q \to Q'$,
  (2) the `reverse' map $u' : Q' \to Q$,
  (3) showing they form an isomorphism,
  and (4) verifying this extends to a cospan isomorphism.
  \begin{description}
    \item[(1) Forward]
      Let us first construct the unique morphism $u : Q \to Q'$ by showing
      the existence of a map $x : G(W) \to Q'$ which coequalizes $s'$ and $t'$.

      By definition, $G(W) = G_l(W) \otimes G_r(W)$,
      and so we can define $x \: \defeq \: \gen{wiring/def-x}$.
      Now we will verify that $s' \cp x = t' \cp x$.
      First, recall that by \Cref{definition:wiring-maps}, we have
      \[
        s' \quad = \quad \gen{wiring/ls-prime-tensor-rs-plus-rs-prime} \qquad
        t' \quad = \quad \gen{wiring/lt-prime-plus-lt-tensor-rt-prime}
      \]
      and so we can calculate
      \begin{align*}
        t' \cp x \qquad
          & = \qquad \gen{wiring/proof-functor-fwd-lhs} \\ \\
          & = \qquad \gen{wiring/proof-functor-fwd-0} \\ \\
          & = \qquad \gen{wiring/proof-functor-fwd-1} \\ \\
          & = \qquad \gen{wiring/proof-functor-fwd-2} \\ \\
          & = \qquad \gen{wiring/proof-functor-fwd-3} \\ \\
          & = \qquad \gen{wiring/proof-functor-fwd-rhs} \qquad = \qquad s' \cp x
      \end{align*}
      Thus since $x$ coequalizes $t'$ and $s'$, by the universal property of
      coequalizers there is a unique morphism $u : Q \to Q'$
      such that $q \cp u = x$.
    \item[(2) Reverse]
      We now construct the reverse morphism $u' : Q' \to Q$.
      First define maps $y_l : G_l(W) \to Q$ and $y_r : G_r(W) \to Q$ as follows.
      \[
        y_l \quad \defeq \quad \gen{wiring/g-yl}
        \qquad \qquad
        y_r \quad \defeq \quad \gen{wiring/g-yr}
      \]
      Now calculate that $y_l$ coequalizes $l'_t$ and $l'_s$.
      \[
        l_t' \cp y_l
        \: = \:
        \gen{wiring/proof-functor-rev-yl-lhs}
        \: = \:
        \gen{wiring/proof-functor-rev-yl-0}
        \: = \:
        \gen{wiring/proof-functor-rev-yl-1}
        \: = \:
        \gen{wiring/proof-functor-rev-yl-rhs}
        \: = \:
        l_s' \cp y_l
      \]
      A similar calculation shows that $r_t' \cp y_r = r'_s \cp y_r$,
      and so by the universal property there must exist unique morphisms
      $u_l : Q_l \to Q$ and $u_r : Q_r \to Q$ such that
      $y_l = q_l \cp u_l$ and $y_r = q_r \cp u_r$.

      Using $u_l$ and $u_r$, we can define
      $y \defeq \gen{wiring/def-y}$
      and calculate that
      \begin{align*}
        l_t \cp q_l \cp \inj_0 \cp y \qquad
          & = \qquad \gen{wiring/proof-functor-rev-y-lhs} \\ \\
          & = \qquad \gen{wiring/proof-functor-rev-y-0} \\ \\
          & = \qquad \gen{wiring/proof-functor-rev-y-1} \\ \\
          & = \qquad \gen{wiring/proof-functor-rev-y-2} \\ \\
          & = \qquad \gen{wiring/proof-functor-rev-y-3} \\ \\
          & = \qquad \gen{wiring/proof-functor-rev-y-rhs} \qquad = \qquad r_s \cp q_r \cp \inj_1 \cp y
      \end{align*}
      Thus since $y$ coequalizes
      $l_t \cp q_l \cp \inj_0$
      and
      $r_s \cp q_r \cp \inj_1$,
      there must exist a unique $u' : Q' \to Q$
      such that $y = q' \cp u'$.

    \item[(3) Isomorphism]
      We have constructed morphisms $u : Q \to Q'$ and $u' : Q' \to Q$, so it
      remains to verify these indeed form an isomorphism.
      Let us first verify $u \cp u' = \id$.
      By definition, $t' \cp q = s' \cp q$.
      Therefore, by the universal property, there is a unique morphism $v$ such that
      $q \cp v = q$.
      Clearly this holds for $v = \id$,
      so if $q \cp u \cp u' = q$, then $u \cp u' = \id$ by uniqueness.
      This holds as follows.
      \[
        \gen{wiring/proof-functor-iso-fwd-lhs}
          = \gen{wiring/proof-functor-iso-fwd-0}
          = \gen{wiring/proof-functor-iso-fwd-1}
          = \gen{wiring/proof-functor-iso-fwd-2}
          = \gen{wiring/proof-functor-iso-fwd-rhs}
      \]

      Similarly, if we can show
      $q' \cp u' \cp u = q'$,
      then $u' \cp u = \id$ by uniqueness.
      \[
        \gen{wiring/proof-functor-iso-rev-lhs}
          = \gen{wiring/proof-functor-iso-rev-0}
          = \gen{wiring/proof-functor-iso-rev-1}
          = \gen{wiring/proof-functor-iso-rev-2}
          = \gen{wiring/proof-functor-iso-rev-rhs}
      \]
      Note that in the third step, we use that $u_l \cp u = \gen{wiring/s-inj0-qprime}$
      and $u_r \cp u = \gen{wiring/s-inj1-qprime}$.
      These follow again by the universal property; we prove the former as
      follows.
      Recall that $q_l$ is the coequalizer of $l'_t$ and $l'_s$ by definition,
      and so clearly $\gen{wiring/proof-functor-iso-rev-yl-rhs}$ must coequalize
      them as well.
      By the universal property, there is a unique morphism $v$ such that
      $q_l \cp v = \gen{wiring/proof-functor-iso-rev-yl-rhs}$.
      It is clear that $v = \gen{wiring/s-inj0-qprime}$ satisfies this equality,
      and so if
      $q_l \cp u_l \cp u = \gen{wiring/proof-functor-iso-rev-yl-rhs}$
      then $u_l \cp u = \gen{wiring/s-inj0-qprime}$ by uniqueness.
      Calculemus:
      \[
        \gen{wiring/proof-functor-iso-rev-yl-lhs}
        = \gen{wiring/proof-functor-iso-rev-yl-0}
        = \gen{wiring/proof-functor-iso-rev-yl-1}
        = \gen{wiring/proof-functor-iso-rev-yl-rhs}
      \]
      And so $u_l \cp u = \gen{wiring/s-inj0-qprime}$ as desired.
      A similar argument yields $u_r \cp u = \gen{wiring/g-yr}$,
      completing the proof that $u$ and $u'$ form an isomorphism.

    \item[(4) Cospan Isomorphism]
      We now show that $u$ extends to an isomorphism of cospans.
      The following calculation equates the cospan source legs
      of \eqref{equation:cd-f-of-tree} and \eqref{equation:cd-f-composed}
      via the universal map $u$.
      \[
        \gen{wiring/proof-functor-cospan-lhs}
        = \gen{wiring/proof-functor-cospan-0}
        = \gen{wiring/proof-functor-cospan-rhs}
      \]
      A similar calculation shows the reverse direction
      and symmetric arguments apply to the target maps,
      completing the proof.
  \end{description}
\end{proof}

We now prove the isomorphism between of $\ToDiagramFast$ and $\ToDiagramSlow$.

\begin{proposition}
  \label{proposition:fast-isomorphic-to-slow}
  Let $e$ be a $\Signature$-term.
  Then there is an isomorphism of diagrams $\ToDiagramFast(e) \cong \ToDiagramSlow(e)$.
\end{proposition}

\begin{proof}
  We proceed by induction.
  In the base case when $e \in \{ \id, \twist \} \cup \Arr$
  we have $\ToDiagramFast(e) = \ToDiagramSlow(e)$.
  For the inductive step, assume that there are isomorphisms
  $u_l : \ToDiagramFast(l) \cong \ToDiagramSlow(l)$
  and
  $u_r : \ToDiagramFast(r) \cong \ToDiagramSlow(r)$.
  We must now show that
  $\ToDiagramFast(e) \cong \ToDiagramSlow(e)$ for
  the cases
  $e =
  \scalebox{0.7}{\begin{tikzpicture}[level distance=7mm, sibling distance=7mm, baseline=-4mm]
    \node {$\tensor$}
      child {node {$l$}}
      child {node {$r$}};
  \end{tikzpicture}}$
  and
  $e =
  \scalebox{0.7}{\begin{tikzpicture}[level distance=7mm, sibling distance=7mm, baseline=-4mm]
    \node {$\cp$}
      child {node {$l$}}
      child {node {$r$}};
  \end{tikzpicture}}$.
  In the former case, observe that
  \[
    \ToDiagramFast\left(\scalebox{0.7}{\begin{tikzpicture}[level distance=7mm, sibling distance=7mm, baseline=-4mm]
      \node {$\tensor$}
        child {node {$l$}}
        child {node {$r$}};
    \end{tikzpicture}}\right)
    \: =     \: \ToDiagramFast(l) \tensor \ToDiagramFast(r) \\
    \: \cong \: \ToDiagramSlow(l) \tensor \ToDiagramSlow(r) \\
    \: =     \: \ToDiagramSlow(l \tensor r)
  \]
  where the isomorphism $u_l \tensor u_r$ is given by inductive hypothesis.
  In the latter case, we derive as follows.
  \[
    \ToDiagramFast\left(\scalebox{0.7}{\begin{tikzpicture}[level distance=7mm, sibling distance=7mm, baseline=-4mm]
      \node {$\cp$}
        child {node {$l$}}
        child {node {$r$}};
    \end{tikzpicture}}\right)
    \: \cong \: \ToDiagramFast(l) \cp \ToDiagramFast(r) \\
    \: \cong \: \ToDiagramSlow(l) \cp \ToDiagramSlow(r) \\
    \: =     \: \ToDiagramSlow(l \cp r)
  \]
  The first two steps of this derivation use
  functoriality of
  $\ToDiagramFast$ (\Cref{proposition:fast-functor-preserves-composition})
  and well-definedness of diagram composition~\cite[Proposition 3.11]{structured-cospans}.
  respectively.

\end{proof}

\begin{corollary}
  \label{proposition:complexity-fast-functor}
  Let $e$ be a $\Signature$-term (binary tree) of
  $N$ leaves.
  Computing $\ToDiagramFast(e)$ has
  $O(N)$ (sequential) and $O(\log N)$ (PRAM CRCW) time complexity.
\end{corollary}
\begin{proof}
  Since $\Signature$ is fixed, we regard the number of wires and edges in the
  diagram to be proportional to $N$
  (cf. \cite[Proposition 6.1, Bounded Sparsity]{coc}).
  The $N$-fold tensor product of leaves has the given complexity by
  \Cref{proposition:complexity-n-fold-tensor-isomorphism},
  and then it remains to compute connected components, which is again $O(N)$
  sequential and $O(\log N)$ parallel.
\end{proof}

%% file: fast-wiring.tex
In \Cref{section:fast-functor} we showed how a diagram can be `wired up' all at
once using the coequalizer and wiring maps (\Cref{definition:wiring-maps}).
We now give a parallel algorithm for computing the wiring maps.
The algorithm is in two parts.
The first part (\Cref{section:fast-wiring-reduce}) reduces the problem of
computing wiring maps to that of finding the `ancestor maps'.
Ancestor maps find, for all nodes in a given $\Signature$-term, a closest
ancestor (if any) labeled $\cp$, for which a given node is
in a specific (left or right) subtree.
The second part of the algorithm (\Cref{section:fast-wiring-predicate}) is a
data-parallel method for computing these ancestor maps efficiently.

\subsection{Reducing the Problem to a Tree Predicate}
\label{section:fast-wiring-reduce}

We now give an alternative, `parallel-friendly' definition of the wiring maps,
which is proven equivalent in
\Cref{proposition:wiring-maps-alternative}.
Efficiently computing the wiring maps in parallel via this alternative
definition will depend on computing the `\emph{ancestor maps}', defined below.

\begin{definition}[Left and Right $\cp$-Ancestors]
  In a binary tree $e$,
  a node $i$ has a \deftext{right $\cp$-ancestor}
  if it is in the \emph{left} subtree of any node $j$ labeled $\cp$.
  Respectively, $i$ has a \deftext{left $\cp$-ancestor} if it is in the
  \emph{right} subtree of any node $j$ labeled $\cp$.
\end{definition}

Computing the left/right ancestors of a node determines how to `wire up' the diagram.
We illustrate this with the following example.

\begin{example}
  Pictured below is a $\Signature$-expression rendered as a tree,
  with the closest left (resp. right) $\cp$-ancestor of each node listed below,
  with $\cdot$ indicating a node has no left (resp. right) $\cp$-ancestor.
  \[ \input{wiring/example-capture-tree.tikz} \]
  Diagrams (leaves) $d_0$ and $d_1$ are `bound' by node $1$, indicating that the
  outputs of $d_0$ will be `connected' to the inputs of $d_1$.
\end{example}

\begin{definition}[Ancestor Maps]
  \label{definition:ancestor-maps}
  Let $e$ be a binary tree of $n$ leaves and $n - 1$ (non-leaf) nodes.
  The \deftext{left and right ancestor maps} of $e$
  are the following functions.
  \begin{align*}
    a_L & : n \to n \\
    a_L & \defeq \begin{cases}
      0   & \qquad \text{if $i$ has no left $\cp$-ancestors} \\
      j+1 & \qquad \text{for the largest $j$ such that $j$ is a left $\cp$-ancestor of $i$}
    \end{cases}
  \end{align*}
  \begin{align*}
    a_R & : n \to n \\
    a_R & \defeq \begin{cases}
      j & \qquad \text{for the smallest $j$ such that $j$ is a right $\cp$-ancestor of $i$} \\
      n   & \qquad \text{if $i$ has no left $\cp$-ancestors}
    \end{cases}
  \end{align*}
\end{definition}

\begin{remark}
  The left ancestor map $a_L$ returns the index of the `closest' ancestor for
  which a given node is in the right subtree- meaning that the ancestor is
  `left' of the node in an inorder traversal.
\end{remark}

We will show how to efficiently compute the ancestor maps in
\Cref{section:fast-wiring-predicate}.

\begin{proposition}[Wiring Maps Alternative Definition]
  \label{proposition:wiring-maps-alternative}
  Let $e$ be a tree with $n$ leaves the diagrams $d_i : A_i \to B_i$ for $i \in \ord{n}$.
  The wiring maps of $e$ can be equivalently defined as follows.
  \begin{align}
    e_s
    \quad = \quad \inj^{S,S'}_0 \cp p_s \cp \bigtensor\limits_{i \in \ord{n}} s_i
    \quad = \quad \gen{wiring/def-es}
    \nonumber
    \\
    \nonumber
    \\
    e_s'
    \quad = \quad \inj^{S,S'}_1 \cp p_s \cp \bigtensor\limits_{i \in \ord{n}} s_i
    \quad = \quad \gen{wiring/def-es-prime}
    \nonumber
    \\
    \label{equation:wiring-maps-alternative}
    \\
    e_t
    \quad = \quad \inj^{T',T}_1 \cp p_t \cp \bigtensor\limits_{i \in \ord{n}} t_i
    \quad = \quad \gen{wiring/def-et}
    \nonumber
    \\
    \nonumber
    \\
    e_t'
    \quad = \quad \inj^{T',T}_0 \cp p_t \cp \bigtensor\limits_{i \in \ord{n}} t_i
    \quad = \quad \gen{wiring/def-et-prime}
    \nonumber
  \end{align}
  Where
  $p_s$ and $p_t$ are the natural transformations defined by
  the stable sorting permutations
  $\sort_{\tl a_L, \id \tr}$ and
  $\sort_{\tl a_R, \id \tr}$, respectively\footnote{
    We define `stable' sorts in \Cref{definition:stable-sort},
    and omit subscripts of these natural transformations to reduce notational burden.
  },
  and the objects $S$, $S'$, $T'$, and $T$ are defined as follows.
  \[
    S  = \bigtensor\limits_{\{i \in \ord{n} | a_L(i) = 0  \}} A_i
    \qquad \qquad
    S' = \bigtensor\limits_{\{i \in \ord{n} | a_L(i) > 0  \}} A_i
  \]
  \[
    T  = \bigtensor\limits_{\{i \in \ord{n} | a_R(i) = n-1\}} B_i
    \qquad \qquad
    T' = \bigtensor\limits_{\{i \in \ord{n} | a_R(i) < n-1\}} B_i
  \]
\end{proposition}
\begin{remark}
  The objects $S$ and $S'$ refer to \emph{unbound} and \emph{bound} source maps, respectively.
  The source map of a diagram is `unbound' if it won't be quotiented with other
  wires in construction of the diagram.
  Notice that the objects in the tensor $S$ are precisely those source objects
  for which there is no left $\cp$-ancestor.
  Thus, these source maps will constitute the left dangling wires of the resulting diagram.
  Meanwhile, $S'$ represent the source maps which \emph{do} have a left
  $\cp$-ancestor, and are therefore part of a composition- these will indicate
  the $\wNode$-nodes in the graph which need to be quotiented in the result.
  Similarly, $T$ and $T'$ are the unbound and bound \emph{target} maps, respectively.
\end{remark}
\begin{proof}
  Induction.
  \begin{description}
    \item[Base Case: Single Leaf]
      Let $e$ be a tree consisting of a single leaf,
      a diagram $d = (s, t, G)$ of type $A \to B$.
      Then $a_L = \id$, so $p_s = \id$,
      and $a_R = \id$ so $p_t = \id$.
      Moreover, $S = A$, $S' = \unit$, $T = B$, and $T' = \unit$,
      and one can calculate as follows.
      \begin{align*}
        e_s
        \quad = \quad \inj^{A,\unit}_0 \cp \id \cp s
        & \quad = \quad \gen{wiring/gt-A-GW-s}
        \\
        e_s'
        \quad = \quad \inj^{A,\unit}_1 \cp \id \cp s
        & \quad = \quad \gen{wiring/gt-GW-zero}
        \\
        e_t
        \quad = \quad \inj^{\unit,B}_1 \cp \id \cp t
        & \quad = \quad \gen{wiring/gt-B-GW-t}
        \\
        e_t'
        \quad = \quad \inj^{\unit,B}_0 \cp \id \cp t
        & \quad = \quad \gen{wiring/gt-GW-zero}
      \end{align*}
    \item[Inductive Step: Tensor]
      Let
      $e = \begin{tikzpicture}[level distance=7mm, sibling distance=7mm, baseline=-4mm]
              \node {$\tensor$}
                child {node {$l$}}
                child {node {$r$}};
            \end{tikzpicture}
      $ be a tree of $n = m_0 + m_1$ leaves.
      By inductive hypothesis, \eqref{equation:wiring-maps-alternative} is equal
      to \Cref{definition:wiring-maps} for $l$ and $r$.
      Then, using $a^l$ and $a^r$ to denote the ancestor maps of $l$ and $r$, we
      have
      \[
        a_L = \gen{wiring/aL-tensor}
        \qquad \qquad
        \text{and}
        \qquad \qquad
        a_R = \gen{wiring/aR-tensor}
      \]
      and so
      \[
        p_s = \gen{wiring/ps-tensor}
        \qquad \qquad
        \qquad \qquad
        p_t = \gen{wiring/pt-tensor}
      \]
      from which we can derive
      \begin{align*}
        \gen{wiring/def-es}
        \quad = \quad \gen{wiring/es-alt-tensor}
        \quad = \quad \gen{wiring/ls-tensor-rs}
        \quad = \quad e_s
        \\ \\
        \gen{wiring/def-es-prime}
        \quad = \quad \gen{wiring/es-prime-alt-tensor}
        \quad = \quad \gen{wiring/ls-prime-tensor-rs-prime}
        \quad = \quad e_s'
        \\ \\
        \gen{wiring/def-et}
        \quad = \quad \gen{wiring/et-alt-tensor}
        \quad = \quad \gen{wiring/lt-tensor-rt}
        \quad = \quad e_s
        \\ \\
        \gen{wiring/def-et-prime}
        \quad = \quad \gen{wiring/et-alt-tensor}
        \quad = \quad \gen{wiring/lt-prime-tensor-rt-prime}
        \quad = \quad e_s
      \end{align*}
      as required.

    \item[Inductive Step: Composition]
      Let
      $e = \begin{tikzpicture}[level distance=7mm, sibling distance=7mm, baseline=-4mm]
              \node {$\cp$}
                child {node {$x$}}
                child {node {$y$}};
            \end{tikzpicture}
      $
      By inductive hypothesis, \eqref{equation:wiring-maps-alternative} is equal
      to \Cref{definition:wiring-maps} for $l$ and $r$.
      Then, using $a^l$ and $a^r$ to denote the ancestor maps of $l$ and $r$, we
      have
      \[
        a_L = \gen{wiring/aL-compose}
        \qquad \qquad
        \qquad \qquad
        a_R = \gen{wiring/aR-compose}
      \]
      and so
      \[
        p_s = \gen{wiring/ps-compose}
        \qquad \qquad
        \qquad \qquad
        p_t = \gen{wiring/pt-compose}
      \]
      from which we can derive
      \begin{align*}
        \gen{wiring/def-es}
        \quad & = \quad \gen{wiring/es-alt-compose}
        \quad = \quad \gen{wiring/ls-inj0}
        \quad = \quad e_s
      \end{align*}
      \begin{align*}
        \gen{wiring/def-es-prime}
        \quad   = \quad \gen{wiring/es-prime-alt-compose}
        \quad & = \quad \gen{wiring/es-prime-alt-compose-mid} \\ \\
        \quad & = \quad \gen{wiring/ls-prime-tensor-rs-plus-rs-prime}
        \quad = \quad e_s'
      \end{align*}
      \begin{align*}
        \gen{wiring/def-et}
        \quad = \quad \gen{wiring/et-alt-compose}
        \quad = \quad \gen{wiring/rt-inj1}
        \quad = \quad e_s
      \end{align*}
      \begin{align*}
        \gen{wiring/def-et-prime}
        \quad = \quad \gen{wiring/et-prime-alt-compose}
        \quad & = \quad \gen{wiring/et-prime-alt-compose-mid} \\ \\
        \quad & = \quad \gen{wiring/lt-prime-plus-lt-tensor-rt-prime}
        \quad = \quad e_t
      \end{align*}
      completing the proof.\qedhere
  \end{description}
\end{proof}

\begin{proposition}
  Let $e$ be a binary tree of $n$ nodes,
  and let $a_L$ and $a_R$ denote the ancestor maps of $e$.
  Computing the wiring maps of $e$ takes $O(n)$ sequential and $O(\log n)$
  parallel time.
\end{proposition}
\begin{proof}
  Assume the ancestor maps $a_L$ and $a_R$ are provided.
  Notice that $p_s$ and $p_t$ are obtained by a stable \emph{integer} sort of
  the array data of $a_L$ and $a_R$.
  This has $O(n)$ sequential and $O(\log n)$ PRAM time complexity: one can use
  e.g., counting sort in the sequential case and radix sort in the parallel
  case.
  The result then follows by the linear (sequential) and $O(1)$ (parallel) time
  complexity of finite function composition.
\end{proof}

\subsection{Computing left and right matching ancestors}
\label{section:fast-wiring-predicate}

In \Cref{section:fast-wiring-reduce}, we reduced the problem of computing the
wiring maps (\Cref{definition:wiring-maps}) to a simple function on trees-
computing the `ancestor maps' (\Cref{definition:ancestor-maps}).
Although in the sequential case there is a straightforward top-down algorithm,
it is not easily parallelizable.

In this section, we give a data-parallel algorithm for computing the closest
ancestor matching a predicate for which a node is in either the left or right
subtree.
This is a generalisation of the problem of computing the ancestor maps.
Our algorithm has $O(\log n)$ PRAM CREW time complexity, which we prove in
\Cref{proposition:lr-p-ancestor-complexity}.

It is first necessary to represent trees in a `parallel friendly' way.
A number of suitable array-based representations of trees exist
(e.g., \cite[p. 128]{introduction-to-algorithms})
having constant-time operations for computing parent and child indices of a given
node.
Thus, instead of fixing a specific array-based representation of a tree, we
simply assume that the following finite functions can be computed in linear
(sequential) and constant (parallel) time.

\begin{definition}
  \label{definition:parent-left-right}
  Let $e$ be a binary tree of size $n = 2m - 1$
  (i.e., $e$ has $m$ leaves and $m - 1$ nodes.)
  Denote by $\parent : n \to n + 1$ the finite function mapping a node to its
  parent, or $n$ if none exists.
  We write $\isRight : n \to 2$ for the predicate function returning $1$ if a
  node is the \emph{right} child of a node, and $0$ otherwise.
  Similarly, write $\isLeft : n \to 2$ for the predicate returning $1$ if a node
  is the \emph{left} child of a node.
\end{definition}

\begin{example}
  The actions of $\parent$, $\isRight$, and $\isLeft$ are illustrated with the
  following example, in which the nodes of the tree are labeled $0 \ldots 4$.
  \[
    \input{wiring/example-parent-left-right.tikz}
  \]
\end{example}

Recall that the purpose of the ancestor maps
(\Cref{definition:ancestor-maps})
is to compute the closest ancestor labeled $\cp$ for each leaf in the tree.
The algorithm presented here solves a slightly more general problem: computing
the closest ancestor matching a predicate $P$ for which a node is in the left or
right subtree.
Inductive definitions of these functions are given below.

\begin{definition}[Left/Right $P$-Ancestor Function]
  \label{definition:left-right-p-ancestor}
  Let $e$ be a binary tree of size $n$, and let $P : n \to 2$ be a predicate on
  nodes.
  The \deftext{Left $P$-Ancestor} function is denoted $\ancestor_L$, and defined below.
  \begin{align*}
    & \ancestor_L : n \to n+1 \\
    & \ancestor_L(i) = \begin{cases}
      n                                     & \qquad \text{$\isRoot(i)$} \\
      \parent(i)                            & \qquad \text{$P(\parent(i)) \lor \isRight(i)$} \\
      \ancestor_L(\parent(i)) + \isRight(i) & \qquad \text{otherwise}
    \end{cases}
  \end{align*}
  Similarly, the \deftext{Right $P$-Ancestor} function is given by
  \begin{align*}
    & \ancestor_R : n \to n+1 \\
    & \ancestor_R(i) = \begin{cases}
      n                                     & \qquad \text{$\isRoot(i)$} \\
      \parent(i)                            & \qquad \text{$P(\parent(i)) \lor \isLeft(i)$} \\
      \ancestor_R(\parent(i)) + \isLeft(i) & \qquad \text{otherwise}
    \end{cases}
  \end{align*}
\end{definition}

The following example illustrates the action of the $\ancestor_L$ map.

\begin{example}
  Pictured below is a tree in which each node $i$ is drawn in
  black if $P(i) = 1$ and white otherwise,
  and the value of $\ancestor_L(i)$ is shown for each node.
  \[ \input{wiring/example-ancestorL.tikz} \]
\end{example}

Our approach to efficiently computing $\ancestor_L$ and $\ancestor_R$
is `bottom-up'.
We first transform the tree of size $n$ to a functional graph with $2n+1$ nodes,
and then iterate its adjacency function until the result converges.
We illustrate this process with the example below.

\begin{example}
  Pictured below is a tree with $n = 5$ nodes,
  and its `left ancestor graph', with $2n + 1 = 11$ nodes.
  \[ \input{wiring/example-tree.tikz} \qquad \qquad \Longrightarrow \qquad \qquad \input{wiring/example-ancestor-graph.tikz} \]
  For each node $i$ in the tree there are two nodes in the graph, $2i$ and $2i + 1$.
  The edges of the graph represent the parent relation, split into two distinct
  cases: when a child is a \emph{left} or \emph{right} child of its parents.
  Left children of the node $i$ are adjacent to the graph node $2i$, while
  right children are adjacent to $2i + 1$.
  When a $P(i)$, its corresponding graph node $2i + 1$ points to itself instead of to its parent.
  Thus, when iterating the adjacency relation, eventually all right-descendents are `captured' by their closest left-ancestor.
\end{example}

This graph is defined more formally below.

\begin{definition}[Ancestor Graph]
  Let $e$ be a tree with $m$ leaves and $n = 2m - 1$ nodes,
  and $P : n \to 2$ a predicate on nodes.
  The \deftext{Left Ancestor Graph} of $e$ has $2n + 1$ vertices, and edges given by the following adjacency function.
  \begin{align*}
    & r_L : 2n+1 \to 2n+1 \\
    & r_L(j) = \begin{cases}
      j                                 & j = 2n \\
      j                                 & \isOdd(j) \land P\left(\left\lfloor \frac{j}{2} \right\rfloor\right) \\
      2 \cdot \parent\left(\left\lfloor \frac{j}{2} \right\rfloor\right)
        + \isRight\left(\left\lfloor \frac{j}{2} \right\rfloor\right)
        & \text{otherwise} \\
    \end{cases}
  \end{align*}
  Symmetrically, the \deftext{Right Ancestor Graph} has adjacency relation
  \begin{align*}
    & r_R : 2n+1 \to 2n+1 \\
    & r_R(j) = \begin{cases}
      j                                 & j = 2n \\
      j                                 & \isEven(j) \land P\left(\left\lfloor \frac{j}{2} \right\rfloor\right) \\
      2 \cdot \parent\left(\left\lfloor \frac{j}{2} \right\rfloor\right)
        + \isRight\left(\left\lfloor \frac{j}{2} \right\rfloor\right)
        & \text{otherwise} \\
    \end{cases}
  \end{align*}
\end{definition}

We can now relate the ancestor graph and the functions $\ancestor_L$ and
$\ancestor_R$.

\begin{proposition}
  \label{proposition:ancestor-as-graph}
  For all nodes $i$ of depth $\depth(i)$ in a tree $e$
  \[
    \ancestor_L(i) = \left\lfloor \frac{r_L^k(2i)}{2} \right\rfloor
    \qquad \qquad \qquad
    \ancestor_R(i) = \left\lfloor \frac{r_R^k(2i)}{2} \right\rfloor
  \]
  for $k > \depth(i)$
  where $r$ is the adjacency function of the ancestor graph.
\end{proposition}
\begin{proof}
  Top-down induction by depth.
  We write $r$ for $r_L$ and prove only the former case, omitting the symmetric
  proof for the $\ancestor_R$ function.
  \begin{description}
    \item[Base Case]
      Suppose $\isRoot(i)$ so $\depth(i) = 0$.
      We have
      $r(2i) = 2 \cdot \parent(i) + \isRight(i) = 2 \cdot \parent(i) = 2n$,
      and so
      \[ \left\lfloor \frac{r(2i)}{2} \right\rfloor = n = \ancestor_L(i) \]
      Moreover, $2n$ is a fixed point for $r$, so this holds for any $r^k$
      where $k \geq 1$.
    \item[Inductive Case 1]
      In the second case, we have $P(\parent(i)) \land \isRight(i)$.
      Suppose $\depth(i) = m$, where $m \geq 1$ because $i$ is not the root node.
      We have
      \[
        r(2i) = 2 \cdot \parent(i) + \isRight(i) = 2 \cdot \parent(i) + 1
      \]
      Now let $k \geq 1$ so that $m + k > \depth(i)$.
      Then
      \[
        \left\lfloor \frac{r^{m+k}(2i)}{2} \right\rfloor
        = \left\lfloor \frac{r^{m+k-1}(2\cdot\parent(i) + 1)}{2} \right\rfloor
        = \left\lfloor \frac{2\cdot\parent(i) + 1}{2} \right\rfloor
        = \parent(i)
        = \ancestor_L(i)
      \]
      as required.
    \item[Inductive Case 2]
      In the final case,
      let $i$ be a node,
      and let $k > \depth(i)$
      so that $k - 1 > \depth(parent(i))$
      and the inductive hypothesis is as below.
      \[
        \left\lfloor \frac{r^{k-1}(2 \cdot \parent(i))}{2} \right\rfloor
        =
        \ancestor_L(\parent(i))
      \]
      Now split into two further cases: (a) when $\neg \isRight(i)$ and
      (b) when $\neg P(\parent(i)) \land \isRight(i)$.
      \begin{description}
        \item[Case (a)]
          Observe that $r(2i) = 2 \cdot \parent(i) + \isRight(i) = 2 \cdot \parent(i)$.
          Then
          \[
            r^k(2i) = r^{k-1}(r(2i)) = r^{k-1}(2\cdot\parent(i))
          \]
          from which we can derive
          \[
            \left\lfloor \frac{r^k(2i)}{2} \right\rfloor
            =
            \left\lfloor \frac{r^{k-1}(2\cdot\parent(i))}{2} \right\rfloor
            =
            \ancestor_L(\parent(i))
            =
            \ancestor_L(i)
          \]
        \item[Case (b)]
          When $\neg P(\parent(i)) \land \isRight(i)$,
          we have
          $r(2i) = 2 \cdot \parent(i) + \isRight(i) = 2 \cdot \parent(i) + 1$
          and
          \begin{align*}
            r(2 \cdot \parent(i) + 1)
            & = 2 \cdot \parent\left(\left\lfloor
                \frac{2 \cdot \parent(i) + 1}{2}
              \right\rfloor\right)
              + \isRight\left(\left\lfloor
                \frac{2 \cdot \parent(i) + 1}{2}
              \right\rfloor\right) \\
            & = 2 \cdot \parent(\parent(i)) + \isRight(\parent(i)) \\
            & = r(2 \cdot \parent(i)) \\
          \end{align*}
          where we use that $\neg P(\parent(i))$ in the first step.
          One can then derive the result immediately as follows.
          \begin{align*}
            \left\lfloor \frac{r^k(2i)}{2} \right\rfloor
            =
            \left\lfloor \frac{r^{k-1}(2 \cdot \parent(i) + 1)}{2} \right\rfloor
            =
            \left\lfloor \frac{r^{k-1}(2 \cdot \parent(i))}{2} \right\rfloor
            & =
            \ancestor_L(\parent(i)) \\
            & =
            \ancestor_L(i)
          \end{align*}
      \end{description}
  \end{description}
\end{proof}

\begin{corollary}
  Let $e$ be a tree of size $n$.
  Then
  \[ \ancestor_L(i) = \left\lfloor \frac{r^n(2i)}{2} \right\rfloor \]
\end{corollary}
\begin{proof}
  The maximum depth of a leaf in a tree of $n$ nodes is $n - 1$.
  Thus $n > \depth(i)$ for all nodes $i$,
  and so the equality holds by \Cref{proposition:ancestor-as-graph}.
\end{proof}

\begin{corollary}
  \label{proposition:lr-p-ancestor-complexity}
  Let $e$ be a tree of size $n$.
  Computing the finite functions $\ancestor_L$ and $\ancestor_R$
  for $e$ has $O(\log n)$ PRAM CREW time complexity.
\end{corollary}
\begin{proof}
  Let $r$ denote either $r_L$ or $r_R$.
  Computing $r$ is $O(1)$ because each entry of the underlying array can be
  computed in parallel in constant time.
  A composition $r \cp r$ is constant PRAM CREW time because each entry of the
  resulting array is a single lookup, with each entry computed in parallel on
  $O(2n+1)$ processors.
  Thus computing $r^n$ is $O(\log n)$ by `repeated squaring'.
  That is, let $f^0 = r$, $f^{n+1} = f^n \cp f^n$, then $f^{\log n} = r^n$.
\end{proof}

\begin{remark}
  The technique of `repeated squaring' a functional graph
  is known as `pointer jumping' in the parallel computing literature.
  See for example \cite[2.2]{introduction-to-parallel-algorithms}.
  Note also that while its time complexity is logarithmic, the algorithm given
  here is not work-efficient and requires $O(n \log n)$ operations.
  We expect this can be improved to $O(n)$, but leave this to future work.
\end{remark}

Finally, note that the ancestor maps $a_L$ and $a_R$ of
\Cref{definition:ancestor-maps} are not exactly the same functions as those in
\Cref{definition:left-right-p-ancestor}.
However, it is a straightforward matter of elementwise arithmetic to transform
between the two finite functions, and so we omit the details.

%% file: applying-functors.tex
We now give an efficient, parallel algorithm for applying strict symmetric
monoidal hypergraph functors to diagrams.
We will hereafter just say `functor' and assume the strict symmetric monoidal
hypergraph structure.
Note that in what follows, we assume arbitrary signatures $\Sigma$ and $\Omega$,
which optionally include the chosen Special Frobenius structure $\Frob$.

The algorithm given here is defined \emph{without} having to first decompose a
diagram into a $\Signature$-term.
The example below illustrates the importance of this for the case of diagrams in
$\DiagramSMC$ which are not equipped with hypergraph structure.

\begin{example}
  \label{example:quadratic-penalty}
  Suppose we are working in the category $\Free_\CMon$.
  The terminal map $\terminal_N : N \to 1$ is represented as the following
  string diagram.
  \[ \gen{examples/example-n2-problem-n-fold-sum} \]
  Represented combinatorially as a diagram, this has $O(N)$ internal wires,
  generators, and edges.
  However, suppose we decompose it by taking `vertical slices',
  then the number of terms in the decomposition is $O(N^2)$.
  The first `slice' is $\sgen{cmon/g-add} \tensor \id \tensor \ldots \tensor \id$,
  having $N - 2$ copies of the identity.  The second `slice' has $N-3$ copies,
  and so on,
  so the average size of each slice is $N/2$ for $N$ slices.
\end{example}

Thus, by naively decomposing a term into a $\Signature$-term, it is possible to
incur a quadratic penalty in the representation size.
Our approach is to avoid this conversion altogether.
That is, we apply functors directly to diagrams, without `roundtripping' through
$\Signature$-terms.
As a bonus, in \Cref{section:optic-composition}, we will show how this can be
used to give a fast algorithm for taking reverse derivatives (and more generally
exploiting hypergraph structure to map diagrams to diagrams of \emph{optics}.)

Our algorithm is based on the Frobenius decompositions defined in
\Cref{proposition:frobenius-decomposition-existence}.
These decompositions put each morphism into a form consisting of Frobenius
spiders and a tensoring of generators.
Given a strict symmetric monoidal hypergraph functor between categories of
diagrams $F : \Diagram_\Sigma \to \Diagram_\Omega$,
our approach is to map Frobenius decompositions in $\Diagram_\Sigma$ to
Frobenius decompositions in $\Diagram_\Omega$.
Thus, in order to give an algorithm for the action of $F$, it suffices to define
two things (1) the action on spiders, and (2) the action on (a tensoring of)
generators.

The section is structured as follows.
In \Cref{section:objects-and-spiders} we define an auxiliary
datastructure (the `segmented finite function').  We use this to encode the
object map of $F$ and give an algorithm for applying $F$ Frobenius spiders of
$\Diagram_\Sigma$.
In \Cref{section:mapping-tensorings} we then show how the same datastructure can
be used to encode the \emph{arrow} map of $F$, and similarly give an algorithm
applying $F$ to a tensoring of generators.
We conclude by showing that applying an encoded functor $F$ has linear
(sequential) and logarithmic (parallel) time complexity.

\subsection{Objects, (Half-)Spiders, and Segmented Finite Functions}
\label{section:objects-and-spiders}

Given signatures $\Sigma$ and $\Omega$ and a functor
$F : \Diagram_\Signature \to \Diagram_\Omega$,
we must encode the data of $F$ in a way that allows us its application to
diagrams to be parallelised.
We make two simplifying assumptions:
(1) that signatures $\Sigma$ and $\Omega$ are finite,
and (2) that the typing relation of both is a \emph{function}.
(This latter assumption is somewhat more restrictive than required; we discuss a
weaker condition in \Cref{section:conclusion}.)

A functor $F : \Diagram_\Sigma \to \Diagram_\Omega$ consists of two maps.
A function on \emph{objects}
$F_0 : \Sigma_0^* \to \Omega_0^*$
and a function on \emph{arrows}
$F_1 : \Diagram_\Sigma(A, B) \to \Diagram_\Omega(F(A), F(B))$.
In this section, we introduce an auxiliary data structure, the `segmented finite
function', which will allow us to encode these maps in a manner suitable for
parallel application.

To motivate this, consider the object map $F_0$.
Since $F$ is a \emph{strict} symmetric monoidal hypergraph functor between
categories of diagrams, we have the following.
\[
  F_0 \left( \sum_{i \in N} A_i \right) = \sum_{i \in N} F_0(A_i)
\]
We may therefore consider the object map of $F$ to have the type
$F_0 : \Sigma_0 \to \Omega^*$
since (by strictness) it is completely defined by its action on generating
objects.
Similarly, we think of the map on arrows as having the type
$F_1 : \Sigma_1 \to \Diagram_\Omega$.\footnote{
  Note that strictly speaking morphisms of $\Diagram_\Omega$ are equivalence classes,
  but the map $F_1$ maps a generator to specific diagram.
  We gloss over this detail to avoid introducing new notation.
}

Now, for each $i \in \Sigma_0$, we can think of $F_0(i) : \Omega_0^*$
as a finite function $F_0(i) : s(i) \to \Omega_0$
where $s : \Sigma_0 \to K$ is a function denoting the length (source) of each
list $s(i) = \len{F_0(i)}$, and $K$ is the maximum such length.
Encoding the object map $F_0$ will mean storing the data of these finite
functions in a `flat' array structure, with the array $s$ encoding their
sources.

\begin{definition}
  Let $f_i : A_i \to B_i$ be a collection of $N$ finite functions.
  A \deftext{segmented finite function} encodes their sources, targets, and array data
  as the following three maps.
  \begin{align*}
    & \sources    : N \to \max_{i \in N} \source(f_i)
    \qquad \qquad
    & \targets    : N \to \max_{i \in N} \target(f_i)
    \\
    & \sources(i) = \source(f_i)
    \qquad \qquad
    & \targets(i) = \target(f_i)
  \end{align*}
  \begin{align*}
    \values     & : \sum_{i \in N} \sources(i) \to \max_{i \in N} \target(f_i) \\
    \values     & = \sum_{i \in N} (f_i \cp \inj_0)
  \end{align*}
\end{definition}

\begin{remark}
  In parallel programming terminology, $\values$ is a \emph{segmented array}
  with segment sizes are given by $\sources$.
  Note that we also store the $\targets$ of each morphism, which will be
  necessary to take \emph{tensor products}.
\end{remark}

Now, given an indexing function, we can then take arbitrarily ordered coproducts
of the functions $f_i$ as follows.

\begin{proposition}
  \label{proposition:indexed-coproduct}
  Let $f_i : A_i \to B$ be a collection of $N$ finite functions sharing a codomain,
  and $x : X \to N$ a finite function.
  Then the $x$-indexed coproduct is calculated as follows.
  \begin{align*}
    \sum_{i \in X} f_{x(i)} & : \sum_{i \in X} A_{x(i)} \to B \\
    \sum_{i \in X} f_{x(i)} & = \left( \sum_{i \in X} \inj_{f(i)} \right) \cp \values
  \end{align*}
\end{proposition}
\begin{proof}
  We may calculate
  $
  \left( \sum_{i \in X} \inj_{x(i)} \right) \cp \values
      = \sum_{i \in X} f_{x(i)} \cp \inj_0
      = \sum_{i \in X} f_{x(i)}
  $
  because $\inj_0 = \id$ since targets of all $f_i$ are equal.
  Note also that the source $\sum_{i \in X} A_{x(i)}$ can be computed
  by the sum of entries of the array
\end{proof}

This kind of indexed coproduct is precisely what we need to apply the object map
$F_0$ of a functor to \emph{Frobenius spiders}.
We begin by defining precisely how the object map is encoded.

\begin{definition}
  Given finite signatures $\Sigma$ and $\Omega$
  and a strict symmetric monoidal hypergraph functor
  $F : \Diagram_\Sigma \to \Diagram_\Omega$,
  the \deftext{object map encoding} of $F$
  is the segmented finite function
  with $\sources = s$, $\values = v$, and $\targets(i) = \Omega$,
  where
  \begin{align*}
    s : \Sigma_0 \to K
    \qquad \qquad \qquad
    &
    v : \sum_{i \in \Sigma_0} s(i) \to \Omega_0
    \\
    s(i) = \len{F_0(i)}
    \qquad \qquad \qquad
    &
    v = \sum_{i \in \Sigma_0} F_0(i)
  \end{align*}
\end{definition}

In the above, the injections of the coproduct
$\source(v) = \sum_{i \in \Sigma_0} s(i)$
have the type
$\inj_x : s(x) \to \sum_{i \in \Sigma_0} s(i)$.
By definition of the coproduct, we have the following.
\begin{align*}
  \inj_x \cp v & : s(x) \to \Omega_0 \\
  \inj_x \cp v & = F_0(x)
\end{align*}
In other words, precomposing the segmented array of values with an injection
amounts to taking a slice of the array corresponding to a specific `segment'.

With the object map $F_0$ suitably represented, we can now show how to map the
data of a Frobenius spider $f = (s, t, \SCL(B))$ in $\Diagram_\Sigma$ to its
corresponding spider in $\Diagram_\Omega$.
Recall that any Frobenius spider can be regarded as the composite of two
\emph{half-spiders} (\Cref{definition:half-spider}):
$f = (s, \id, \SCL(B)) \cp (\id, t, \SCL(B))$.
It therefore suffices to determine the action of $F$ on half-spiders.

Since half-spiders are in the image of the functor $\HalfSpider$,
we may regard them simply as morphisms of $\Wires$.
The action of $F$ on such morphisms is as follows.

\begin{proposition}
  Let $f : A \to B$ be a morphism of $\Wires$,
  and suppose $F$ is a strict symmetric monoidal hypergraph functor
  whose object encoding is the segmented finite function
  with $\sources = s$ and $\values = v$.
  Then $F(\HalfSpider(f)) = \HalfSpider(f') : A' \to B'$, where
  \begin{align*}
    B'(\swn) & : \sum_{b \in B(W)} s(B(\swn)(b)) \to \Omega_0
    \qquad
    &
    \qquad
    f'     & : \sum_{a \in A(W)} s(A(\swn)(a)) \to \sum_{b \in B(W)} s(B(\swn)(b))
    \\
    B'(\swn) & \defeq \left(\sum_{b \in B(W)} \inj_{B(\swn)(b)} \right) \cp v
    \qquad
    &
    \qquad
    f'     & \defeq \sum_{a \in A(W)} \inj_{f_W(a)}
  \end{align*}
  and $A'(\swn) = f' \cp B'(\swn)$ is fixed by naturality.
\end{proposition}
\begin{proof}
  The definition of $B'$ is required by the strictness assumption.
  We must then verify that $A'(\swn) = f' \cp B'(\swn)$,
  which follows by naturality of $f$:
  \begin{align*}
    f'_W \cp B'(\swn)
      & : \sum_{a \in A(\swn)} s(A(\swn)(a)) \to \Omega_0 \\
      & : \sum_{a \in A(\swn)} s(B(\swn)(f(a))) \to \Omega_0 \\
    f'_W \cp B'(\swn)
      & =
        \left(\sum_{a \in A(W)} \inj_{f_W(a)}\right)
        \cp \left(\sum_{b \in B(W)} \inj_{B(\swn)(b)}\right)
        \cp v \\
      & = \left(\sum_{a \in A(W)} \inj_{B(\swn)(f(a))}\right) \cp v \\
      & = \left(\sum_{a \in A(W)} \inj_{A(\swn)(a)}\right) \cp v
  \end{align*}
  Lastly, we must check that the hypergraph structure is preserved, and that $F$
  satisfies the functor laws.
  This is done by induction: in the base case, it is straightforward to verify
  that the definition holds for generating operations in the image of
  $\HalfSpider$.
  The inductive step for tensor product is similarly straightforward.
  One can verify the case for composition as follows.
  Suppose $A\overset{f}{\to}B\overset{g}{\to}C$ are morphisms of $\Wires$.
  Then
  \begin{align*}
    f' \cp g' & : \sum_{a \in A(W)} s(B(\swn)(f(a))) \to \sum_{c \in C(W)} s(C(wn)(c)) \\
    f' \cp g'
      & = \left( \sum_{a \in A(W)} \inj_{f(a)} \right)
          \cp \left( \sum_{b \in B(W)} \inj_{g(b)} \right)
        = \sum_{a \in A(W)} \inj_{g(f(a))}
        = (f \cp g)'
  \end{align*}
  Finally, since $F$ is a hypergraph functor, we have $F(\HalfSpider(f)) =
  \HalfSpider(f')$ by uniqueness of $\HalfSpider$.
\end{proof}
\begin{remark}
  Note carefully that the injections in the definition of $B'$ and $f'$ above
  are not the same.
  The indices of injections used in the definition of $B'$ range over the
  segments of $v$,
  whereas those in the definition of $f'$ range over elements of $B(W)$.
\end{remark}

The object $B'(\swn)$ of the mapped spider is precisely a $B(\swn)$-indexed
coproduct of the values table.
Similarly, the value of $f'$ is given by `indexed injections'.
In order to compute both, we use the following algorithm.
Given a finite function $s : N \to K$
and a map $x : X \to N$ thought of as an array of indices,
\texttt{injections(s, x)} computes a coproduct of injections
$\inj_{x(0)} \ldots \inj_{x(A-1)}$.

\begin{minted}{python}
def injections(s, x):
    p = prefix_sum(s)
    r = segmented_arange(x >> s)
    return FiniteFunction(sum(s), r.table + repeat(x >> p, x >> s).table)
\end{minted}

Note carefully that the $+$ operation in the above algorithm denotes pointwise
sum of integer arrays, not the coproduct.
Given the object map encoding of $F$,
one can then compute
$f' = \mathtt{injections(s, f >> B(\swn))}$
and
$B'(\swn) = \mathtt{injections(s, B(\swn)) >> v}$.

\begin{proposition}
  \label{proposition:complexity-indexed-injections}
  Let $s : N \to K$
  and $x : X \to N$ be finite functions.
  The sequential time complexity of computing $\mathtt{injections}$ is
  $O(N) + O(X) + O(\mathtt{sum}(x \cp s))$
  and its PRAM CREW complexity is
  $O(\log N) + O(\log X)$.
\end{proposition}
\begin{proof}
  In the sequential case, both function compositions are $O(X)$,
  prefix sum is $O(N)$, and remaining operations are $O(\mathtt{sum}(x \cp s))$.
  In the parallel case, the same operations are
  $O(\log N)$, $O(1)$, and $O(\log X)$.
\end{proof}

Note that when computing $f'$ and $B'(\swn)$ using the $\mathtt{injections}$
function, $\mathtt{sum}(x \cp s)$ is proportional to $B'(W)$, and so the
computing spiders is linear (sequential) in the size of the result.

\subsection{Mapping Tensorings of Generators}
\label{section:mapping-tensorings}

Having defined the action of $F$ on the \emph{spiders} of a Frobenius decomposition,
it now remains to define its action on a \emph{tensoring of generators}.
Such a tensoring is a diagram of the form
$g = (\inj_0, \inj_1, G)$
whose data is defined up to isomorphism by the list $G(\sxn) : G(X) \to \Sigma_1$.
In other words, the isomorphism class of diagrams represented by $g$ is in the
image of the $\Sigma$-term
$\bigtensor_{i \in G(X)} G(\sxn)(i)$ under $F$,
and so by strictness it is required that
$F(g) \cong \bigtensor_{i \in G(X)} F_1(G(\sxn)(i))$.

As with the object map, the action of $F$ on arrows is then completely defined
by its action on operations $\Sigma_1$.
We can therefore consider it to have the type
$F_1 : \Sigma_1 \to \Diagram_\Omega$.
This leads to a straightforward implementation in the sequential case.
Encode the data of $F_1$ as a list of length $\len{\Sigma_1}$
whose $i^{\text{th}}$ entry is the diagram $F_1(i)$,
then in pseudocode, we might write

\begin{minted}{python}
  def apply_functor(F1: List[Diagram], g: Diagram):
    diagrams = [ F1[g.G.xn(i)] for i in range(0, g.X) ]
    Diagram.tensor_list(diagrams)
\end{minted}

The $N$-fold tensor of diagrams is then a straightforward extension of the
binary case.
Recall that a diagram is essentially a collection of finite maps, upon which
coproducts and tensor products are pointwise.
Thus, taking the tensor product of a list of diagrams amounts to taking
\emph{finite} coproducts instead of binary coproducts.
The parallel case is `morally' the same, but obtaining a proper logarithmic time
algorithm takes special care: computing the finite tensor product of a list of
finite functions does not immediately translate to the parallel case.
Consider the following naive implementation.

\begin{minted}{python}
def coproduct_list(fs: List[FiniteFunction]):
    return FiniteFunction(fs[0].cod, concatenate([f.table for f in fs]))
\end{minted}

The problem is that constructing the argument to $\mathtt{concatenate}$ takes
time linear in the size of $\mathtt{fs}$.
Since the length of $\mathtt{fs}$ equals $G(X)$, this takes time linear in the
number of operations in the diagram, and therefore does not enjoy a speedup in
the parallel case.
The same problem exists for the tensor product of finite maps.

The solution is to encode the data of $F_1$ not as a list of diagrams, but as
a `diagram of lists'.
More precisely, we instead encode each diagram component as a separate segmented
finite function.
From this representation we will be able to extract arbitrary coproducts and
tensors of the `segments' to construct the result.

The encoding of $F_1$ is a collection of segmented finite functions: one for
each of the source and target maps, and one for each of the components $G(f)$
for $f \in \SchemaBPMG$.

\begin{definition}
  Let $\Sigma$ and $\Omega$ be finite signatures,
  and
  $F : \Diagram_\Sigma \to \Diagram_\Omega$
  a strict symmetric monoidal hypergraph functor,
  and denote by
  $d_i = (s_i, t_i, G_i)$ the collection of diagrams
  $d_i = F_1(i)$
  for $i \in \Sigma_1$.
  The \deftext{arrow map encoding} of $F$
  consists of segmented finite functions for
  $s_i$, $t_i$, and $G_i(f)$ for each component $f$ in the schema $\SchemaBPMG$.
\end{definition}

In order to compute tensor products of diagrams, it will be necessary to compute
indexed tensor products from segmented finite functions.

\begin{proposition}
  Let $f_i : A_i \to B_i$ be a collection of $N$ finite functions, and
  $x : X \to N$ a finite function indexing the collection.
  Then the indexed tensor product is given by the array
  \begin{align*}
    \bigotimes\limits_{i \in X} f_{x(i)}
      & : \sum_{i \in X} A_{x(i)} \to \sum_{i \in X} B_{x(i)} \\
      & = \sum_{i \in X} f_{x(i)} \oplus \mathsf{repeat}(x \cp \sources, p)
  \end{align*}
  where $\oplus$ denotes elementwise addition of natural numbers
  and $p$ denotes the partial sums of codomains
  $
    p(i) = \sum_{j \in \ord{i}} \targets(x(j)) %
  $.
\end{proposition}
\begin{proof}
  Recall that the tensor product can be written in terms of the coproduct as below.
  \[ \bigotimes_{i \in X} f_{x(i)} = \sum_{i \in X} (f_{x(i)} \cp \inj_i) \]
  Here, each injection has type
  $
    \inj_i : B_{x(i)} \to \sum_{j \in X} B_{x(j)} %
  $
  with array data given by
  $
    \bl p(i), p(i)+1, \ldots, p(i) + B_i - 1 \br
  $
  The array data of the composition $f_{x(i)} \cp \inj_i$ is therefore
  that of $f_{x(i)}$ plus the constant $p(i)$,
  which we may rewrite as below
  \begin{align*}
    \sum_{i \in X} (f_{x(i)} \cp \inj_i)
      & = \sum_{i \in X} f_{x(i)} \oplus \mathsf{repeat}(x \cp \sources, p)
  \end{align*}
\end{proof}

Computing $F(h)$ for a tensoring $h = (\inj_0, \inj_1, H)$ is then
straightforward.
Since coproducts and tensor products are pointwise, each finite function con the
data of $F(h)$ is the $H(\sxn)$-indexed coproduct (or tensor product) of its
corresponding segmented finite function.

\begin{proposition}
  \label{proposition:complexity-arrow-map}
  Complexity of applying the arrow map to an $N$-fold tensoring of generators
  $(\inj_0, \inj_1, g)$ is $O(N)$ (sequential) and $O(\log N)$ (PRAM) in the
  number of operations of $g$.
\end{proposition}
\begin{proof}
  As in \Cref{proposition:complexity-fast-functor}, the number of wires and
  edges in the diagram is proportional to $N$.
  Applying $F_1$ consists of computing the $G(\sxn)$-indexed co- and tensor
  products for each component morphism of $G$,
  and so complexity is the same as these operations,
  which are $O(N)$ (sequential) and $O(\log N)$ (PRAM) by
  \Cref{proposition:complexity-indexed-injections}.
\end{proof}

\begin{corollary}
  Let $F : \Diagram_\Sigma \to \Diagram_\Omega$
  be a strict symmetric monoidal hypergraph functor,
  and $d = (s, t, G)$ be a diagram in $\Diagram_\Sigma$ of type $d : A \to B$.
  Given object and arrow map encodings,
  the sequential time complexity of computing $F(d)$ is linear
  in the number of wires $G(W)$, edges $G(\Ei)$ and $G(\Eo)$, operations
  $G(X)$, and boundaries $A(W)$ and $B(W)$ of $d$.
  The parallel time complexity is logarithmic in the same.
\end{corollary}

%% file: optic-composition.tex
We now show how the Hypergraph structure of $\DiagramFrob$ can be used to encode
\emph{optics}.
Our approach is inspired by the string diagrams of~\cite{optic-string-diagrams}
and~\cite{bruno}.
Along with the parallel algorithm for functor application described in
\cref{section:applying-functors}, this allows us to define an algorithm for
mapping diagrams into diagrams of optics.
In general, this allows for modelling systems with bidirectional information
flow.
An example of such systems are neural networks~\cite{cfgbl} viewed as morphisms
in categories with \emph{reverse derivatives}~\cite{rdc}.
We explore this example specifically, and define a scalable method for taking
reverse derivatives of large diagrams.
Moreover, this addresses an efficiency issue with the naive approach to taking
reverse derivatives as pointed out in~\cite{bruno}.

As in \Cref{section:applying-functors}, we assume arbitrary signatures $\Sigma$
and $\Omega$, which are now assumed to include the additional chosen Special
Frobenius structure.
We begin by informally recalling optics.

\begin{definition}[Informal, \cite{bruno,categories-of-optics}]
  \label{definition:optic}
  An \deftext{optic} of type
  $\OpticObj{\fwd{A}}{\rev{A}} \to \OpticObj{\fwd{B}}{\rev{B}}$
  is a triple
  \[
    M \in \cat{C}
    \qquad \qquad
    \fwd{f} : \fwd{A} \to \fwd{B} \tensor M
    \qquad \qquad
    \rev{f} : M \tensor \rev{B} \to \rev{A}
  \]
  where we call $\fwd{f}$ and $\rev{f}$ the \deftext{forward} and \deftext{reverse}
  map, respectively.
\end{definition}

The whole optic is considered as a system modelling forward and backward
information flow.
The $\fwd{f}$ morphism captures the `forward'
information flow of the model, mapping data $\fwd{A}$ to $\fwd{B}$.
The object $M$ is some `memory' of the input $\fwd{A}$,
which is passed to the `reverse' map $\rev{f}$,
which maps `output-like' data and \emph{memory} back into `input-like' data.
In the specific case of reverse derivatives
(see \Cref{section:reverse-derivatives}),
`input-like' and 'output-like' will mean
\emph{changes} in input and output, respectively.

Note that, strictly speaking, morphisms in categories of optics are
\emph{equivalence clases} of triples $(M, f, f')$.
This will not concern us.
Instead, the focus of this section is on how representatives of such
classes--specific triples $(M, f, f')$--can be encoded as diagrams.
The main contribution of this section is to recognise that optic composition can
be `simulated' using the hypergraph structure present in
$\DiagramFrob$.
In order to state this formally, it is first necessary to define the
\emph{interleave} morphism.
\begin{definition}
  Assume for each generator $X \in \Sigma_0$ a pair of objects
  $\fwd{X} \in \Omega_0^*$ and $\rev{X} \in \Omega_0^*$.
  Then let $A = A_0 \tensor A_1 \tensor \ldots \tensor A_{N-1}$ be a tensoring of
  generating objects in $\Sigma_0$.
  The \deftext{interleaving} at $A$ is the permutation with the following type.
  \[
    \smcinterleave_A
      :               \bigtensor_{i \in N} \fwd{A_i} \tensor \bigtensor_{i \in N} \rev{A_i}
      \longrightarrow \bigtensor_{i \in N} (\fwd{A_i} \tensor \rev{A_i})
  \]
\end{definition}

Given a chosen forward and reverse map for each generating operation $f : A \to B$,
this interleaving makes it possible to define a strict monoidal functor
taking each $f$ to its chosen optic.

\begin{definition}
  Let $\Sigma$ and $\Omega$ be monoidal signatures, and assume the following
  data is given.
  \begin{itemize}
    \item For each generating object $A \in \Sigma_0$, a pair of objects (lists)
        $\fwd{A} \in \Omega_0^*$
        and
        $\rev{A} \in \Omega_0^*$
    \item For each generating operation
      $f : \bigtensor_{i \in N_0} A_i \to \bigtensor_{i \in N_1} B_i$
      in $\Sigma_1$,
      a pair of morphisms
      \begin{itemize}
        \item $\fwd{f} : \bigtensor_{i \in N_0} \fwd{A_i}
                     \longrightarrow
                         \left(\bigtensor_{i \in N_1} \fwd{B_i}\right) \tensor M$,
        \item $\rev{f} : M \tensor \left(\bigtensor_{i \in N_1} \rev{B_i} \right)
                     \longrightarrow
                         \bigtensor_{i \in N_0} \rev{A_i}$
      \end{itemize}
  \end{itemize}
  Denote by
  $\ToOptic : \Diagram_\Sigma \to \Diagram_\Omega$
  the strict symmetric monoidal hypergraph functor defined inductively
  whose action on generating objects is
  $F(A) = \fwd{A} \tensor \rev{A}$
  and on arrows is given below.
  \[
    F(f) \qquad \defeq \qquad \gen{optics/interleaved-optic-f}
  \]
\end{definition}

Note that diagrams in the image of $\ToOptic$ compose in the same way as optics.
\[
  \ToOptic(f \cp g) \quad = \quad
  \ToOptic(f) \cp \ToOptic(g) \quad = \quad
  \gen{optics/interleaved-optic-f-cp-g}
\]
Intuitively, the `flow of information' in wires labeled $\rev{A}$ and $\rev{B}$ is
right-to-left: the output of $\rev{g}$ is connected to the input of $\rev{f}$, which
itself connects to the left boundary.

Compare this to the graphical syntax of~\cite{optic-string-diagrams},
wherein the author gives a diagrammatic language for Tambara modules using
oriented wires.
\emph{Optics} are then diagrams in this language having the following type.
\[ \gen{optics/example-optic} \]
These diagrams then compose analogously as below.
\[ \gen{optics/optic-compose-f-g} \]
Diagrams in the image of $\ToOptic$ are not in general of this type, but it is
straightforward to adapt them by pre- and post-composing with $\smcinterleave$.

This is useful because for a given $d : A \to B$ the diagram $\ToOptic(d)$
is \emph{not} monogamous acyclic
(Definitions \ref{definition:monogamous} and \ref{definition:acyclic}),
since the \emph{outputs} of each $\rev{f}$ connect to the \emph{left} boundary.
However, if for all $f \in \Sigma_1$, every $\fwd{f}$ and $\rev{f}$ is
monogamous acyclic, then the resulting diagram can be converted to a monogamous
acyclic one using the $\smcinterleave$ morphism as below.
\begin{equation}
  \label{equation:adapted-monogamous-acyclic}
  \input{optics/adapted-monogamous-acyclic.tikz}
\end{equation}
This will be important to our case study because it means that we can extract a
morphism which can be interpreted in a \emph{symmetric monoidal} category.
This means it can be regarded as a \emph{function} which simultaneously computes
the forward and reverse maps of reverse derivatives regarded as optics.

\subsection{Case Study: Efficient Reverse Derivatives via Optic Composition}
\label{section:reverse-derivatives}

To conclude the paper, we now discuss how the constructions given so far can be
applied to the setting of machine learning.
In particular, we will show how the encoding of optics using hypergraph
structure allows for computing reverse derivatives of large diagrams in an
efficient way.

We begin by informally recalling reverse derivatives, a key component of the
formulation of gradient-based learning given in \cite{cfgbl}.
A (cartesian) reverse derivative category (RDC), equips morphisms
$f : A \to B$
with a \emph{reverse derivative}
$\R[f] : A \times B' \to A'$
satisfying various axioms.
Intuitively, this maps an input $A$ and a \emph{change in output} $B'$ to a
change in \emph{input} $A'$.

Reverse derivatives can also be thought of as \emph{lenses}: pairs of morphisms
$(f, \R[f])$ which compose as in the diagram below.
\begin{equation*}
  (f, \R[f]) \cp (g, \R[g]) \qquad
    = \qquad \left( f \cp g \quad , \quad \input{optics/r-composition.tikz} \right)
\end{equation*}

This definition leads to two kinds of inefficiency.
First, we are required to represent two distinct diagrams for each morphism.
This means that it is not possible to apply a functor of the kind defined in
\Cref{section:applying-functors}.
Thus, in order to compute the reverse derivative of a map $f$ as a lens, it
seems necessary to first extract a $\Signature$-term, and then incrementally
build a pair of diagrams.
As we have seen in \Cref{example:quadratic-penalty}, extracting such a
$\Signature$-term can cause a quadratic penalty.
The second problem, pointed out in~\cite{bruno}, is that the
definition of lens composition leads to a space-for-time tradeoff which leads to
multiple redundant `copies' of the forward morphism of the lens.

To avoid this redundant computation, one can use the observation of~\cite[Figure
10]{bruno} which shows how lenses can be composed as \emph{optics}.
This relies on lenses being a special case of optic whose base category is
cartesian.
In terms of optics, the forward maps of lenses are of the form
$\sgen{s-copy-f-id}$, where $\gen{g-copy}$ is the diagonal map.
We may therefore map a morphism $f$ in an RDC to the following optic.
\[
  \ToOptic(d) \qquad = \qquad \input{optics/interleaved-reverse-derivative.tikz}
\]

Having done so, we can then use \eqref{equation:adapted-monogamous-acyclic} to
extract a morphism of the RDC which runs the forward and reverse passes simultaneously.
This gives two specific benefits:
(1) we can efficiently compute a large symbolic description of a reverse derivative in parallel,
and (2) decomposition to $\Signature$-terms is completely avoided.

Note that what has been presented here is more general than just reverse
derivatives.
Since we have already given an algorithm for applying \emph{hypergraph}
functors, we can apply it to those traced symmetric monoidal categories whose
trace is given by the hypergraph structure.
We conjecture that this will allow for modelling systems of bidirectional
information flow \emph{with feedback}.
More precisely, the algorithm given here will allow us to map a morphism with
feedback into an \emph{optic with feedback}.

%% file: conclusion.tex
The datastructures and algorithms presented here improve on previous
work~\cite{coc} in a number of ways.
We handle the case of diagrams equipped with Special Frobenius structure, extend
to arbitrary sets of generating objects, and eliminate the dependency on an
underlying sparse array implementation.
In addition, `natively' allowing for diagrams equipped with Special Frobenius
structure, allows for the representation of diagrams of optics.
This in turn enables us to give a parallel algorithm for taking reverse
derivatives.

However, as in~\cite{coc}, a number of improvements remain.
Aside from algorithms for matching and rewriting, future work should also
include algorithms for \emph{evaluating} diagrams.
While we expect this can be added with little effort, it may be the case that
alternative representations of internal diagram wiring affect the efficiency of
evaluation algorithms.
For instance, we currently store the operations of a diargam in a single array
$G(\sxn) : G(X) \to \Sigma_1$.
However, it may be advantageous for performance reasons to store \emph{multiple}
such arrays so that applying operations in parallel on GPU hardware is more
performant.

It may also be possible to weaken the assumption of the PRAM CRCW model assumed
in the paper.
Since in most cases the CREW model is sufficient, this should require only a few
modifications.
For example, the algorithm to compute the universal morphism of coequalizers
requires the CRCW assumption, but it may be possible to obtain a PRAM CREW
algorithm by replacing concurrent writes with an integer sort.
In addition, we may also be able to improve the work-efficiency of some of the
parallel algorithms presented here, in particular the ancestor maps
(\Cref{proposition:lr-p-ancestor-complexity}).

Finally, when defining algorithms to apply a functor $F$ to a diagram in
\Cref{section:applying-functors}, we required that the typing relation $\tau$ be
a function.
However, one of the advantages of our representation is in allowing for
\emph{polymorphic} generators.
Thus, a natural extension is to allow functors between categories with
polymorphic generators.
To do so would require using \emph{type} information when mapping a given
operation, since there is no longer a unique diagram $F(g)$ which can be
determined simply from an operation label $g \in \Sigma_1$.

%% file: sorting.tex
\begin{proposition}
  $p : A \to A$ is an isomorphism in $\FinFun$
  if and only if $p$ is a permutation.
\end{proposition}
\begin{proof}
  If $p$ is a permutation it is certainly an isomorphism since it has an inverse.
  It therefore suffices to show that if $p$ is an isomorphism, then it is a permutation.
  Since $p$ is an isomorphism, it must be a bijection, and thus a permutation.
\end{proof}

\begin{definition}
  Let $f : A \to B$ be a function.
  An isomorphism $s : A \to A$
  is said to \deftext{sort $A$ by $f$}
  if $s \cp f$ is non-decreasing, i.e.
  for all $i, j \in A$,
  $i < j \implies f(s(i)) \leq f(s(j))$
\end{definition}

\begin{proposition}
  \label{proposition:mono-equivalence}
  Suppose $f : A \to B$ is mono, and $s : A \to A$ sorts by $f$.
  Then $s \cp f$ is strictly increasing, i.e.,
  $i < j \Longleftrightarrow f(s(i)) < f(s(j))$.
\end{proposition}
\begin{proof}
  We have $i < j \implies f(s(i)) \leq f(s(j))$ by definition.
  We prove two cases.
  \begin{description}
    \item[$i < j \implies f(s(i)) < f(s(j))$ (Case 1)]
        In the first case,
        we have $i < j \implies f(s(i)) \leq f(s(j))$
        because $s$ sorts by $f$.
        Then suppose that $f(s(i)) = f(s(j))$.
        Both $f$ and $s$ are monomorphisms, so their composition is as well.
        We can therefore use that $i \neq j \implies f(s(i) \neq f(s(j))$
        to conclude that $i < j \implies f(s(i)) < f(s(j))$.

    \item[$f(s(i)) < f(s(j)) \implies i < j$ (Case 2)]
        In the second case, we prove the contrapositive:
        $ i \geq j \implies f(s(i)) \geq f(s(j)) $
        which we can rewrite as
        $ j \leq i \implies f(s(j)) \leq f(s(i))$.
        This holds because when $j = i$, we have $f(s(j)) = f(s(j))$,
        and when $j < i$, the hypothesis applies.
  \end{description}
\end{proof}

\begin{proposition}
  \label{proposition:unique-sort}
  Let $N$ be a finite set, and $f : N \to B$ a monomorphism.
  There is a unique isomorphism $\sort_f : N \to N$
  that sorts by $f$.
\end{proposition}
\begin{proof}
  Let $r, s : N \to N$ be two permutations sorting by $f$.
  Since $f$ is a monomorphism,
  for all $i, j \in \ord{N}$
  we have
  $f(r(i)) < f(r(j)) \Longleftrightarrow i < j \Longleftrightarrow f(s(i)) < f(s(j))$
  by \Cref{proposition:mono-equivalence}.
  We may then calculate that
  \begin{align*}
    r^{-1}(i) < r^{-1}(j)
      & \Longleftrightarrow f(r(r^{-1}(i))) < f(r(r^{-1}(j)))  \\
      & \Longleftrightarrow f(i) < f(j)  \\
      & \Longleftrightarrow f(s(s^{-1}(i))) < f(s(s^{-1}(j)))  \\
      & \Longleftrightarrow s^{-1}(i) < s^{-1}(j)
  \end{align*}
  We can now use this to show that $r^{-1} = s^{-1}$.

  Recall that $\ord{i}$ denotes the set of natural numbers less than $i$,
  so that $|\ord{i}| = i$ by definition,
  and consequently $|\ord{i}| = |\ord{j}| \Longleftrightarrow i = j$.
  Moreover, for any permutation $p : N \to N$ and $j \in \ord{N}$,
  we have
  $|\ord{j}| = |\{ i \in \ord{N} \mid i < j \}| = |\{ i' \in \ord{N} \mid p(i') < j \}|$
  because for each $i \in \ord{N}$ there is exactly one $i'$ such that $p(i') =
  i$ since $p$ is a bijection. %
  Using this and the fact that
  $r^{-1}(i) < r^{-1}(j) \Longleftrightarrow s^{-1}(i) < s^{-1}(j)$
  we can calculate
  \begin{align*}
    r^{-1}(i) & = \left|\ord{r^{-1}(i)}\right| \\
              & = \left| \left\{ i \in \ord{N} \mid i < r^{-1}(j)         \right\} \right| \\
              & = \left| \left\{ i \in \ord{N} \mid r^{-1}(i) < r^{-1}(j) \right\} \right| \\
              & = \left| \left\{ i \in \ord{N} \mid s^{-1}(i) < s^{-1}(j) \right\} \right| \\
              & = \left| \left\{ i \in \ord{N} \mid i < s^{-1}(j)         \right\} \right| \\
              & = \left|\ord{s^{-1}(i)}\right| \\
              & = s^{-1}(i)
  \end{align*}
  and therefore by uniqueness of inverses, we have that $r = s$.
\end{proof}

\begin{definition}[Stable Sort]
  \label{definition:stable-sort}
  The \deftext{stable sort} by a key $f : A \to B$ is the morphism
  $\sort_k$
  where $k = \tl f, \id \tr : A \to B \times A$
\end{definition}
\begin{remark}
  In the above definition, $s'$ is forced to distinguish elements of $A$ by both
  their `key' $f$ and their `position', given by the identity function.
  In this way, the function $\tl f, \id \tr$ becomes a monomorphism.
  By \Cref{proposition:unique-sort} there is a unique sorting morphism, and so
  the above definition is well-defined.
\end{remark}

%% file: proofs-isomorphism.tex
\subsection{$\FromDiagram$ Proofs}
\label{section:proofs-isomorphism-from-diagram}

In this section, we verify that $\FromDiagram$ (\Cref{proposition:from-diagram})
defines a strict symmetric monoidal hypergraph functor.
We begin with the deferred proof that a Frobenius decomposition exists in any
hypergraph category.

\begin{proof}[Proof of \Cref{proposition:frobenius-decomposition-existence}]
  We show inductively that any morphism in a Hypergraph category has a Frobenius
  decomposition.
  \begin{description}
    \item[Base Case: Generator]
      A Frobenius decomposition of a generator $g$ is its \emph{singleton diagram}:
     \[ g \qquad = \qquad \input{proof/isomorphism-inverse-on-generators-lhs.tikz} \]
    \item[Inductive Step: Tensor $f_0 \tensor f_1$]
      Let $f_0$ and $f_1$ be morphisms having a Frobenius decomposition.
      Then $f_0 \tensor f_1$ has a Frobenius decomposition.
      \begin{align*}
        \gen{proof/frobenius-decomposition-existence-tensor-lhs}
        = \qquad
        & \gen{proof/frobenius-decomposition-existence-tensor-mid} \\
        = \qquad
        & \gen{proof/frobenius-decomposition-existence-tensor-rhs}
      \end{align*}
    \item[Inductive Step: Composition $f_0 \cp f_1$]
      Let $f_0$ and $f_1$ be morphisms having a Frobenius decomposition.
      Then $f_0 \cp f_1$ has a Frobenius decomposition.
      \begin{align*}
        & \gen{spider-lemmas/complete-functor-lhs} \\ \\
        & = \qquad
        \gen{proof/frobenius-decomposition-existence-compose-mid} \\ \\
        & = \qquad
        \gen{proof/frobenius-decomposition-existence-compose-rhs}
        \qedhere
      \end{align*}
  \end{description}
\end{proof}

We now complete the proof sketch of
\Cref{proposition:diagram-frobenius-decomposition},
verifying that the Frobenius decomposition given is indeed isomorphic to the
original diagram.
To do so relies on the following lemma which shows that a diagram can be
immediately factorised into `almost-Frobenius' form.

\begin{proposition}[Almost-Frobenius Decomposition]
  \label{proposition:almost-frobenius-decomposition}
  Let $(s, t, G)$ be a diagram.
  Then there is a choice of coequalizer so that the following `almost Frobenius'
  decomposition is isomorphic to $(s, t, G)$.
  \[ \input{almost-frobenius-decomposition.tikz} \]
  where
  \[
    \hat{s} = \HalfSpider(s)
    \qquad
    \hat{t} = \HalfSpider(t)
    \qquad
    e_s = \HalfSpider(G(\swi))
    \qquad
    e_t = \HalfSpider(G(\swo))
    \qquad
    h = (\inj_0, \inj_1, H)
  \]
  are the half-spiders with the labeling $G(\swn)$
  and $H$ is defined in terms of $G$ as follows.
  \[ \input{almost-frobenius-Hprime.tikz} \]
\end{proposition}
\begin{proof}
  The proof is by direct calculation of the composite with suitable chosen
  coequalizers.
  We will first focus on the inner part of the diagram, ignoring $\hat{s}$ and $\hat{t}$
  until the end of the proof.

  We first compute the diagrams corresponding to following composites.
  \[
    x \quad = \quad \input{proof/almost-frob-x.tikz}
    \qquad \qquad
    y \quad = \quad \input{proof/almost-frob-y.tikz}
    \qquad \qquad
    z \quad = \quad \input{proof/almost-frob-z.tikz}
  \]
  Now, $e_s^\dagger = (\id, G(\swi), \SCRL{G})$ by \Cref{proposition:dagger-swaps-source-target},
  and then \Cref{proposition:frobenius-split-id-f} gives
  $x = \left(\gen{g-identity}, \gen{s-id-plus-gswi}, \SCRL{G}\right)$.
  The same propositions also yield that
  $z = \left(\gen{s-id-plus-gswo}, \gen{g-identity}, \SCRL{G}\right)$.
  Finally, it is straightforward by \Cref{proposition:tensor}
  that $y = ({\gen{g-identity} \atop \gen{g-inj0}}, {\gen{g-identity} \atop \gen{g-inj1}}, \SCRL{G} \tensor H)$.

  We now compute the composite $(x \cp y) \cp z$, in the order indicated by
  bracketing.
  The composite
  \[ (x \cp y)
    =   \left(\gen{g-identity}, \gen{s-id-plus-t}, \SCRL{G}\right)
    \cp \left({\gen{g-identity} \atop \gen{g-inj0}}, {\gen{g-identity} \atop \gen{g-inj1}}, \SCRL{G} \tensor H\right)
  \]
  can be computed by choosing the following coequalizer.
  \[
    \coeq\left(
      \input{proof/almost-frob-xy-coeq-x.tikz}
      \quad
      ,
      \quad
      \input{proof/almost-frob-xy-coeq-y.tikz}
    \right)
    \qquad
    =
    \qquad
    \input{proof/almost-frob-xy-coeq.tikz}
  \]
  This yields a cospan with the source map $\inj_0$ and target map
  $\id \tensor \id$.

  Next, we compute $(x \cp y) \cp z$
  using the following choice of coequalizer.
  \[
    \coeq\left(
      \input{proof/almost-frob-xyz-coeq-xy.tikz}
      \quad
      ,
      \quad
      \input{proof/almost-frob-xyz-coeq-z.tikz}
    \right)
    \qquad
    =
    \qquad
    \input{proof/almost-frob-xyz-coeq.tikz}
  \]
  This yields source and target maps $s = \id$ and $t = \id$, respectively, with
  the apex of the cospan evaluating to $G$.

  It then remains to show that pre- and post-composing with $\hat{s}$ and
  $\hat{t}^\dagger$ yields the diagram $(s, t, G)$.
  Since the source and target maps of $(x \cp y) \cp z$ are identities, this
  follows by \Cref{proposition:coequalizer-f-f-is-id}, completing the proof.
\end{proof}

\begin{remark}
  Note that the above is \emph{not yet a Frobenius decomposition} because $h$ is not
  a tensoring of generators: its boundaries are not necessarily in the desired
  ordering.
  Obtaining a Frobenius decomposition will simply be a matter of decomposing
  $h'$ into the form $p \cp h \cp q$, where $h$ is a tensoring of generators and
  $p$ and $q$ are permutations.
\end{remark}

From this proposition follows the proof that the specific Frobenius
decomposition given in \Cref{proposition:diagram-frobenius-decomposition} indeed
composes to recover the original morphism.

\begin{proof}[Proof of \Cref{proposition:diagram-frobenius-decomposition}]
  First decompose the diagram $(s, t, G)$
  into the `almost-Frobenius' form
  $(\hat{s}, e_s^\dagger, h', e_t, \hat{t}^\dagger)$
  described in \Cref{proposition:almost-frobenius-decomposition},
  To obtain a true Frobenius decomposition, it then suffices to factor the
  diagram $h' = (\inj_0, \inj_1, H')$
  into the composite $p \cp g \cp q$, where $p$ and $q$ are
  permutations, and $g$ is a diagram isomorphic to the $n$-fold tensor product
  of generators.

  To obtain $g$, we will simply permute the wires of $h'$
  to obtain a diagram $(p, q, h)$, where $h$ is a tensoring of generators.
  More concretely, we first construct an isomorphism of structured cospans
  $\alpha : h' \to h$, whose components are defined as follows.
  \[
    \alpha_W = \sort^{-1}_{k_i} \tensor \sort^{-1}_{k_o}
    \qquad
    \alpha_{\Ei} = \sort^{-1}_{k_i}
    \qquad
    \alpha_{\Eo} = \sort^{-1}_{k_o}
    \qquad
    \alpha_X = \id
  \]
  Here, $\sort_f$ denotes the unique permutation sorting by a monomorphism $f$
  (see \Cref{section:sorting})
  while $k_i = \tl G(\sxi), G(\sporti) \tr : X \times \Nat$
  and $k_o = \tl G(\sxo), G(\sporto) \tr : X \times \Nat$
  are the `sorting key' functions.
  Note that $k_i$ and $k_o$ are monomorphisms by well-formedness of $H'$.

  The target of $\alpha$ is the structured cospan
  $h = (\sort_{k_i} \cp \inj_0 \: , \: \sort_{k_o} \cp \inj_1 \: , \: H)$,
  which by
  \Cref{proposition:cospan-with-factorable-isomorphism-legs}
  factors as $p \cp g \cp q$,
  where
  \begin{itemize}
    \item $g = (\inj_0, \inj_1, H)$
    \item $p = (\sort_{k_i}, \id, \MkD{H})$
    \item $q = (\id, \sort_{k_o}, \MkD{H})$
  \end{itemize}

  It therefore remains to show that $g$ is isomorphic to the $n$-fold tensor
  product of generators.
  To do so, we appeal to \Cref{proposition:n-fold-tensor-isomorphism}, from
  which it is immediate that the legs of $g$ are already $\inj_0$ and $\inj_1$
  as required.
  It then suffices to verify that $H$ is of the form given in
  \eqref{equation:n-fold-tensor-isomorphism-h}.

  On objects, $H$ is given by
  \[
    H(W) = G(\Ei) + G(\Eo)
    \qquad
    H(\Ei) = G(\Ei)
    \qquad
    H(\Eo) = G(\Eo)
    \qquad
    H(X) = G(X)
  \]
  and on morphisms, we must have the following by naturality of $\alpha$.
  \[
    H(\swi)    %
                 = \inj_0
   \qquad \qquad
    H(\swo)    %
                 = \inj_1
   \qquad \qquad
    H(\sxi)    %
                 = \sort_{k_i}  \cp G(\sxi)
   \qquad \qquad
    H(\sxo)    %
                 = \sort_{k_o}  \cp G(\sxo)
  \]
  \[
    H(\sporti) %
                 = \sort_{k_i} \cp G(\sporti)
   \qquad
    H(\sporto) %
                 = \sort_{k_o} \cp G(\sporto)
   \qquad
    H(\sxn)    %
                 = G(\sxn)
   \]
   \[
    H(\swn)    %
               = \input{s-frob-decomp-h-swn.tikz}
   \]
  First observe that $H(\swi) = \inj_0$, $H(\swo) = \inj_1$, and $H(\sxn) =
  G(\sxn)$ as in \eqref{equation:n-fold-tensor-isomorphism-h}.
  Write $a : H(\sxn) \to \Obj^*$ and $b : H(\sxn) \to \Obj^*$ for the chosen
  typings of generators in $H$, respectively.
  Now regard the remaining morphisms in the image of $H$ as finite arrays.
  $H(\sxi) = \sort_{k_i} \cp G(\sxi) = \sort_{k_i} \cp \tl G(\sxi), G(\sporti) \tr \cp \proj_0$
  is a non-decreasing array of the form
  $
    0 \overset{\len{a(0)}}{\ldots} 0,
    1 \overset{\len{a(1)}}{\ldots} 1,
    \ldots
  $
  and is therefore equal to
  \[ \bigotimes_{i \in G(X)} \terminal_{\len{a(i)}} \]
  as required.
  A similar equality holds for $H(\sxo)$.
  In a similar way, $H(\sporti)$ is `piecewise increasing', i.e. a concatenation of arrays
  \[
    \tl 0, 1, \ldots, \len{a(0)}-1 \tr
    \: + \:
    \tl 0, 1, \ldots, \len{a(1)}-1 \tr
    \: + \:
    \ldots
    \: + \:
    \tl 0, 1, \ldots, \len{a(G(X)-1)}-1 \tr
  \]
  which is clearly in the form described in
  \eqref{equation:n-fold-tensor-isomorphism-h}.
  Similar reasoning applies to $H(\sporto)$.

  Finally, since $H'$ is well-formed,
  we have
  \[
    \forall e \in \Ei \quad \swn(\swi(e)) = a(\sxi(e))_{\sporti(e)}
    \qquad \qquad
    \forall e \in \Eo \quad \swn(\swo(e)) = b(\sxo(e))_{\sporto(e)}
  \]
  and since $\sort_{k_i}$ is a bijection, in $H$ we have
  \[
    \forall e \in \Ei \quad
      G(\swn)(G(\swi)(\sort_{k_i}(e))) = a(G(\sxi)(\sort_{k_i}(e)))_{G(\sporti)(\sort_{k_i}(e))}
  \]
  and
  \[
    \forall e \in \Eo \quad
      G(\swn)(G(\swo)(\sort_{k_o}(e))) = b(G(\sxo)(\sort_{k_o}(e)))_{G(\sporto)(\sort_{k_o}(e))}
  \]
  so we may think of the function
  $H(\swn) = \input{s-frob-decomp-h-swn.tikz}$
  as the concatenation of finite arrays
  $ a(0) + a(1) + \ldots a(G(X) - 1) + b(0) + b(1) + \ldots + b(G(X)-1) $
  as required by \Cref{proposition:n-fold-tensor-isomorphism}.
\end{proof}

We can now use Frobenius decompositions to define the functor
$\FromDiagram : \DiagramFrob \to \FreeFrob$.

\begin{proposition}
  \label{proposition:from-diagram-is-well-defined}
  $\FromDiagram$ is well-defined.
\end{proposition}
\begin{proof}
  Let $d_i = (s_i, t_i, G_i)$ be diagrams for $i \in \{0, 1\}$,
  and $\alpha : d_0 \to d_1$ an isomorphism of diagrams.
  Assume that the `Almost-Frobenius' decomposition
  (\Cref{proposition:almost-frobenius-decomposition}) of $d_0$ and $d_1$ is a
  bona fide Frobenius decomposition.
  As witnessed by the proof of
  \Cref{proposition:diagram-frobenius-decomposition},
  this assumption is without loss of generality, since any diagram is isomorphic
  to one of the desired form.
  We may therefore regard $\FromDiagram$ as being defined on these
  representatives.
  Then we can calculate as below, where we write $x$ in place of
  $\HalfSpider(x)$ to avoid notational burden.
  However, the reader should be clear that these diargams represent morphisms in
  the category $\FreeFrob$.
  \begin{align*}
    \gen{proof/from-diagram-well-defined-lhs} \qquad
    = \qquad
    & \gen{proof/from-diagram-well-defined-0} \\
    = \qquad
    & \gen{proof/from-diagram-well-defined-1} \\
    = \qquad
    & \gen{proof/from-diagram-well-defined-rhs}
  \end{align*}
  In the first step of the calculation, we use the isomorphism $\alpha$ to
  express components of $d_1$ in terms of $d_0$.
  The second step follows by \Cref{proposition:spider-sliding}.
  The final step is to show that
  $\alpha_{\Ei} \cp g_1 \cp \alpha_{\Ei}^{-1} = g_0$.
  Since $d_0$ and $d_1$ are isomorphic, and $g_0$ and $g_1$ are simply a tensoring of generators,
  there must be a permutation $\pi$ such that
  $g_1 = \pi \cp g_0 \cp \pi^{-1}$
  and $pi$ is `blockwise', i.e. natural in the input types of generators $g_0$.
  It then suffices to show that $\alpha_{\Ei}^{-1} \cp \pi = \id$,
  and therefore $\alpha_{\Ei}^{-1} = \pi$.
  This follows from well-formedness.
  Observe that by naturality $G_1(\sxi) = \alpha_{\Ei}^{-1} \cp G_0(\sxi) \cp \alpha_X$,
  which we may rewrite as $G_1(\sxi) = \alpha_{\Ei}^{-1} \cp \pi \cp G_0(\sxi)$.
  Since $G_1(\sxi)$ and $G_0(\sxi)$ are both monotonic,
  it must be that $\alpha_{\Ei} \cp \pi = (p_0 \tensor p_1 \tensor \ldots \tensor p_X)$
  is a tensoring of permutations.
  That is, each $p_i$ permutes the order of inputs wires for a single generator.
  We now show that each $p_i = \id$.
  By naturality, we can write
  $G_1(\sporti) = \alpha_{\Ei}^{-1} \cp \pi \cp \pi^{-1} \cp G_0(\sporti)$.
  Since we know that $\alpha_{\Ei}^{-1} \cp \pi$ is a tensoring of permutations
  and $\pi$ is blockwise, we must have that
  $G_1(\sporti) = (p_0 \cp \inj_0) + \ldots + (p_X \cp \inj_0)$.
  Finally, we can use that $(f + g) = (\id + h) \implies f = \id$ to conclude
  that every $p_i = \id$, and so $\alpha_{\Ei} = \pi^{-1}$.

  By similar reasoning one can check that $\alpha_{\Eo}$ corresponds to the
  blockwise permutation of generator outputs, and so
  $g_1 = \alpha_{\Ei}^{-1} \cp g_0 \cp \alpha_{\Ei}$,
  completing the proof.
\end{proof}

\begin{proposition}
  \label{proposition:from-diagram-is-functor}
  $\FromDiagram$ is a functor.
\end{proposition}
\begin{proof}
  $\FromDiagram(\id) = \id$ by definition,
  so it remains to show that
  $\FromDiagram(f \cp g) = \FromDiagram(f) \cp \FromDiagram(g)$.

  Let $d_i = (s_i, t_i, G_i)$ be diagrams for $i \in \{0, 1\}$,
  and let $(s_i, e_{s_i}, g_i, e_{t_i}, t_i)$ be the Frobenius decomposition of
  $d_i$ specified in \Cref{proposition:diagram-frobenius-decomposition}.
  For $\FromDiagram$ to be a hypergraph functor,
  Frobenius spiders must be preserved, and moreover we will map the tensoring of
  generators $g_0$ and $g_1$ to a corresponding tensoring in $\FreeFrob$.
  We therefore must show the following holds in $\FreeFrob$
  where $q = \coeq(t_0 \cp \inj_0, s_1 \cp \inj_1)$.
  \begin{align*}
    \FromDiagram(d_0) \cp \FromDiagram(d_1)
    & \cong \gen{spider-lemmas/complete-functor-lhs} \\
    & \cong \gen{spider-lemmas/complete-functor-rhs} \\
    & \cong \FromDiagram(d_0 \cp d_1)
  \end{align*}
  By \cite[Theorem 1.2]{fong-thesis}, we have that
  $ \gen{spider-lemmas/coequalizer-vs-legs-lhs} \cong \gen{spider-lemmas/coequalizer-vs-legs-rhs} $.
  Moreover, since $q$ is a coequalizer, it is a surjective function, meaning
  that its corresponding morphism in $\FreeFrob$ is constructed by tensor and
  composition from the generators
  $\gen{g-identity}$,
  $\gen{g-twist}$,
  and
  $\gen{frob/g-join}$.
  We can therefore prove the following lemma for arbitrary arrows $f_0$ and $f_1$.
  \begin{align*}
    \gen{spider-lemmas/functor-lhs}
      & \derivecong{Axiom} \qquad
      \gen{spider-lemmas/functor-0} \\ \\
      \derivecong{\Cref{proposition:spider-sliding}} \qquad
      & \gen{spider-lemmas/functor-1} \\ \\
      \derivecong{\Cref{proposition:spider-functions}} \qquad
      & \gen{spider-lemmas/functor-2} \\ \\
      \derivecong{\Cref{proposition:spider-sequential-bubbles}} \qquad
      & \gen{spider-lemmas/functor-3} \\ \\
      \derivecong{\Cref{proposition:spider-sliding}} \qquad
      & \gen{spider-lemmas/functor-4} \\ \\
      \derivecong{\Cref{proposition:bubble-through-surjectives}} \qquad
      & \gen{spider-lemmas/functor-rhs}
  \end{align*}
  which completes the proof.
\end{proof}

\begin{proposition}
  \label{proposition:from-diagram-is-strict-monoidal}
  $\FromDiagram$ is a strict monoidal functor
\end{proposition}
\begin{proof}
  We must show that $\FromDiagram(c_0 \tensor c_1) = \FromDiagram(c_0) \tensor \FromDiagram(c_1)$.
  This follows as in the inductive step for the proof of
  \Cref{proposition:frobenius-decomposition-existence}.
\end{proof}

%% file: cospan-proofs.tex
We can now prove the correspondence between diagrams representing Frobenius
Spiders and those whose apexes lie in the image of $\SCL$.
Note that the proof relies on \Cref{proposition:dagger-swaps-source-target},
which is given in \Cref{section:daggers-of-cospans}.

\begin{proof}[Proof of \Cref{proposition:frobenius-spiders-in-diagram}]
  We first show that every diagram of the form
  $d = (s, t, \SCL(Y))$.
  is a Frobenius spider.
  We have
  $
    d
    = (s, \id, \SCL(Y)) \cp (\id, t, \SCL(Y))
    = \HalfSpider(s) \cp \HalfSpider(t)^\dagger
  $
  by \Cref{proposition:dagger-swaps-source-target}.
  Then by \Cref{definition:half-spider},
  $\HalfSpider(s)$ and $\HalfSpider(t)$ are constructed by tensor and composition
  from generators $\gen{g-identity}, \gen{g-twist}, \gen{frob/g-join}$, and
  $\gen{frob/g-unit}$.

  We now verify that every Frobenius spider is a diagram of the form
  $(s, t, \SCL(Y))$.
  The generators of Special Frobenius monoids are the structured cospans as
  given in \cite[Theorem 3.12]{structured-cospans},
  from which it is immediate that each generator has an apex of the form $\SCL(Y)$.
  The apexes of tensor products also have this form, since $\SCL$ preserves
  colimits, so so it remains to check that the apex of a composite
  $(s_0, t_0, \SCL(Y_0)) \cp (s_1, t_1, \SCL(Y_1))$
  is also in the image of $\SCL$.
  Composites are computed by coequalizer.
  Let $q = \coeq(\SCL(t_0) \cp \inj_0, \SCL(s_1) \cp \inj_1)$.
  Then $q = \coeq(\SCL(t_0 \cp \inj_0), \SCL(s_1 \cp \inj_1))$
  because colimits of $\ACSet$s are pointwise.
  For the same reason, we have
  $q_W = \coeq(t_0 \cp \inj_0, s_1 \cp \inj_1)$
  and for components $Z \neq W$
  we have
  $\SCL(x)_Z = \initial : 0 \to 0$
  and thus
  $q_Z = \initial : 0 \to 0$.
  We can then see that
  $q = \SCL(\coeq(t_0 \cp \inj_0, s_1 \cp \inj_1))$,
  completing the proof.
\end{proof}

\subsection{Daggers of Cospans}
\label{section:daggers-of-cospans}
We now prove that the dagger of a diagram swaps the `legs' of the cospan,
so $(s, t, G)^\dagger = (t, s, G)$.
We state this formally in \Cref{proposition:dagger-swaps-source-target}.
To do so requires several lemmas which we first prove.

Note that the proofs in this section are general statements about cospans in
categories having finite colimits and coequalizers (and thus all finite
colimits).

\begin{proposition}
  \label{proposition:coequalizer-f-f-is-id}
  $\coeq(f, f) = \id$
\end{proposition}
\begin{proof}
  It is immediate that $f \cp \id = f \cp \id$,
  and then all $g$ make co-forks $f \cp g$
  with the universal morphism $u = g$.
\end{proof}

\begin{proposition}
  \label{proposition:frobenius-unit-zeroes-source}
  Let $f = \cospan{A}{s}{G}{t}{B}$ be a cospan.
  Then there is a choice of coequalizer so that
  \[
    \input{frob/s-unit-f.tikz}
    \qquad
     =
    \qquad
    \cospan{A}{\gen{g-counit}}{G}{\gen{g-t}}{B}
  \]
\end{proposition}
\begin{proof}
  First, observe that $\sgen{s-id-plus-s}$ is the coequalizer
  $\coeq\left(\sgen{g-inj0}, \sgen{s-s-inj1}\right)$.
  It is immediate by the counit axiom that this makes the co-fork commute,
  then given some $q$ which also makes the co-fork commute, one can check that
  $u = \sgen{s-inj1-cp-q}$ is the universal morphism.
\end{proof}

\begin{proposition}
  \label{proposition:frobenius-counit-zeroes-target}
  Let $f = \cospan{A}{s}{G}{t}{B}$ be a cospan.
  Then there is a choice of coequalizer so that
  \[
    \input{frob/s-f-unit.tikz}
    \qquad
     =
    \qquad
    \cospan{A}{\gen{g-s}}{G}{\gen{g-counit}}{B}
  \]
\end{proposition}
\begin{proof}
  As in \Cref{proposition:frobenius-unit-zeroes-source}.
\end{proof}

\begin{proposition}
  \label{proposition:frobenius-split-id-f}
  Let $f = \cospan{A}{s}{G}{t}{B}$ be a cospan.
  Then there is a choice of coequalizer so that
  \[
    \input{frob/s-split-id-f.tikz}
    \qquad
     =
    \qquad
    \cospan{A}{\gen{g-s}}{G}{\gen{frob/s-s-plus-t}}{A + B}
  \]
\end{proposition}
\begin{proof}
  Factor the cospan into the composite
  $
    (\sgen{g-identity}, \sgen{cmon/g-add}, A)
    \cp
    (\sgen{frob/s-s-tensor-id}, \sgen{frob/s-t-tensor-id}, G)
  $.
  Now choose
  $\coeq\left(\sgen{frob/s-tensor-plus-zero-zero}, \sgen{frob/s-tensor-zero-s-id}\right)
    = \sgen{frob/s-coeq-split-id-f}
  $,
  which commutes by the unit axiom.
  Then given some $q$ which also makes the co-fork commute,
  one can check that there is a unique
  $u = \sgen{frob/s-universal-split-id-f}$.
  The result then follows again by the unit axiom.
\end{proof}

Finally, we can prove the main proposition.

\begin{proposition}
  \label{proposition:dagger-swaps-source-target}
  Let $f = \cospan{A}{s}{G}{t}{B}$ be a cospan.
  Then there is a choice of coequalizer so that
  \[
    \input{frob/s-dagger-f.tikz}
    \qquad
    =
    \qquad
    \cospan{B}{t}{G}{s}{A}
  \]
\end{proposition}
\begin{proof}
  By Propositions
  \ref{proposition:frobenius-unit-zeroes-source}
  and
  \ref{proposition:frobenius-split-id-f}
  we have that
  \[
    \input{frob/s-cup-id-f.tikz} \qquad = \qquad \cospan{\counit}{\gen{g-counit}}{G}{\gen{frob/s-s-plus-t}}{A + B}
  \]
  and so by definition of tensor product we obtain the following.
  \begin{equation*}
    \label{equation:dagger-swaps-source-target-lhs}
    \input{frob/s-cup-id-f-tensor-id.tikz} \qquad = \qquad \cospan{B}{\gen{g-inj1}}{G + B}{\gen{frob/s-s-plus-t-tensor-id}}{A + B + B}
  \end{equation*}
  Similarly, we have
  \begin{equation*}
    \label{equation:dagger-swaps-source-target-rhs}
    \input{frob/s-id-tensor-cap.tikz} \qquad = \qquad \cospan{A + B + B}{\gen{s-id-tensor-plus}}{A + B}{\gen{g-inj0}}{A}
  \end{equation*}
  so it remains to find a coequalizer for the following composite.
  \[ \input{frob/s-cup-id-f-tensor-id.tikz} \quad \cp \quad \input{frob/s-id-tensor-cap.tikz} \]
  Such a coequalizer is given by
  \[
    \input{frob/s-coeq-full.tikz}
  \]
  and the result then follows by the unit axiom.
\end{proof}

\subsection{Frobenius Spiders and Permutations}

Diagrams whose legs are composites with permutations can be factored
into a composite with a spider.
We state this formally with the following pair of propositions.

\begin{proposition}
  \label{proposition:composite-of-permutation}
  Let $c_0 \defeq (s_0, t_0, \SCL(B)) : A \to B$
  and $c_1 \defeq (s_1, t_1, G) : B \to C$
  be diagrams
  where $t_0$ is an isomorphism.
  Then there is a choice of coequalizer such that
  \[ c_0 \cp c_1 = (s_0 \cp t_0^{-1} \cp s_1 \: , \: t_1 \: , \: G) \]
  Similarly, if
  $c_2 \defeq (s_2, t_2, \SCL(C)) : C \to D$
  and $s_2$ is an isomorphism, then
  \[ c_1 \cp c_2 = (s_1 \: , \: t_2 \cp s_2^{-1} \cp t_1 \: , \: G) \]
\end{proposition}
\begin{proof}
  We show the first case, when $t_0$ is an isomorphism, and rely on an
  symmetric argument for the second case, which is omitted.
  We begin by showing that the relevant coequalizer has the following closed form.
  \[ \coeq(t_0 \cp \inj_0, s_1 \cp \inj_1) = \input{proof/composite-of-permutation-coequalizer.tikz} \]
  Observe that the co-fork commutes:
  \[
    \input{proof/composite-of-permutation-commutes-0.tikz}
      \quad = \quad s_1
      \quad = \quad \input{proof/composite-of-permutation-commutes-1.tikz}
  \]
  Now suppose there is some $q'$ such that
  \[
    \input{proof/composite-of-permutation-suppose-qprime-0.tikz} \quad = \quad \input{proof/composite-of-permutation-suppose-qprime-1.tikz}
  \]
  There exists a $u$ such that
  \[
    \input{proof/composite-of-permutation-exists-u-lhs.tikz} \quad = \quad \input{proof/composite-of-permutation-exists-u-rhs.tikz}
  \]
  which we may calculate as follows:
  \begin{align*}
    \input{proof/composite-of-permutation-exists-u-lhs.tikz}
      & = \input{proof/composite-of-permutation-exists-u-0.tikz} \\ \\
      & = \input{proof/composite-of-permutation-exists-u-1.tikz} \\ \\
      & = \input{proof/composite-of-permutation-exists-u-2.tikz} \\ \\
      & = \input{proof/composite-of-permutation-exists-u-3.tikz} \\ \\
      & = \input{proof/composite-of-permutation-exists-u-rhs.tikz}
  \end{align*}
  So we have
  \[
    u \defeq \input{proof/composite-of-permutation-u.tikz}
  \]
  and moreover this $u$ is unique.
  Suppose there is some $v$ with
  \[
    \input{proof/composite-of-permutation-exists-u-lhs.tikz} \quad = \quad \input{proof/composite-of-permutation-v-rhs.tikz}
  \]
  Then we can derive
  \[
    u
    \quad = \quad \input{proof/composite-of-permutation-u.tikz}
    \quad = \quad \input{proof/composite-of-permutation-v-capped.tikz}
    \quad = \quad v
  \]

  Finally, we must also show that the apex of the composite is equal to $G$.
  By \Cref{proposition:composition}, the coequalizer
  $q : \SCL(B) + G \to G'$
  is the $\ACSet$ morphism whose components are $\id$ except for
  $q_W = \coeq(t_0 \cp \inj_0, s_1 \cp \inj_1)$.
  On objects, we have $G'(Y) = G(Y)$ which can be checked in two cases.
  First,
  $(\SCL(B) + G)(Y) = G(Y) = G'(Y)$
  for objects $Y \neq W$, in which case $q_Y$ is the identity.
  Second, observe that the codomain of $q_W$ is $G(W)$,
  and so $G'(W) = G(W)$.
  On morphisms, it suffices to check only those morphisms having source or target
  $G'(W)$, since other all morphisms must be immediately equal by naturality.
  We begin with $G'(\swi)$.
  By naturality, we must have
  $G'(\swi) = (\SCL(B)(\swi) \tensor G(\swi)) \cp q$,
  and since $\SCL(B) = 0$, we also have
  $\SCL(B)(\swi) = \SCL(B)(\swo) = \gen{g-counit}$.
  We can then calculate
  \begin{align*}
    G'(\swi)
      & \quad = \quad \input{proof/composite-of-permutation-swi-0.tikz} \\
      & \quad = \quad \input{proof/composite-of-permutation-swi-1.tikz} \\
      & \quad = \quad G(\swi)
  \end{align*}
  Similarly, we have $G'(\swo) = G(\swo)$.
  Lastly, we use that $G'(\swn)$ is the canonical morphism of the coequalizer to
  calculate that
  \[ G'(\swn)
    \quad = \quad \input{proof/composite-of-permutation-swn-0.tikz}
    \quad = \quad G(\swn) \]
  and so $G' = G$.
\end{proof}

\begin{corollary}
  \label{proposition:cospan-with-factorable-isomorphism-legs}
  Let $c \defeq \cospan{\SCL(A)}{\SCL(p \cp s_0) \cp \iW_G}{G}{\SCL(q \cp t_0) \cp \iW_G}{\SCL(B)}$
  be a diagram.
  Then $c = \hat{p} \cp c \cp \hat{q}$
  where
  $\hat{p} = (p, \id, \SCL(A))$,
  $\hat{q} = (\id, q, \SCL(B))$,
  and $c = (s_0, t_0, G)$.
\end{corollary}

%% file: spider-lemmas.tex
The following lemmas are generic to any hypergraph category.
We use them in particular to prove functoriality of $\FromDiagram$.

\begin{proposition}[Bialgebra]
  \label{proposition:spider-bialgebra}
  \[
    \input{spider-lemmas/bialgebra-lhs.tikz}
    \qquad = \qquad
    \input{spider-lemmas/bialgebra-rhs.tikz}
  \]
\end{proposition}
\begin{proof}
  Immediate by \cite[Theorem 1.2]{fong-thesis}.
  One can also derive this equality directly using the associativity,
  commutativity, and snake axioms.
\end{proof}

\begin{proposition}[Partial Naturality]
  \label{proposition:spider-sliding}
  Let $f : A \to B$ be a morphism constructed by tensor and composition from the
  generators
  $\gen{g-identity}$,
  $\gen{g-twist}$,
  $\gen{frob/g-split}$,
  and
  $\gen{frob/g-join}$.
  Then
  \begin{equation}
    \label{equation:spider-sliding}
    \input{spider-lemmas/split-natural-lhs.tikz} \quad = \quad \input{spider-lemmas/split-natural-rhs.tikz}
    \qquad \qquad
    \input{spider-lemmas/join-natural-lhs.tikz} \quad = \quad \input{spider-lemmas/join-natural-rhs.tikz}
  \end{equation}
\end{proposition}
\begin{proof}
  We prove the first equation by induction, with the second case following in a
  symmetric way.
  Let $f$ be a morphism constructed by tensor and composition from the generators given.
  We first check the base case of the induction.
  \begin{description}
    \item[Base Case $f = \id$]
      Immediate:
      $\gen{frob/g-split} = \gen{frob/g-split}$
    \item[Base Case $f = \gen{g-twist}$]
      \[
        \gen{spider-lemmas/sliding-twist-0}
        \quad = \quad
        \gen{spider-lemmas/sliding-twist-1}
        \quad = \quad
        \gen{spider-lemmas/sliding-twist-2}
      \]
    \item[Base Case $f = \gen{frob/g-split}$]
      \[
        \gen{spider-lemmas/sliding-split-0}
        \quad = \quad
        \gen{spider-lemmas/sliding-split-1}
        \quad = \quad
        \gen{spider-lemmas/sliding-split-2}
        \quad = \quad
        \gen{spider-lemmas/sliding-split-3}
      \]
    \item[Base Case $f = \gen{frob/g-join}$]
      Immediate by \Cref{proposition:spider-bialgebra},
      \[
        \input{spider-lemmas/bialgebra-lhs.tikz}
        \qquad = \qquad
        \input{spider-lemmas/bialgebra-rhs.tikz}
      \]
    \item[Inductive Step]
      Now assume $f_0$ and $f_1$ satisfy \eqref{equation:spider-sliding}.
      Then their tensor product does too:
      \[
        \input{spider-lemmas/sliding-tensor-lhs.tikz}
        \quad = \quad
        \input{spider-lemmas/sliding-tensor-rhs.tikz}
      \]
      And so does their composite:
      \[
        \input{spider-lemmas/sliding-compose-lhs.tikz}
        \quad = \quad
        \input{spider-lemmas/sliding-compose-mid.tikz}
        \quad = \quad
        \input{spider-lemmas/sliding-compose-rhs.tikz}
        \qedhere
      \]
  \end{description}
\end{proof}

\begin{proposition}
  \label{proposition:spider-functions}
  Let $f$ be a morphism formed by tensor and composition from the generators
  $\gen{g-identity}$,
  $\gen{g-twist}$,
  $\gen{frob/g-unit}$,
  and
  $\gen{frob/g-join}$.
  Then we have the following equalities.
  \begin{align}
    \label{equation:spider-functions}
    \gen{spider-lemmas/functions-lhs} \quad = \quad \gen{spider-lemmas/functions-rhs}
    \qquad \qquad & \qquad \qquad
    \gen{spider-lemmas/functions-flip-lhs} \quad = \quad \gen{spider-lemmas/functions-flip-rhs}
    \nonumber
    \\
    \\
    \gen{spider-lemmas/functions-dagger-lhs} \quad = \quad \gen{spider-lemmas/functions-dagger-rhs}
    \qquad & \qquad
    \gen{spider-lemmas/functions-dagger-flip-lhs} \quad = \quad \gen{spider-lemmas/functions-dagger-flip-rhs}
    \nonumber
  \end{align}
\end{proposition}
\begin{proof}
  We first prove the top left equation of \eqref{equation:spider-functions} by
  induction.
  As in \Cref{proposition:spider-sliding}, one can derive base cases directly,
  or observe that they hold via \cite[Theorem 1.2]{fong-thesis}.
  Then the inductive step holds as follows.
  Assume $f_0$ and $f_1$ satisfy the top-left equation of
  \eqref{equation:spider-functions}.
  Then
  \[
    \gen{spider-lemmas/functions-tensor-lhs}
    \qquad = \qquad
    \gen{spider-lemmas/functions-tensor-rhs}
  \]
  and
  \[
    \gen{spider-lemmas/functions-compose-lhs}
    \qquad = \qquad
    \gen{spider-lemmas/functions-compose-mid}
    \qquad = \qquad
    \gen{spider-lemmas/functions-compose-rhs}
  \]
  as required.
  We can now show the remaining equations hold as well.
  The top right equation holds by commutativity
  \[ 
    \gen{spider-lemmas/functions-flip-lhs}
    \quad = \quad
    \gen{spider-lemmas/functions-flip-0}
    \quad = \quad
    \gen{spider-lemmas/functions-flip-1}
    \quad = \quad
    \gen{spider-lemmas/functions-flip-2}
    \quad = \quad
    \gen{spider-lemmas/functions-flip-rhs}
  \]
  the bottom left equation is simply the dagger of the top left:
  \[
    \gen{spider-lemmas/functions-dagger-lhs}
    \quad = \quad
    \left[ \gen{spider-lemmas/functions-lhs} \right]^\dagger
    \quad = \quad
    \left[ \gen{spider-lemmas/functions-rhs} \right]^\dagger
    \quad = \quad
    \gen{spider-lemmas/functions-dagger-rhs}
  \]
  and the bottom right equation holds by commutativity as in the proof of the
  top right equation.
\end{proof}

\begin{proposition}
  \label{proposition:spider-sequential-bubbles}
  \[
    \gen{spider-lemmas/spider-sequential-bubbles-lhs}
    \qquad = \qquad
    \gen{spider-lemmas/spider-sequential-bubbles-rhs}
  \]
\end{proposition}
\begin{proof}
  Follows by \cite[Theorem 1.2]{fong-thesis}.
\end{proof}

\begin{proposition}
  \label{proposition:spider-inverse-surjective}
  Let $q = \HalfSpider(f)$ for some surjective $f$,
  so that $q$ is constructed by tensor and composition from generators
  $\gen{g-identity}$,
  $\gen{g-twist}$,
  and
  $\gen{frob/g-join}$.
  Then $q^\dagger \cp q = \id$.
\end{proposition}
\begin{proof}
  Induction. In the base case, it is easy to see that the equation holds when
  $q = \gen{g-identity}$ or $q = \gen{g-twist}$.
  Further,
  $\input{frob/s-split-join.tikz} = \gen{g-identity}$
  holds as an axiom.
  It therefore remains to prove the inductive step.
  Let $q_0$ and $q_1$ be morphisms such that
  $q_i^\dagger \cp q_i = \id$ for $i \in \{0, 1\}$.
  Then the equation holds for their tensor product
  \[
    \gen{spider-lemmas/inverse-surjective-tensor-lhs}
    \qquad = \qquad
    \gen{spider-lemmas/inverse-surjective-tensor-rhs}
  \]
  and also for their composite.
  \[
    \gen{spider-lemmas/inverse-surjective-compose-lhs}
    \qquad = \qquad
    \gen{spider-lemmas/inverse-surjective-compose-0}
    \qquad = \qquad
    \gen{spider-lemmas/inverse-surjective-compose-1}
    \qquad = \qquad
    \gen{g-identity}
    \qedhere
  \]
\end{proof}

Note that in the theory of \emph{Extra} Special Frobenius Monoids,
\Cref{proposition:spider-inverse-surjective}
holds for \emph{any} morphism in the image of $\HalfSpider$; the surjectivity
requirement is dropped.
This is due to the presence of the additional axiom
\[ \gen{frob/g-unit} \cp \gen{frob/g-counit} = \gen{g-empty} \]
which proves the remaining base case.

\begin{proposition}
  \label{proposition:split-join-swap}
  \begin{equation}
    \label{equation:split-join-swap}
    \gen{spider-lemmas/split-join-swap-lhs}
    \qquad = \qquad
    \gen{spider-lemmas/split-join-swap-rhs}
  \end{equation}
\end{proposition}
\begin{proof}
  \[
    \gen{spider-lemmas/split-join-swap-lhs}
    \quad = \quad \gen{spider-lemmas/split-join-swap-0}
    \quad = \quad \gen{spider-lemmas/split-join-swap-1}
    \quad = \quad \gen{spider-lemmas/split-join-swap-rhs}
    \qedhere
  \]
\end{proof}

The final lemma forms a key step in the proof that $\FromDiagram$ is functorial.
However, its derivation is somewhat fiddly, so we omit it from the main part of
the proof.

\begin{proposition}
  \label{proposition:bubble-through-surjectives}
  Let $f_0$ and $f_1$ be arbitrary morphisms, and
  suppose $q = \HalfSpider(q')$ for some surjective function $q'$.
  That is, assume $q$ is constructed by tensor and composition from generators
  $\gen{g-identity}$,
  $\gen{g-twist}$,
  and
  $\gen{frob/g-join}$.
  Then
  \[
    \gen{spider-lemmas/bubble-through-surjectives-lhs}
    \qquad = \qquad
    \gen{spider-lemmas/bubble-through-surjectives-rhs}
  \]
\end{proposition}
\begin{proof}
  We derive as follows.
  \begin{align*}
    \gen{spider-lemmas/bubble-through-surjectives-lhs} \qquad
    & = \qquad
    \gen{spider-lemmas/bubble-through-surjectives-0} \\
    & = \qquad
    \gen{spider-lemmas/bubble-through-surjectives-1} \\
    & = \qquad
    \gen{spider-lemmas/bubble-through-surjectives-2} \\
    & = \qquad
    \gen{spider-lemmas/bubble-through-surjectives-3} \\
    & = \qquad
    \gen{spider-lemmas/bubble-through-surjectives-4} \\
    & = \qquad
    \gen{spider-lemmas/bubble-through-surjectives-rhs}
  \end{align*}
  In the above,
  the first step uses \Cref{proposition:split-join-swap},
  the penultimate step uses \Cref{proposition:spider-inverse-surjective},
  and all other steps use \Cref{proposition:spider-functions}.
\end{proof}